\documentclass[12pt]{article}
\usepackage{amscd}
\usepackage[all]{xy}
\usepackage{CJK}
\usepackage{CJKnumb}
\usepackage{graphicx}
\usepackage{amsfonts}
\usepackage{amsmath}
\usepackage{amssymb}
\usepackage{fancyhdr}
\usepackage{indentfirst}
\usepackage{titlesec}
\usepackage{mathrsfs}
\usepackage{bbm}
\usepackage{dsfont}
\usepackage{amsthm} 
\usepackage{mathrsfs}

\newcommand{\cM}{\mathcal{M}}

\newcommand{\cR}{\mathcal{R}}
\newcommand{\cE}{\mathcal{E}}

\newcommand{\cN}{\mathcal{N}}
\newcommand{\cA}{\mathcal{A}}
\newcommand{\cD}{\mathcal{D}}
\newcommand{\cC}{\mathcal{C}}

\newcommand{\cB}{\mathcal{B}}
\newcommand{\cK}{\mathcal{K}}

\theoremstyle{definition}  \newtheorem{Def}{Definition}
\theoremstyle{plain}  \newtheorem{Thm}{Theorem} \theoremstyle{plain}
\newtheorem{cor}{Corollary}
 \theoremstyle{remark} \newtheorem{Rek}{Remark}
 \theoremstyle{plain}\newtheorem{Lemma}{Lemma}
 \theoremstyle{plain}
\newtheorem{Prop}{Proposition}
\theoremstyle{remark}

\begin{document}

\title{\textbf{2-Modules and the Representations of 2-Rings}}
\date{}
\author{Fang Huang, Shao-Han Chen, Wei Chen, Zhu-Jun Zheng\thanks{Supported in part by NSFC with grant Number
10971071
 and Provincial Foundation of Innovative
Scholars of Henan.} }
\maketitle

\begin{center}
\begin{minipage}{5in}
{\bf  Abstract}: In this paper, we develop 2-dimensional algebraic
theory which closely follows the classical theory of modules. The
main results are giving definitions of 2-modules and the
representations of 2-rings. Moreover, for a 2-ring $\cR$, we prove
that its modules form a 2-abelian category.

 \\
\end{minipage}
\end{center}
\section{Introduction}
Ordinary algebra (which we call it 1-dimensional algebra) is the
study of algebraic structures on sets, i.e. of sets equipped with
certain operations satisfying certain equations. The work on
2-dimensional algebra is to study the algebraic structures on
groupoids\cite{18}.

One of the goals of higher-dimensional algebra is to categorify
mathematical concepts, a very example is (symmetric) 2-groups \cite
{2,3,5,6,10,11,12,16,18,20}, which play the role similar as
(abelian) groups in 1-dimensional algebra. The higher-dimensional
algebra was studied and used in many fields of mathematics such as
algebraic geometry \cite{2,3}, topological field theory\cite
{13,14}, etc.

In {\cite 8}, M.Jibladze and T.Pirashvili introduced categorical
rings(We call them 2-rings). As 1-dimensioal algebra, it's natural
to define 2-modules and the representations of 2-rings. M. Dupont in
his PhD. thesis \cite {18}, mentioned the 2-modules, as the additive
Gpd-functors, also proved that the 2-category formed by these
2-modules is a 2-abelian $Gpd$-category.

In this paper, we give a definition of the 2-module following the
1-dimensional case. When we finish this paper, we find that
V.Schmitt gave another definition similar as 2-module(\cite{19}).
But our definition is much closer to the classical case. The
equivalence between our definition and M.Dupont's definition will be
given in our coming paper.

Abeian category plays an important role in homology theory,
localization theory, representation theory, etc\cite{1,4,17}. So
2-abelian category should be an important tool in studying higher
dimensional algebraic theory.

In our coming papers, we will prove that for any a 2-abelian
category $\cA$, there exists a 2-ring $\cR$, and the embedding of
$\cA$ into the $Gpd$-category $(\cR$-2-Mod), extending from
Freyd-Mitchell Embedding Theorem\cite{4}, and study the (co)homology
theory in 2-abelian category similar as in 1-dimensional case.

This paper is organized as follows:

After recalling basic facts on (symmetric)2-groups, 2-rings, we give
the definitions of $\cR$-2-modules, $\cR$-homomorphism between them,
and morphism between $\cR$-homomorphisms, where $\cR$ is a 2-ring.
As an application of $\cR$-2-modules, we also give the definition of
representation of 2-rings. For our next work, we concretely
construct the 2-category structure of all $\cR$-2-modules. In
section 3, we prove ($\cR$-2-Mod) is an additive, and also a
2-abelian $Gpd$-category by using the similar methods as M.Dupont
discussing the symmetric 2-groups(\cite{18}). This section is the
main part of this paper.

\section{Basic Results on 2-Modules}
Our goal in this section is to give the definition of
$\cR$-2-modules, and the representations of 2-rings. We will begin
by reviewing some definitions about (symmetric)2-groups and 2-rings,
also called categorical groups and categorical rings in
\cite{5,6,8,11,16}.
\begin{Def}\cite{8}
A 2-group $\cA$ is a groupoid equipped with a monoidal structure,
i.e. a bifunctor $+ : \cA\times\cA\rightarrow \cA$, an unit object
$0\in\cA$, and natural isomorphisms:
\begin{align*}
&<a,b,c>:(a+b)+c\rightarrow a+(b+c),\\
&\hspace{1.5cm}l_{a}: 0+a\rightarrow a,\\
&\hspace{1.5cm}r_{a}: a+0\rightarrow a
\end{align*}
satisfying the Mac Lane coherence conditions, i.e. the following
diagrams commute:
 \begin{center}
 \scalebox{0.9}[0.85]{\includegraphics{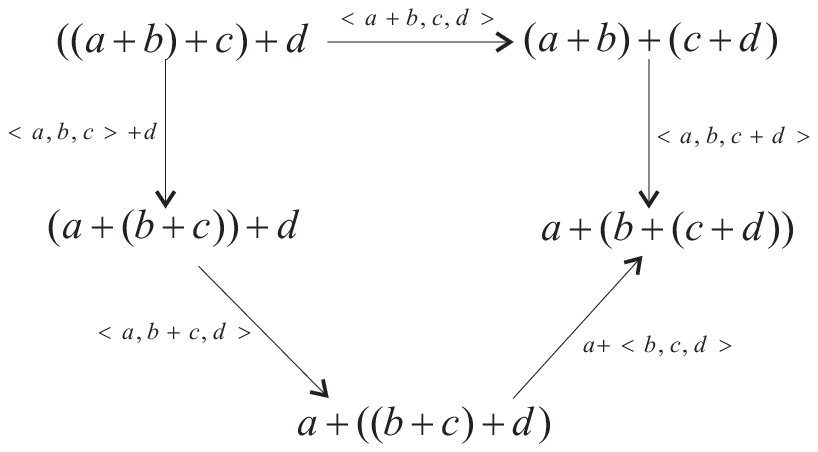}}
{\footnotesize Fig.1.}
\end{center}
 \begin{center}
 \scalebox{0.9}[0.85]{\includegraphics{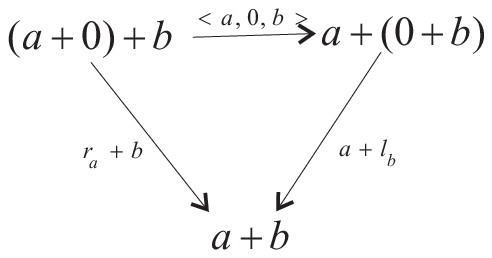}}
{\footnotesize Fig.2.}
\end{center}
\end{Def}
Moreover, for each object $a\in\cA$, there exists an object
$a^{*}\in\cA$, and an isomorphism $\eta_{a}: a^{*}+a\rightarrow 0$.

\begin{Def}\cite {18}
Symmetric 2-group is a 2-group $\cA$, together with the natural
isomorphism $c_{a,b}: a+b\rightarrow b+a$ satisfies $c_{a,b}\circ
c_{b,a}=id$. Also, $c_{a,b}$ is compatible with $<-,-,->$ of the
monoidal structure $+$ .
\end{Def}
\begin{Rek}
In our paper, when we say $\cA$ is a symmetric 2-group, it means
$(\cA,0,+, <-,-,->,l_-,r_-,\eta_-,c_{-,-})$.
\end{Rek}
For any $a,b,c,d\in obj(\cA)$(\cite 8),
$$\langle^{a\ b}_{c\ d}\rangle: (a+b)+(c+d)\longrightarrow (a+c)+(b+d)$$
will be denoted the composite canonical isomorphism in the following
commutative diagram:
\begin{center}
 \scalebox{0.9}[0.85]{\includegraphics{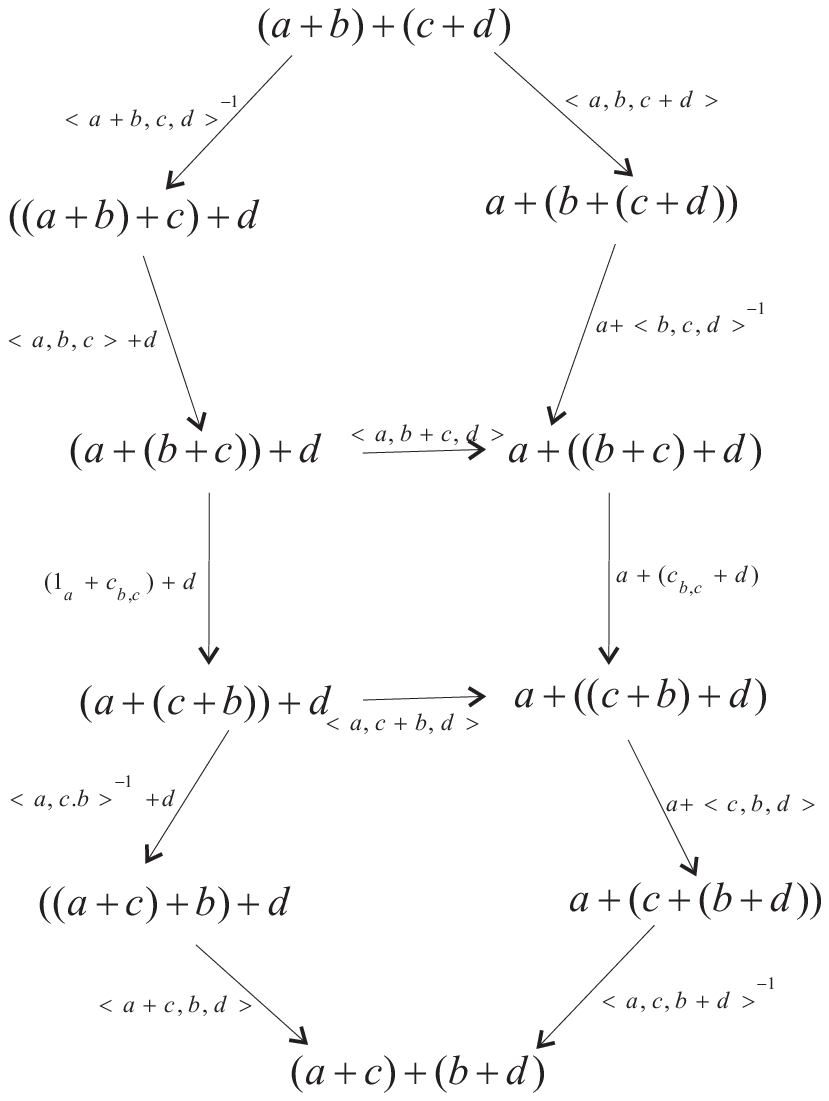}}
 {\footnotesize Fig.3.}
\end{center}

\begin{Def}\cite {18} Let
$\cA, \cB$ be 2-groups. A homomorphism ${\bf F}=(F,F_{+},F_{0}):\cA
\rightarrow \cB$ consists of a functor $F: \cA\rightarrow \cB$ and
two natural morphisms:
\begin{align*}
&F_{+}(a,b): F(a+b)\rightarrow F(a)+F(b),\\
&\hspace{1.3cm}F_{0}:F(0)\rightarrow 0
\end{align*}
such that the following diagrams commute:
\begin{center}
 \scalebox{0.9}[0.85]{\includegraphics{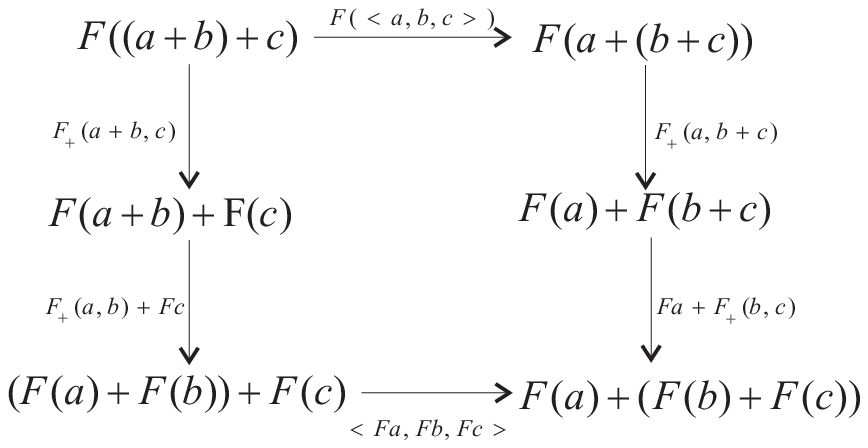}}
{\footnotesize Fig.4.}
\end{center}
\begin{center}
 \scalebox{0.9}[0.85]{\includegraphics{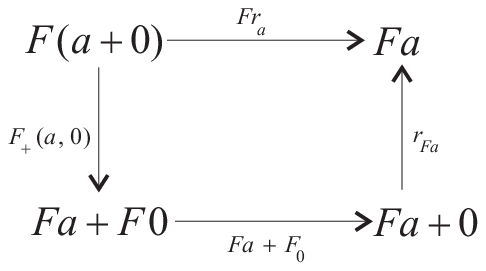}}
{\footnotesize Fig.5.}
\end{center}
\begin{center}
 \scalebox{0.9}[0.85]{\includegraphics{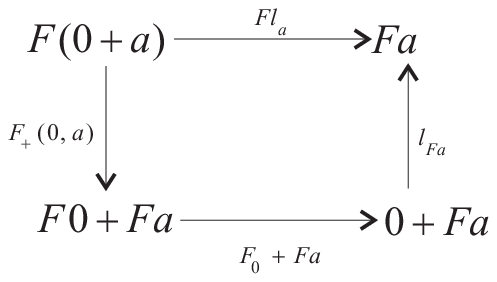}}
{\footnotesize Fig.6.}
\end{center}
\end{Def}

\begin{Rek}
(i)\cite {18} The homomorphism ${\bf F}=(F,F_{+},F_{0})$ is called
strict when isomorphisms $F_{+}, F_{0}$ are identities.

(ii) If $\cA,\ \cB$ are two symmetric 2-groups, the homomorphism of
symmetric 2-groups is the homomorphism ${\bf
F}=(F,F_+,F_0):\cA\rightarrow\cB$, together with the following
commutative diagram:
\begin{center}
 \scalebox{0.9}[0.85]{\includegraphics{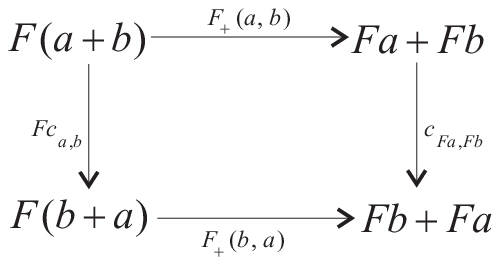}}
{\footnotesize Fig.7.}
\end{center}

(iii) Given homomorphisms $\cA\xrightarrow[]{{\bf
F}=(F,F_+,F_0)}\cB\xrightarrow[]{{\bf G}=(G,G_+,G_0)} \cC$ of
2-groups, their composition is ${\bf H}=(H,H_+,H_0):\cA \rightarrow
\cC$, where $H=G\circ F: \cA\rightarrow \cC$ is a composition of
functors, and $H_+,H_0$ are the following compositions:
\begin{align*}
&H_{+}(a,c): H(a+c)=(GF)(a+c)=G(F(a+c))\\
&\hspace{2cm}\xrightarrow[]{G(F_+(a,c))} G(Fa+Fc)\xrightarrow[]{G_+(Fa,Fc)} G(Fa)+G(Fc)=Ha+Hc,\\
&H_{0}:(GF)(0)=G(F(0))\xrightarrow[]{G(F_{0})}G(0)\xrightarrow[]{G_{0}}0.
\end{align*}
It is easy to check ${\bf H}=(H, H_+, H_0)$ is a homomorphism from
$\cA$ to $\cC$.
\end{Rek}
\noindent\textbf{Notation.} The homomorphism ${\bf F}=(F,F_+,F_0)$
of 2-groups is in fact a functor satisfies some compatible
conditions(Fig.4.-6.). So we will only write $F$ for abbreviation.
\begin{Def}\cite {18}
Given homomorphisms $F,\ G: \cA\rightarrow \cB$ of 2-groups, a
morphism from $F$ to $G$ is a natural transformation $\varepsilon:
F\Rightarrow G$ such that, for any objects $a,\ b\in\cA$, the
following diagrams commute:
\begin{center}
\scalebox{0.9}[0.85]{\includegraphics{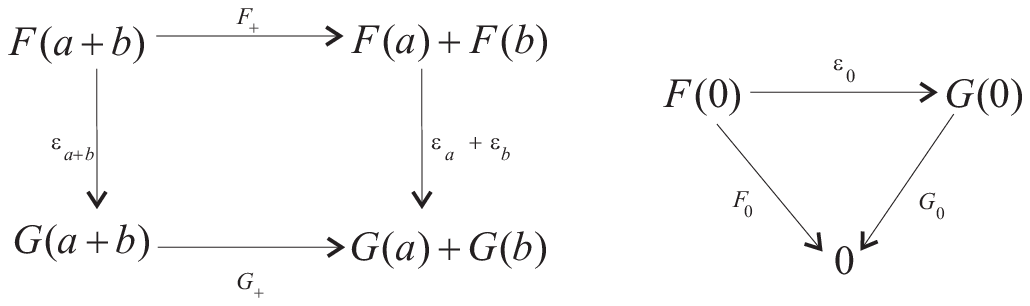}}
\end{center}
\end{Def}

\begin{Rek}
The morphism between two homomorphisms $F,\ G$ of 2-groups $\cA,\cB$
is a natural isomorphism. In fact, for any object $a\in\cA$,
$\varepsilon_{a}: F(a)\rightarrow G(a)$ is a morphism in groupoid
$\cB$, so $\varepsilon_{a}$ is invertible.
\end{Rek}

\begin{Prop}\cite{3}
There is a 2-category (2-Gp) with 2-groups as objects, homomorphisms
of 2-groups as 1-morphisms, and morphisms of homomorphisms as
2-morphisms. If the 2-groups are symmetric 2-groups, we denote this
2-category by (2-SGp).
\end{Prop}
\begin{Def}\cite 8
A 2-ring is a symmetric 2-group $\cR$, together with a bifunctor
$\cdot:\cR\times\cR\rightarrow \cR\ $(denoted by multiplication), an
object $1\in \cR$, and natural isomorphisms:

$[r,s,t]: (r \cdot s)\cdot t\rightarrow r\cdot(s\cdot t)$
(associativity),

$\hspace{0.7cm}\lambda_{r}: 1\cdot r\rightarrow r$ (left unitality),

$\hspace{0.7cm}\rho_{r}: r\cdot 1\rightarrow r$ (right unitality),

$[r_{s_{0}}^{s_{1}}>: r\cdot(s_{0}+s_{1})\rightarrow r\cdot
s_{0}+r\cdot s_{1}$ (left distributivity),

$<_{r_{1}}^{r_{0}}s]: (r_{0}+r_{1})\cdot s\rightarrow r_{0}\cdot
s+r_{1}\cdot s$ (right distributivity).

It is required that the $[-,-,-]$, together with $\lambda_{-}$ and
$\rho_{-}$ constitute a monoidal structure (Fig.1.-2. commute).
Moreover, the following diagrams commute for all possible objects of
$\cR$:
\begin{center}
\scalebox{0.9}[0.85]{\includegraphics{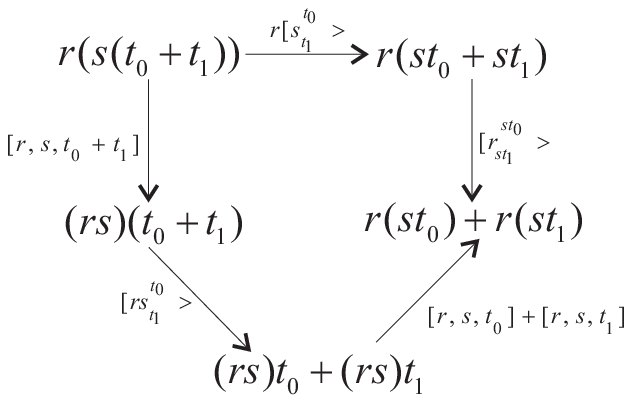}} {\footnotesize
Fig.8.}
\scalebox{0.9}[0.85]{\includegraphics{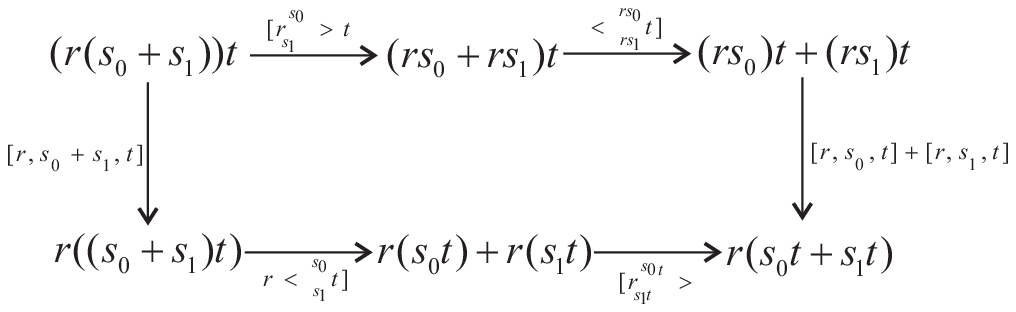}}{\footnotesize
Fig.9.}

\scalebox{0.9}[0.85]{\includegraphics{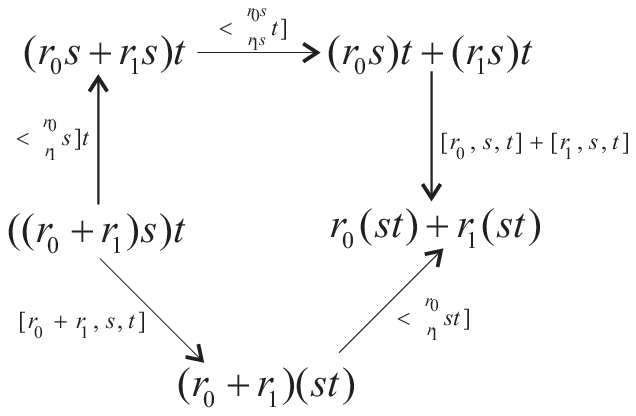}}{\footnotesize
Fig.10.}

\scalebox{0.9}[0.85]{\includegraphics{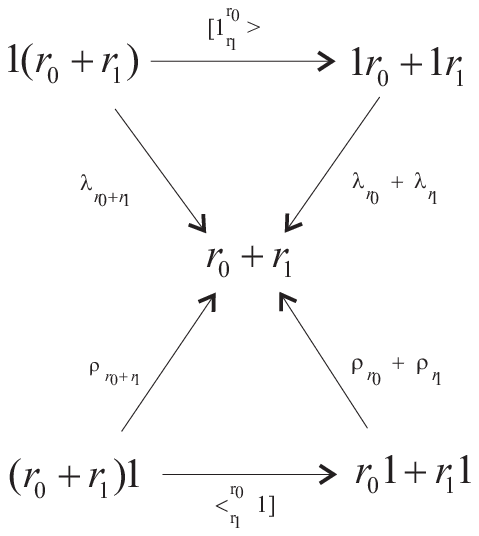}}{\footnotesize
Fig.11.}

\scalebox{0.9}[0.85]{\includegraphics{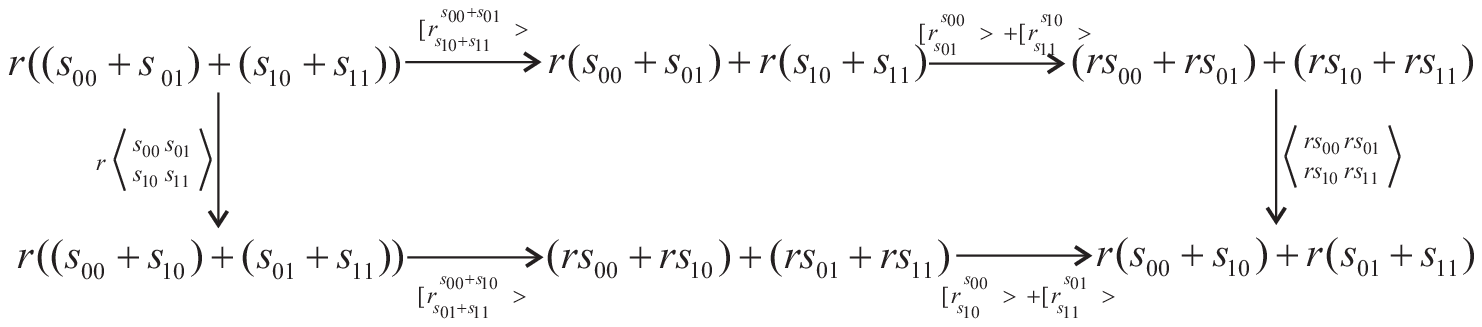}}{\footnotesize
Fig.12.}

\scalebox{0.9}[0.85]{\includegraphics{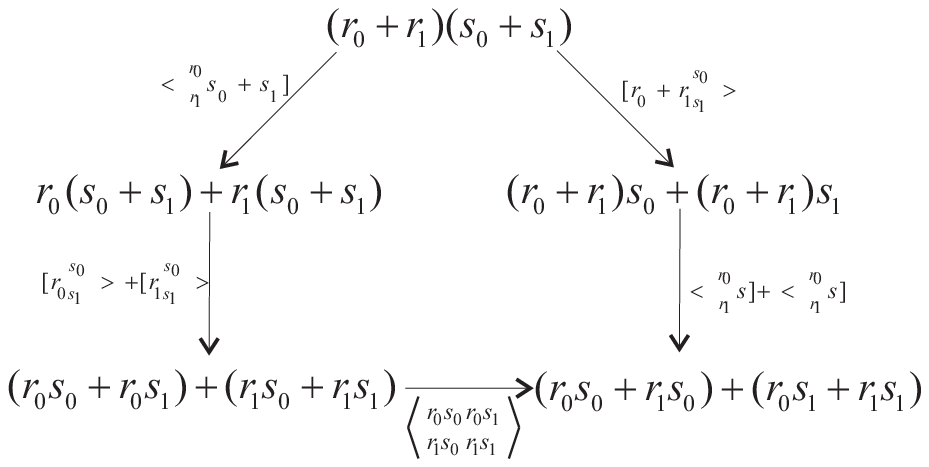}}{\footnotesize
Fig.13.}

\scalebox{0.9}[0.85]{\includegraphics{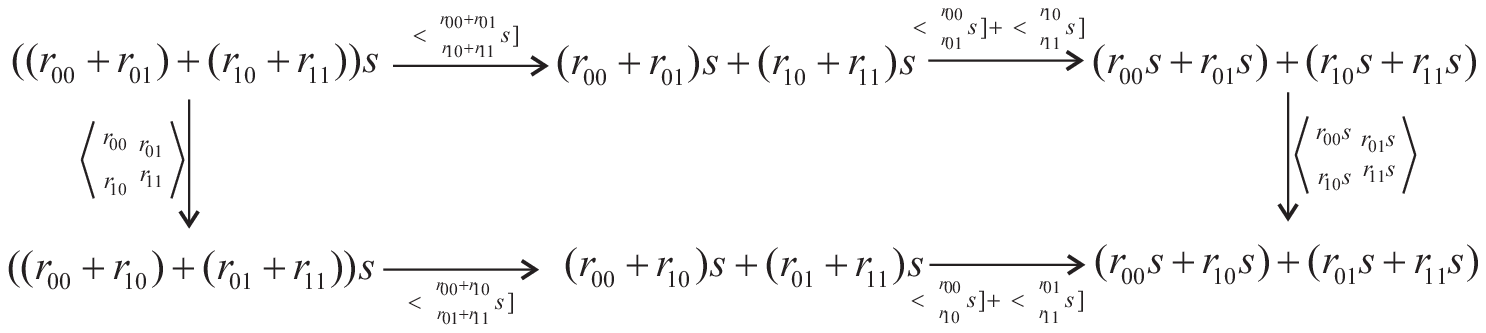}} {\footnotesize
Fig.14.}
\end{center}
\end{Def}

\begin{Rek}
1. N.T.Quang in \cite{14} discussed the relations between these
2-rings and Ann-categories.

2. If the above natural isomorphisms are identities, $\cR$ is called
a strict 2-ring.
\end{Rek}

\begin{Def}\cite 8
A homomorphism of 2-rings is a quadruple
${\bf F}=(F,F_{+},F_{\cdot},F_{1})$, where $F$ is a functor from $\cR_{1}$
to $\cR_{2}$, $F_{+}, F_{\cdot}$ are natural morphisms of the forms
$$F_{+}(r_{0},r_{1}): F(r_{0}+r_{1})\rightarrow F(r_{0})+F(r_{1})$$
$$F_{\cdot}(r_{0},r_{1}): F(r_{0}\cdot r_{1})\rightarrow F(r_{0})\cdot F(r_{1})$$
and $F_{1}: F(1)\rightarrow 1,\ F_{0}:F(0)\rightarrow 0$ are
morphisms, such that $(F,F_{+},F_{0})$, $(F,F_{\cdot},F_{1})$ are
monoidal functors with respect to the monoidal structures
corresponding to the $+$ and $\cdot$, respectively. Moreover the
diagrams commute:
\begin{center}
\scalebox{0.9}[0.85]{\includegraphics{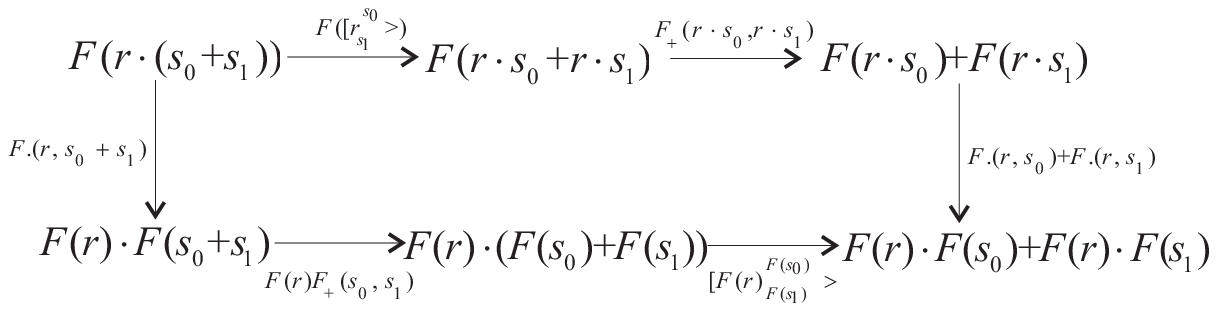}}{\footnotesize
Fig.15.}
\scalebox{0.9}[0.85]{\includegraphics{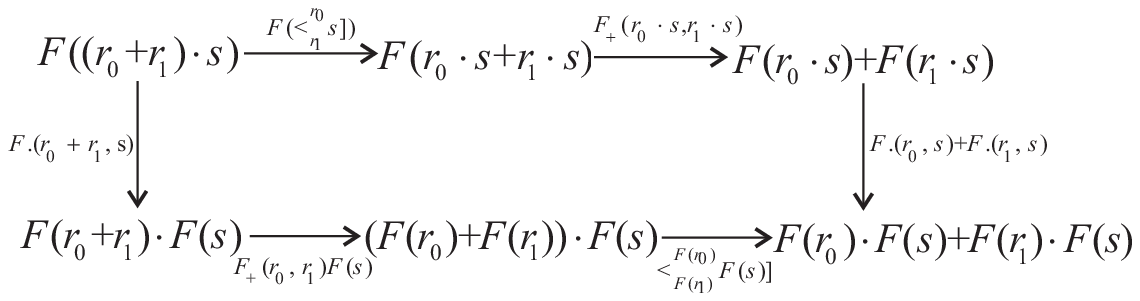}}{\footnotesize
Fig.16.}
\end{center}
\end{Def}
\noindent\textbf{Notation.} The homomorphism ${\bf
F}=(F,F_+,F_{\cdot},F_0)$ of 2-rings is in fact a functor $F$
satisfies some compatible conditions(Fig.4.,5.,15.,16.). So, we will
only write $F$ for abbreviation.
\begin{Def}
The morphism of homomorphisms $F,G: \cR_{1}\rightarrow \cR_{2}$ of
2-rings is a natural transformation $\varepsilon: F\Rightarrow G$,
such that, for any objects $r,\ s\in\cR_{1}$, the following diagrams
commute:
\begin{center}
\scalebox{0.9}[0.85]{\includegraphics{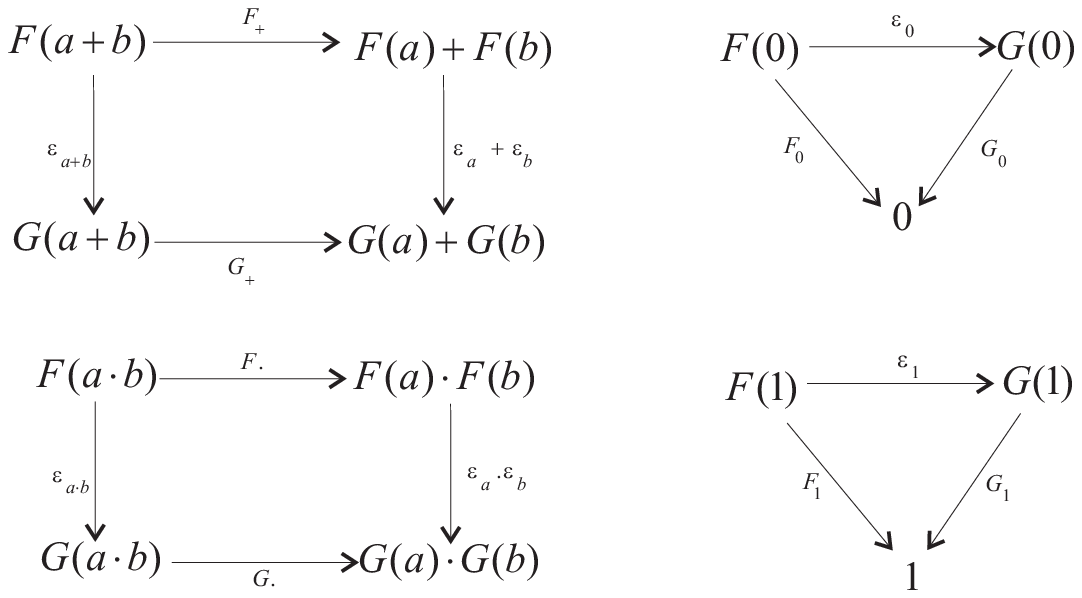}}{\footnotesize
Fig.17.}
\end{center}
\end{Def}

\begin{Prop}\cite {18}There is a 2-category with 2-rings as objects, homomorphisms
of 2-rings as 1-morphisms, and morphisms of homomorphisms as
2-morphisms.
\end{Prop}
\begin{Def}
Let $\cR$ be a 2-ring. An $\cR$-2-module is a symmetric 2-group
$\cM$ equipped with a bifunctor $\cdot: \cR\times\cM\rightarrow \cM$
(called operation of $\cR$ on $\cM$) denoted by $(r,m)\mapsto r\cdot
m$ and natural isomorphisms:
\begin{align*}
&a_{m,n}^{r}: r\cdot(m+n)\rightarrow r\cdot m+r\cdot n,\\
&b_{m}^{r,s}:(r+s)\cdot m\rightarrow r\cdot m+r\cdot n,\\
&b_{r,s,m}:(rs)\cdot m\rightarrow r\cdot(s\cdot m),\\
&i_{m}:I\cdot m\rightarrow m,\\
&z_{r}:r\cdot 0\rightarrow 0
\end{align*}
such that the following diagrams commute:
\begin{center}
\scalebox{0.9}[0.85]{\includegraphics{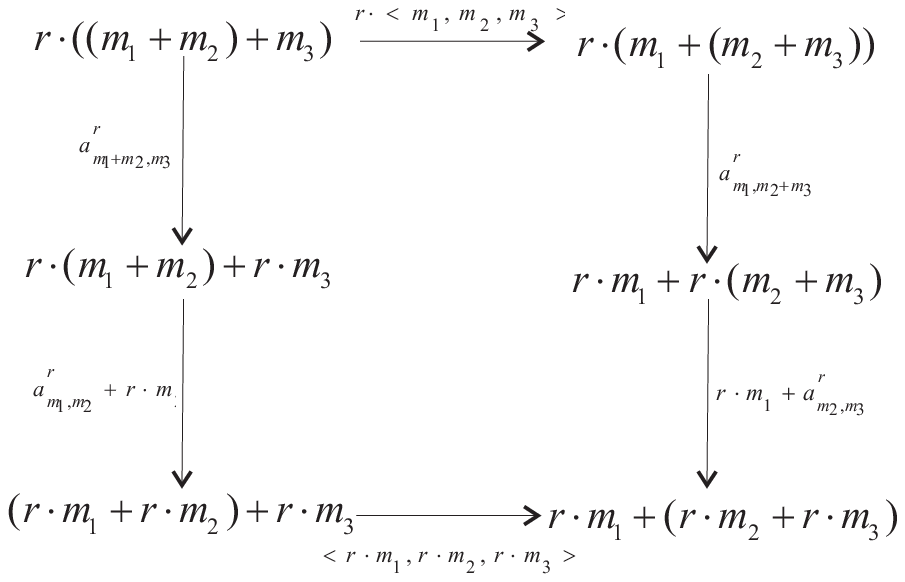}}{\footnotesize
Fig.18.}
\scalebox{0.9}[0.85]{\includegraphics{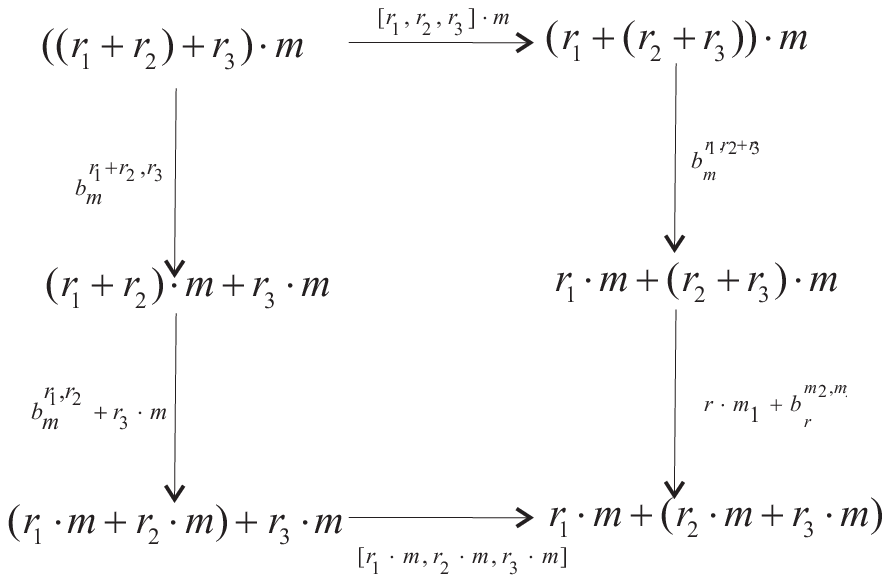}}{\footnotesize
Fig.19.}
\scalebox{0.9}[0.85]{\includegraphics{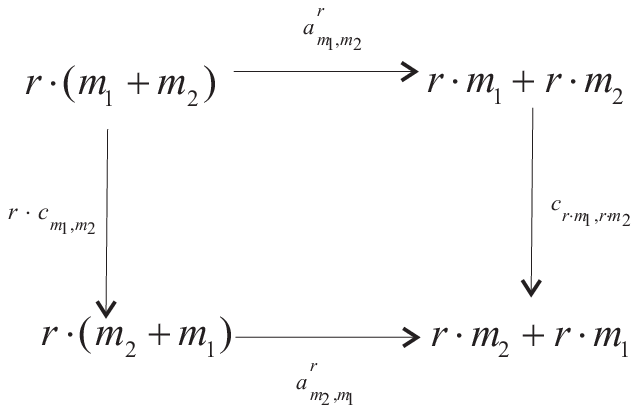}}{\footnotesize
Fig.20.}
\scalebox{0.9}[0.85]{\includegraphics{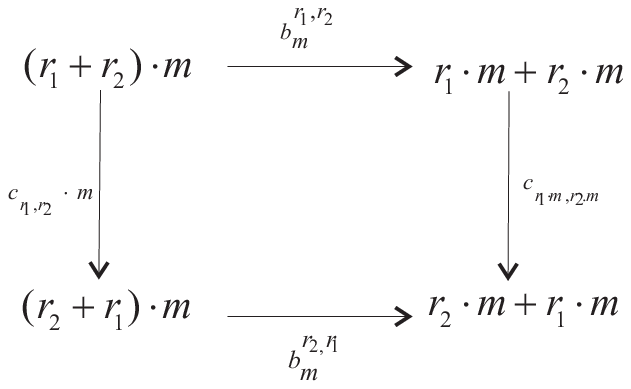}}{\footnotesize
Fig.21.}
\scalebox{0.9}[0.85]{\includegraphics{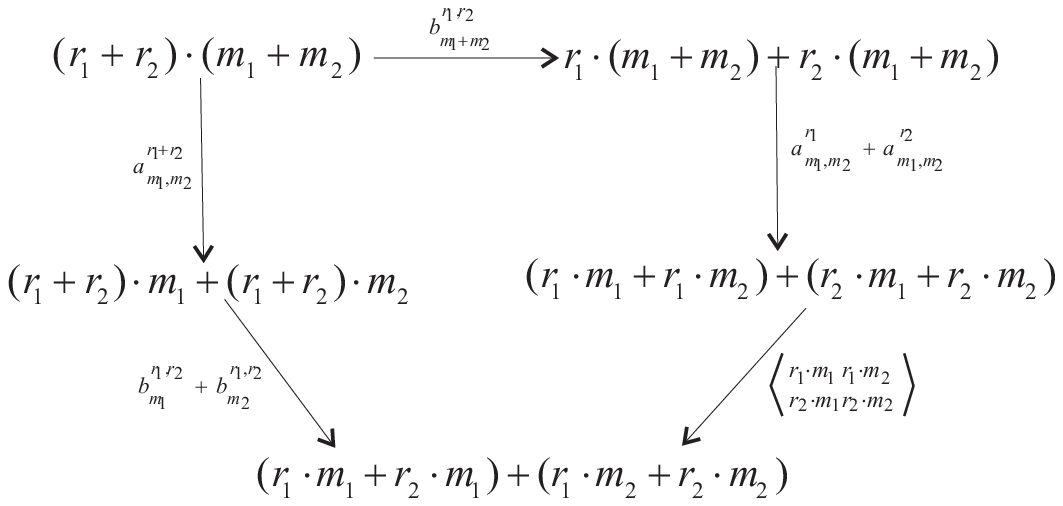}}{\footnotesize
Fig.22.}
\scalebox{0.9}[0.85]{\includegraphics{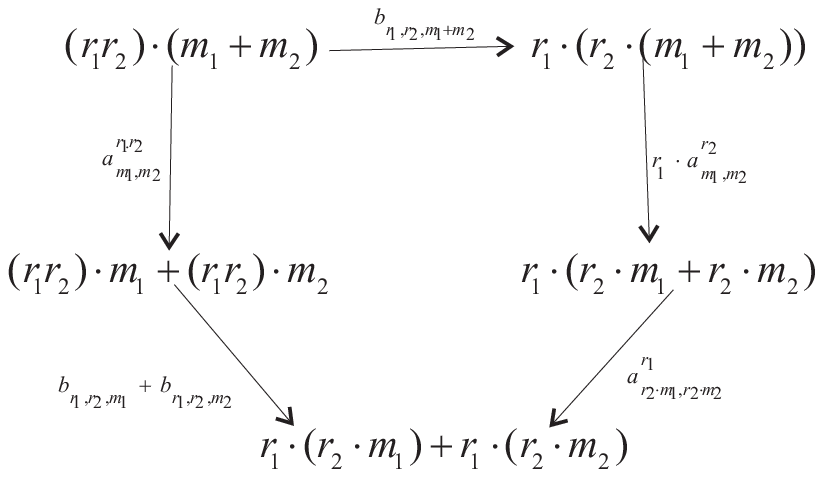}}{\footnotesize
Fig.23.}
\scalebox{0.9}[0.85]{\includegraphics{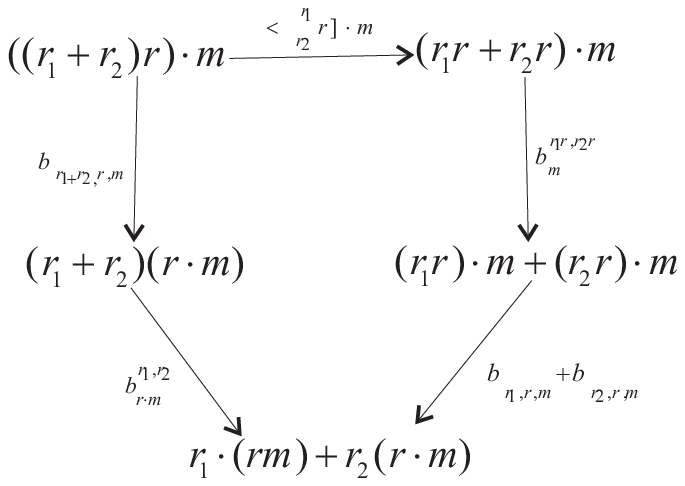}}{\footnotesize
Fig.24.}
\scalebox{0.9}[0.85]{\includegraphics{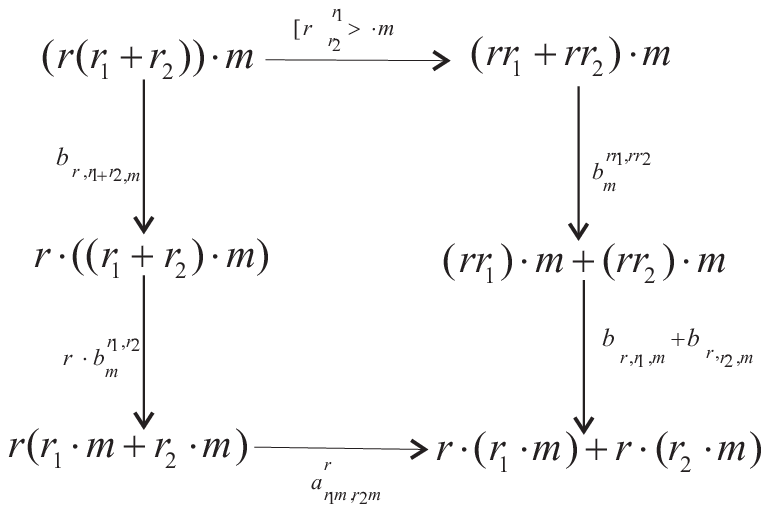}}{\footnotesize
Fig.25.}
\scalebox{0.9}[0.85]{\includegraphics{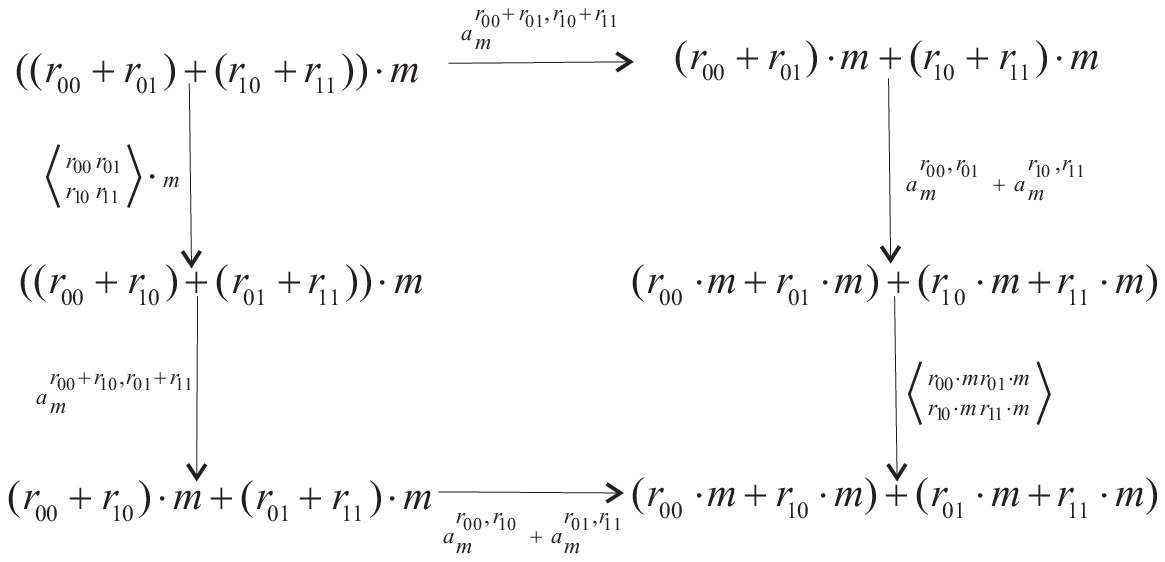}}{\footnotesize
Fig.26.}
\scalebox{0.9}[0.85]{\includegraphics{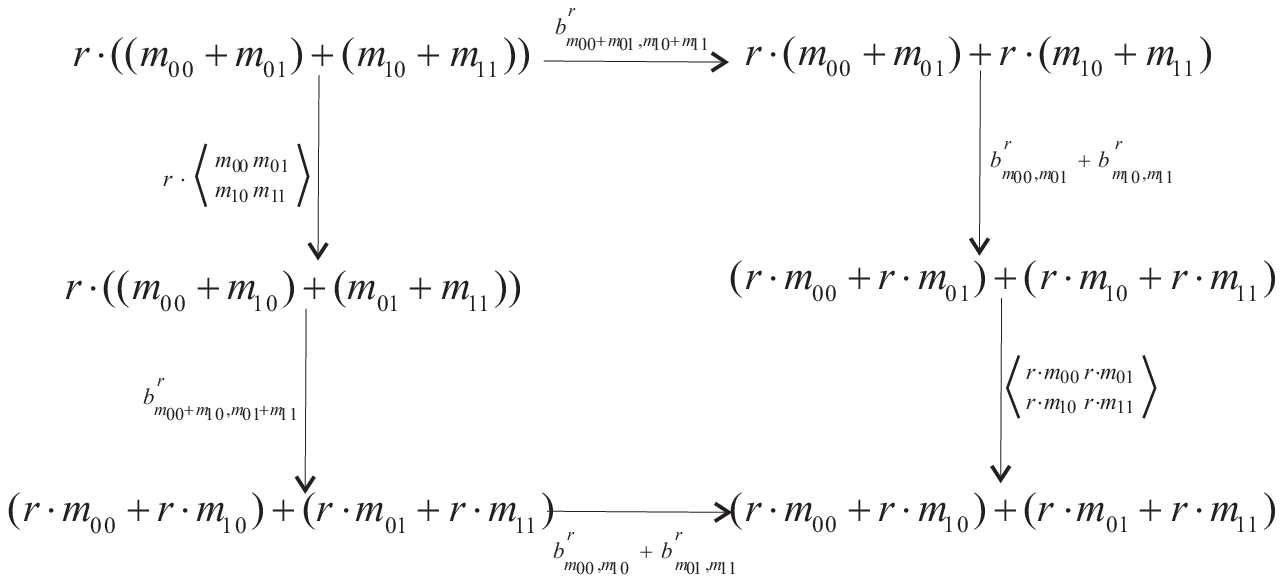}}{\footnotesize
Fig.27.}
\scalebox{0.9}[0.85]{\includegraphics{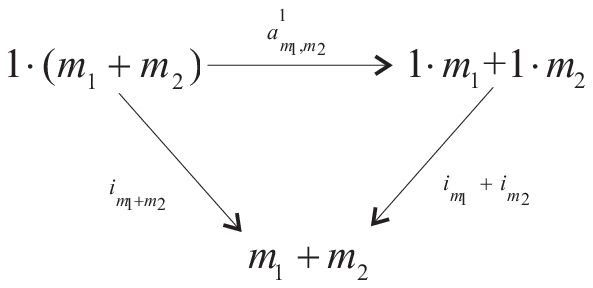}}{\footnotesize
Fig.28.}
\scalebox{0.9}[0.85]{\includegraphics{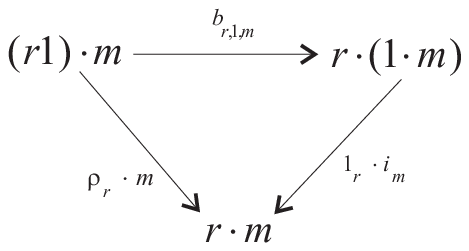}}{\footnotesize
Fig.29.}
\scalebox{0.9}[0.85]{\includegraphics{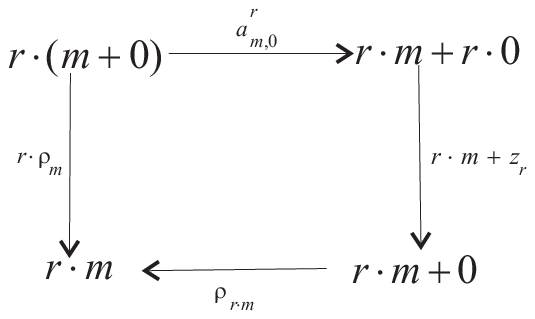}}{\footnotesize
Fig.30.}
\scalebox{0.9}[0.85]{\includegraphics{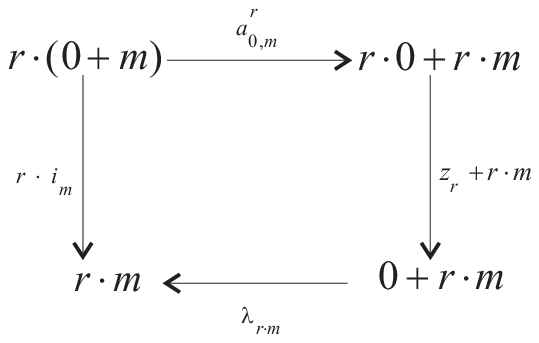}}{\footnotesize
Fig.31.}
\end{center}
\end{Def}

\begin{Def}
For two $\cR$-2-modules $\cM,\ \cN$, the $\cR$-homomorphism between
them is a quadruple ${\bf F}=(F,F_{+},F_{0},F_{2})$, where
$(F,F_{+},F_{0}):\cM\rightarrow \cN$ is a homomorphism of symmetric
2-groups $\cM$, $\cN$, $F_{2}$ is a natural isomorphism of the form
$$
F_{2}(r,m):F(r\cdot m)\rightarrow r\cdot F(m),
$$
for $r\in\cR,\ m\in\cM$, such that the following diagrams commute
for all possible objects of $\cR$ and $\cM$, respectively.
\begin{center}
\scalebox{0.9}[0.85]{\includegraphics{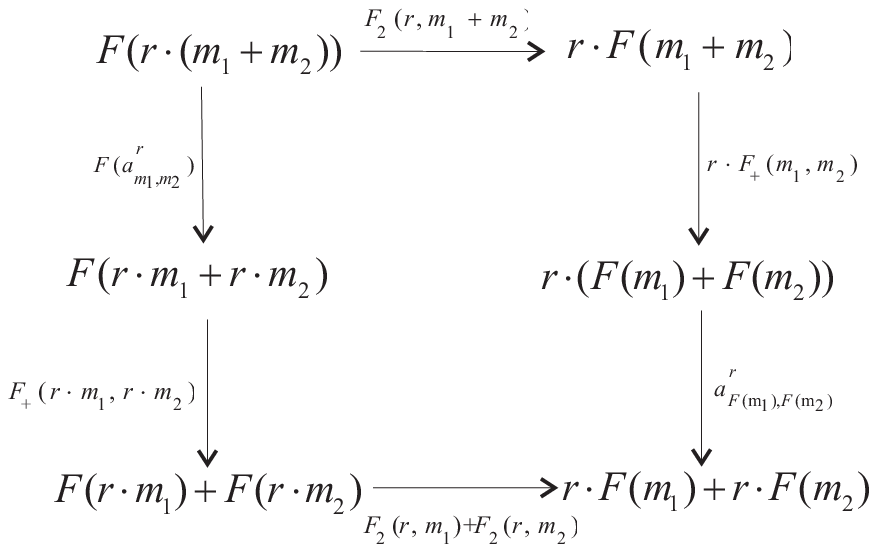}}{\footnotesize
Fig.32.}
\scalebox{0.9}[0.85]{\includegraphics{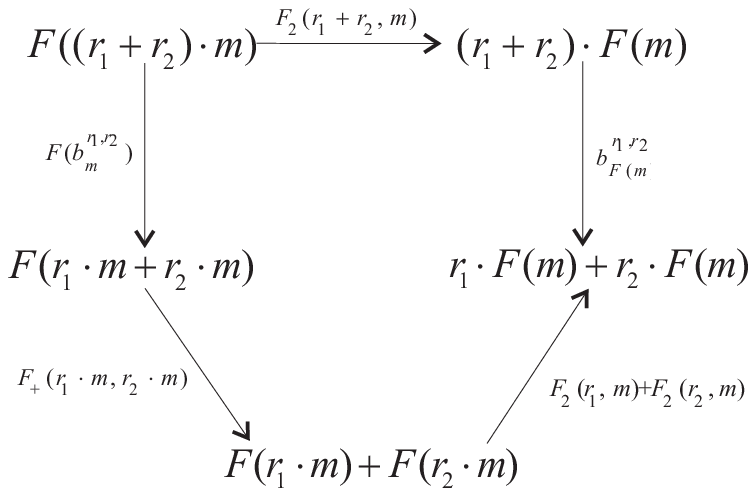}}{\footnotesize
Fig.33.}
\scalebox{0.9}[0.85]{\includegraphics{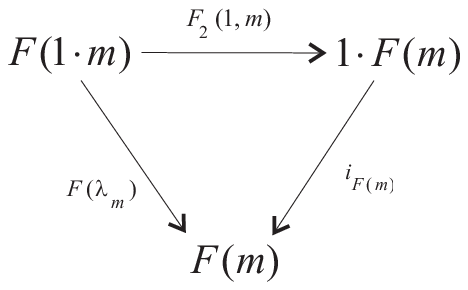}}{\footnotesize
Fig.34.}
\scalebox{0.9}[0.85]{\includegraphics{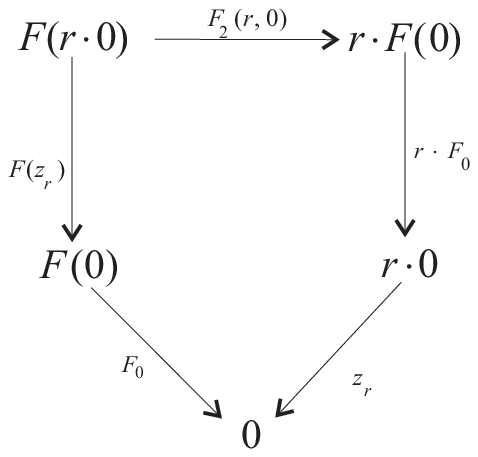}}{\footnotesize
Fig.35.}
\scalebox{0.9}[0.85]{\includegraphics{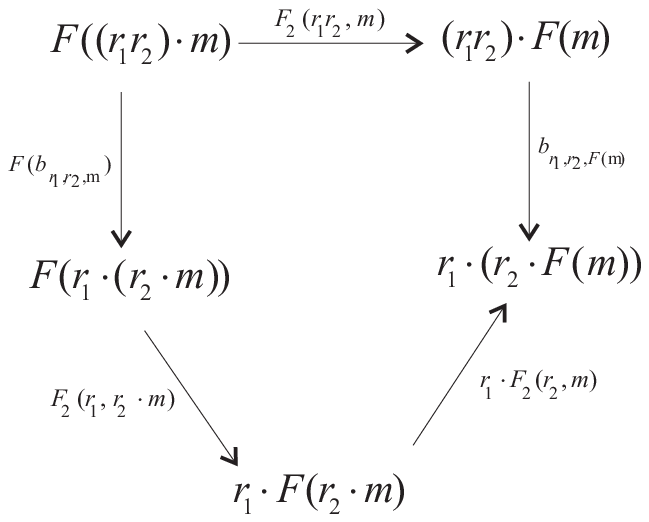}}{\footnotesize
Fig.36.}
\end{center}
\end{Def}

\begin{Def}
The morphism of two $\cR$-homomorphisms ${\bf
F}=(F,F_{+},F_{0},F_{2}),\ {\bf
G}=(G,G_{+},G_{0},G_{2}):\cM\rightarrow \cN$ of $\cR$-2-modules is a
natural transformation $\tau: F\Rightarrow G$ such that $\tau$ is
the morphism from $(F,F_{+},F_{0})$ to $(G,G_{+},G_{0})$ as
homomorphisms of symmetric 2-groups, moreover, the following diagram
commutes, for $r\in\cR,m\in\cM$:
\begin{center}
\scalebox{0.9}[0.85]{\includegraphics{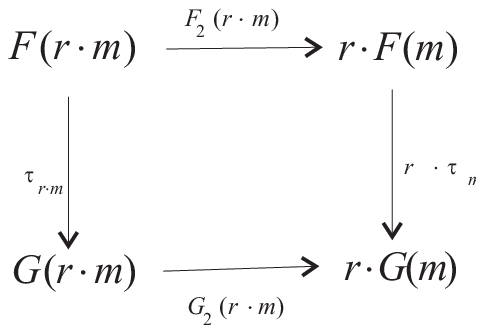}}{\footnotesize
Fig.37.}
\end{center}
\end{Def}
\noindent\textbf{Notation.} A homomorphism ${\bf F}=(F,F_+,F_0,F_2)$
of $\cR$-2-modules is in fact a functor $F$ satisfying some
compatible conditions(Fig.32.-36.). So, we always write $F$ for
abbreviation.

\begin{Thm}
Let $\cR$ be a 2-ring. The $\cR$-2-modules is a 2-category
 with objects $\cR$-2-modules, denoted by ($\cR$-2-Mod). Its 1-morphisms
 are $\cR$-homomorphisms between $\cR$-2-modules,
2-morphisms are morphisms of $\cR$-homomorphisms.
\end{Thm}
\begin{proof}
To prove this theorem we need to check the homomorphisms satisfy the
definition of a 2-category \cite 2.

\noindent\textbf{1)} 1-morphisms can be composed as a category.

Let $ F:\cM_{1}\rightarrow \cM_{2}$, $G:\cM_{2}\rightarrow \cM_{3}$
be $\cR$-homomorphisms. The composition of them is a composition
functor $H\triangleq G\circ F: \cM_{1}\rightarrow \cM_{3}$, together
with natural isomorphisms:
\begin{align*}
&H_{+}(m_{1},m_{2}):H(m_{1}+m_{2})=G(F(m_{1}+m_{2}))\xrightarrow[]{G(F_{+})}G(F(m_{1})+F(m_{2}))\\
&\hspace{4.6cm}\xrightarrow[]{G_{+}}G(F(m_{1}))+G(F(m_{2}))=H(m_{1})+H(m_{2}).\\
&H_{0}:H(0)=G(F(0))\xrightarrow[]{G(F_{0})}G(0)\xrightarrow[]{G_{0}}
  0.\\
&H_{2}(r,m):H(r\cdot m)=G(F(r\cdot m))
  \xrightarrow[]{G(F_{2})} G(r\cdot F(m))\xrightarrow[]{G_{2}}r\cdot G(F(m))=r\cdot H(m).
\end{align*}
\noindent\textbf{2)} 2-morphisms can be composed in two distinct
ways.
\begin{itemize}
\item Vertical composition of 2-morphisms.

Let $F,\ F^{'},\ F^{''}:\cM_{1}\rightarrow\cM_{2}$ be
$\cR$-homomorphisms of $\cR$-2-modules. $\tau:F\Rightarrow F^{'}$,
$\sigma:F^{'}\Rightarrow F^{''}$ are morphisms between them, the
composition of $\tau$ and $\sigma$ is
$\varepsilon=\sigma\circ\tau:F^{'}\Rightarrow F^{''}$ given by the
family of morphisms
$\{\varepsilon_{m}\triangleq\sigma_{m}\circ\tau_{m}:F(m)\rightarrow
F^{''}(m)\ \  \forall m\in\cM_{1}\}$ and such that the following
diagram commutes:
\begin{center}
\scalebox{0.9}[0.85]{\includegraphics{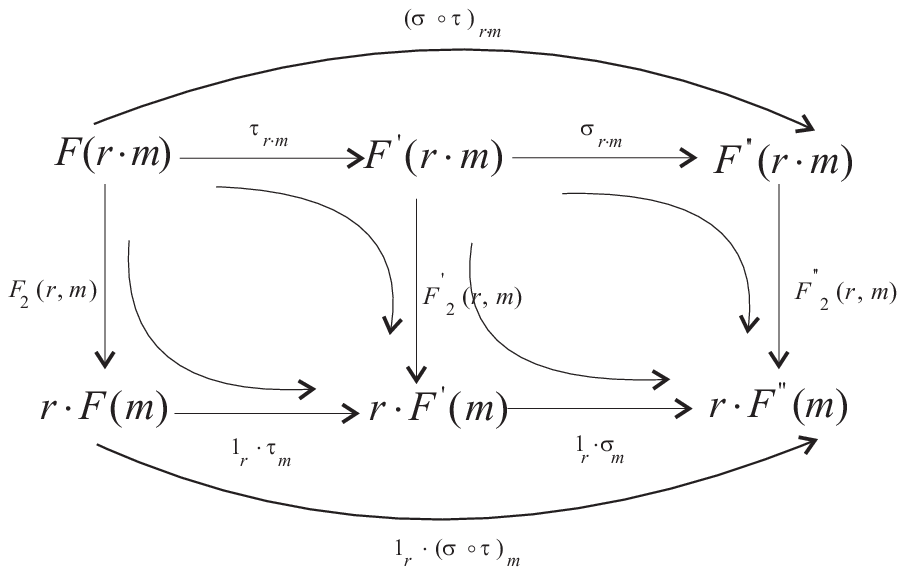}}
\end{center}
\item Horizontal composition of 2-morphisms.

Let $\alpha:F\Rightarrow F^{'}:\cM_{1}\rightarrow \cM_{2}$, $\beta:
G\Rightarrow G^{'}:\cM_{2}\rightarrow \cM_{3}$ be morphisms of
$\cR$-homomorphisms. The composition of $\alpha$ and $\beta$ is
$\tau=\beta\ast\alpha:( G\circ F)\Rightarrow ( G^{'}\circ F^{'})$
given by the family of morphisms
$\{\tau_{m}\triangleq\beta_{F^{'}(m)}\circ G(\alpha_{m}):(G\circ
F)(m)\rightarrow (G^{'}\circ F^{'})(m),\ \ \forall m\in\cM_{1}\}$,
such that the following diagrams commute:
\begin{center}
\scalebox{0.9}[0.85]{\includegraphics{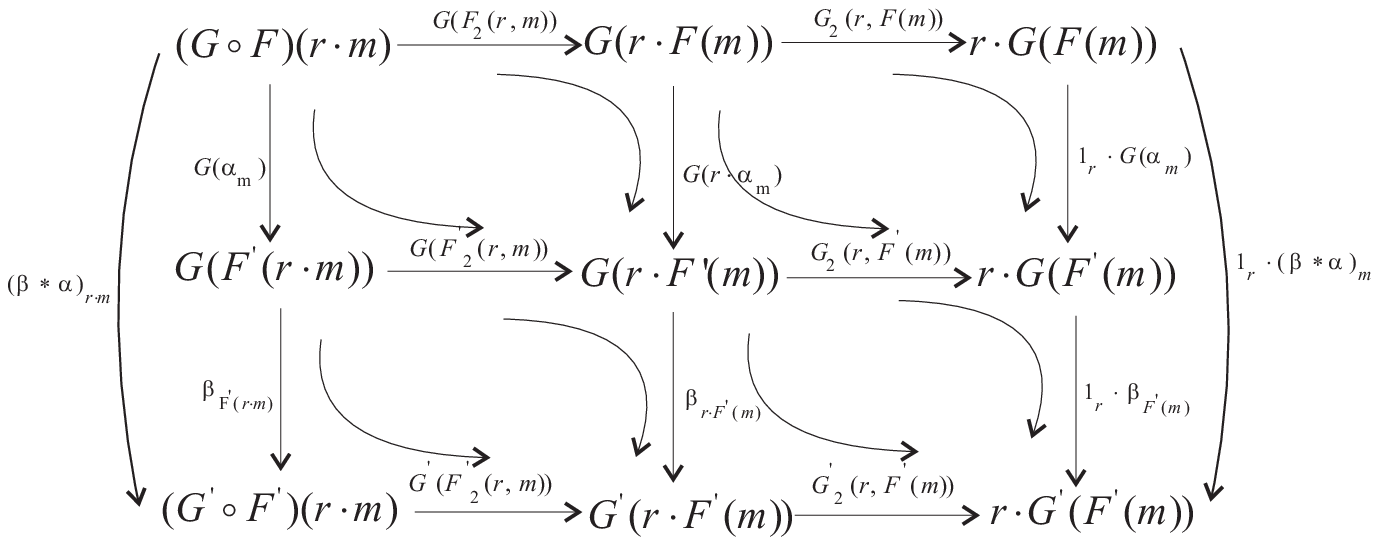}}
\end{center}
\end{itemize}
The above $\textbf{1)., 2).}$ satisfy \\
(i) Composition of 1-morphisms is associative, i.e. for
$\cR$-homomorphisms $\cM_{1}\xrightarrow[]{F} \cM_{2}
\xrightarrow[]{G} \cM_{3}\xrightarrow[]{H}\cM_{4}$, we have
$$
(H\circ G)\circ F=H\circ(G\circ F).
$$
In fact, there is a identity morphism $\varphi:( H\circ G)\circ
F\Rightarrow H\circ(G\circ F)$, defined by, $\varphi_m\triangleq
1_{m}:((H\circ G)\circ F)(m)=(H\circ G)(F(m))=H(G(F(m)))=H((G\circ
F)(m))=(H\circ(G\circ F))(m)$, for any $m\in\cM_1$, also $\varphi$
such that the Fig.4-6. and Fig.37. commute.

Moreover, for any $\cR$-2-module$\cM$, there is a $\cR$-homomorphism
${\bf1_{\cM}}=(1_{\cM},id,id):\cM\rightarrow \cM$ given by
$1_{\cM}(m)=m,\ 1_{\cM}(f)=f,$ for any object $m$ and morphism $f$
in $\cM$, such that, for any $ F:\cM\rightarrow \cN$, and $
G:\cK\rightarrow \cM$, we have $ F\circ 1_{\cM}= F$ and
$1_{\cM}\circ G=G$.\\
(ii) Vertical composition is associative, i.e. for any 2-morphisms
as in the diagram
\begin{center}
\scalebox{0.9}[0.85]{\includegraphics{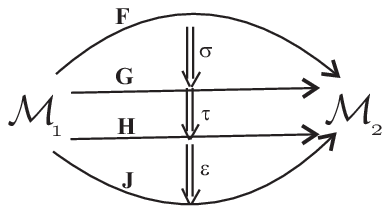}}
\end{center}
From the composition of morphisms in $\cM_{2}$, we have
$$
((\varepsilon\circ\tau)\circ\sigma)_{m}=(\varepsilon_{m}\circ\tau_{m})\circ\sigma_{m}=
\varepsilon_{m}\circ(\tau_{m}\circ\sigma_{m})=\varepsilon_{m}\circ(\tau\circ\sigma)_{m}.
$$
Moreover, for any $F:\cM_{1}\rightarrow \cM_{2}$, there exists
$1_{F}:F\Rightarrow F$ given by $1_{F}\triangleq1_{F}$, where
$(1_{F})_m=1_{Fm}$, such that, for any $\sigma:F\Rightarrow G$, we
have $\sigma\circ 1_{F}=\sigma$, since
$(\sigma\circ1_{F})_{m}=\sigma_{m}\circ
(1_{F})_{m}=\sigma_{m}\circ1_{Fm}=\sigma_{m}$, for any object
$m\in\cM_{1}$.\\
(iii) Horizontal composition is associative, i.e. for any
2-morphisms as in the following diagram:
\begin{center}
\scalebox{0.9}[0.85]{\includegraphics{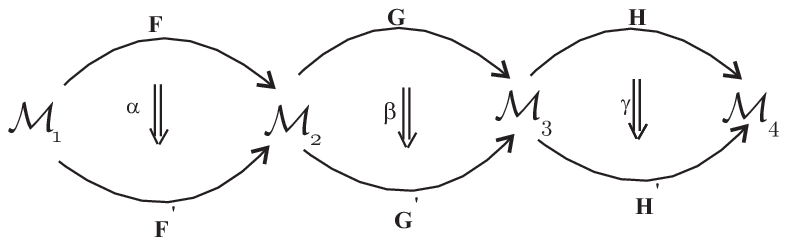}}
\end{center}
We need to check
$$
(\gamma\ast\beta)\ast\alpha=\gamma\ast(\beta\ast\alpha).
$$
In fact, for any object $m\in\cM_{1}$, from the definition of the
horizontal composition and natural transformation, the following
commutative ensures the associativity.
\begin{center}
\scalebox{0.9}[0.85]{\includegraphics{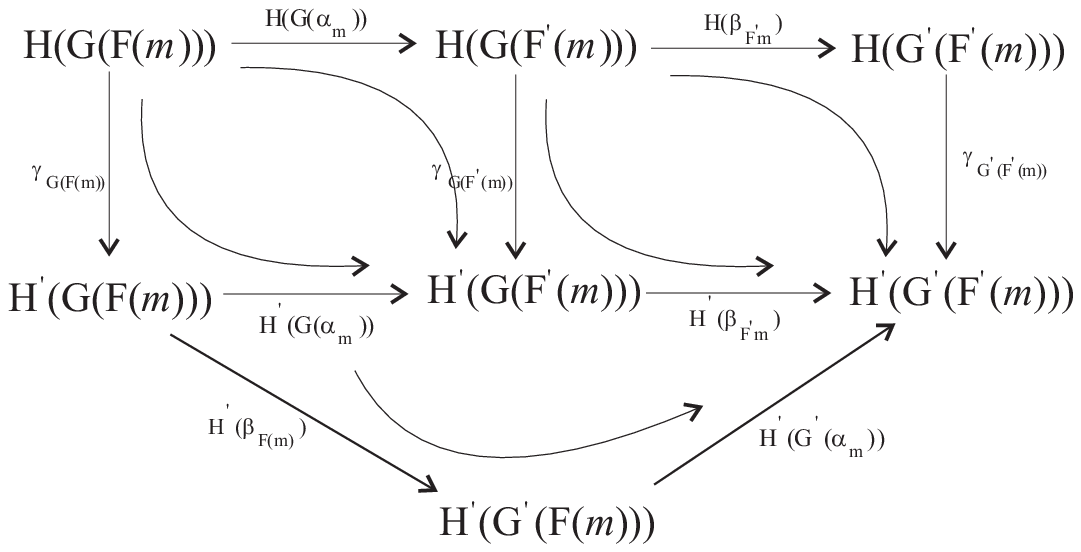}}
\end{center}
(iv) Vertical composition and horizontal composition of
2-morphisms satisfy the exchange law. \\
For any 1-morphisms and 2-morphisms as in the following diagram
\begin{center}
\scalebox{0.9}[0.85]{\includegraphics{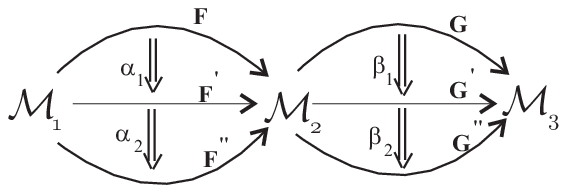}}
\end{center}
We need to check
$$
(\beta^{'}\circ\beta)\ast(\alpha^{'}\circ\alpha)=(\beta\ast\alpha)\circ(\beta^{'}\ast\alpha^{'}).
$$
In fact, for any $m\in\cM_{1}$, we have the following commutative
diagram
\begin{center}
\scalebox{0.9}[0.85]{\includegraphics{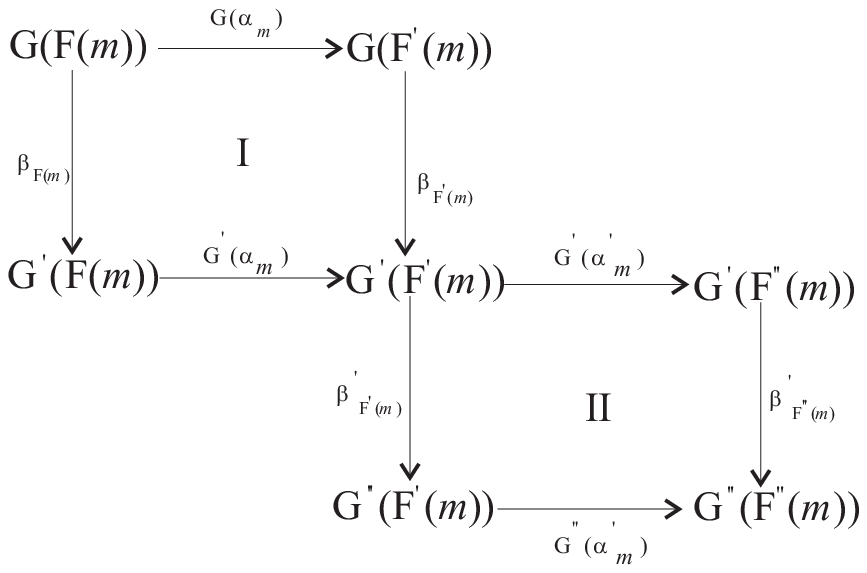}}
\end{center}
Note that, the above diagram commutes, since the small diagrams I
and II commute from the properties of natural
transformations $\beta,\ \beta^{'}.$\\
The above ingredients satisfy the definition of 2-category, so all
$\cR$-2-modules form a 2-category.
\end{proof}

\begin{Prop}
Let $\cA$ be a symmetric 2-group, $\cE nd \cA=\{ F:\cA\rightarrow
\cA\  is\ an\\
 endomorphism\ of\  \cA\}$. Then $\cE nd \cA$ is a
2-ring.
\end{Prop}

\begin{proof}
\noindent\textbf{Step 1.} $\cE nd\cA$ is a category consists of the
following data:

$\cdot$ Object is an endomorphism $F:\cA\rightarrow \cA$.

$\cdot$ Morphism is a morphism between two $\cR$-endomorphisms.
Composition of morphisms:
$$Hom(F,G)\times Hom(G,H)\longrightarrow Hom(F,H)$$
$$\hspace{2.1cm}(\tau,\sigma)\mapsto \sigma\circ \tau$$
is defined by,
$(\sigma\circ\tau)_{A}\triangleq\sigma_{A}\circ\tau_{A}$ as a
composition of morphisms in $\cA$, for any $A\in \cA$, which is a
morphism from $F$ to $H$, for $\tau, \sigma$ such that Fig.7.
commutes.

The above ingredients are subject to the following axioms:

(1) For $F\in obj(\cE nd\cA)$, there is an identity $ 1_{F}\in
Hom(F,F)$, defined by $(1_{F})_{A}\triangleq1_{F(A)}$, for any
$A\in\cA$, such that for any $\tau:F\Rightarrow G$, $\tau\circ
1_{F}=\tau$, since $(\tau\circ
1_{F})_{A}=\tau_{A}\circ1_{FA}=\tau_{A}$.

(2) Given morphisms $\tau:F\Rightarrow G$, $\sigma: G\Rightarrow H$,
and $\epsilon:H\Rightarrow J$, from the associativity of composition
of morphisms in $\cA$, we have
$$
(\epsilon\circ\sigma)\circ\tau=\epsilon\circ(\sigma\circ\tau).
$$
\noindent\textbf{Step 2.} $\cE$nd$\cA$ is a groupoid, i.e. for any
morphism $\tau:F\Rightarrow G:\cA\rightarrow\cA$ has an inverse. For
any object $a\in\cA$, $\tau_{a}:F(a)\rightarrow G(a)$ ia a morphism
in $\cA$, and $\cA$ is a groupoid, $\tau_{a}$ has an inverse
$(\tau_{a})^{*}$, so the inverse of $\tau$ is defined by
$(\tau^{*})_{a}\triangleq(\tau_{a})^{*}$.

\noindent\textbf{Step 3.} $\cE$nd$\cA$ is a monoidal category.
\begin{itemize}
\item There is a monmoidal structure on $\cE$nd$\cA$, i.e. there is a
 bifunctor
\begin{align*}
&\hspace{2.5cm}+:\cA\times\cA\longrightarrow\cA\\
&\hspace{3cm}( F,\ \ G)\mapsto F+G,\\
&(\tau: F\Rightarrow G,\ \tau^{'}:F^{'}\Rightarrow G^{'})\mapsto
\tau+\tau^{'}:F+F^{'}\Rightarrow G+G^{'}
\end{align*}

where $F+G$ and $\tau+\tau^{'}$ are
given as follows:\\
$(F+G)(a)\triangleq F(a)+G(a),\ (\tau+\tau^{'})_{a}\triangleq
\tau_{a}+\tau^{'}_{a}$ under the monoidal structure of
$\cA$.\\
$(F+G)_{+}:(F+G)(a+b)\rightarrow (F+G)(a)+(F+G)(b)$ is the
composition
\begin{align*}
 &(F+G)(a+b)=F(a+b)+G(a+b)\xrightarrow[]{F_{+}(a,b)+G_{+}(a,b)}
(Fa+Fb)+(Ga+Gb)\\
             &\hspace{0.5cm} \xrightarrow[]{\langle^{Fa\ Fb}_{Ga\ Gb}\rangle}(Fa+Ga)+(Fb+Gb)=(F+G)(a)+(F+G)(b).
\end{align*}
$(F+G)_{0}:(F+G)(0)\rightarrow 0$ is the composition
$$
(F+G)(0)=F0+G0\xrightarrow[]{F_{0}+G_{0}} 0+0\xrightarrow[]{l_0}0.
$$
Since $F$ and $G$ are endomorphisms of symmetric 2-group $\cA$, and
$\tau$, $\tau^{'}$ are morphisms of homomorphisms, so $F+G\in\cE
nd\cA$, and $\tau+\tau^{'}$ is morphism of homomorphism.
\item There is a unit object ${\bf 0}=(0,1_0,1_0):\cA\rightarrow \cA$, given by
$0(a)\triangleq 0$, $\forall a\in\cA$, where 0 is the unit object of
$\cA$.
\item There are natural isomorphisms:
\begin{align*}
&<F,G,H>:(F+G)+H\Rightarrow F+(G+H),\\
&\hspace{1.8cm}l_{F}: 0+F\rightarrow F,\\
&\hspace{1.8cm}r_{F}:F+0\rightarrow F
\end{align*}
given by:
\begin{align*}
&<F,G,H>_{a}\triangleq
<Fa,Ga,Ha>:((F+G)+H)(a)=(Fa+Ga)+Ha\\
&\hspace{3cm}\rightarrow
Fa+(Ga+Ha)=(F+(G+H))(a),\\
&(l_{F})_{a}:(0+F)(a)=0(a)+Fa=0+Fa\xrightarrow[]{l_{Fa}}
Fa,\\
&(r_{F})_{a}:(F+0)(a)=Fa+0(a)=Fa+0\xrightarrow[]{r_{Fa}}Fa,
\end{align*}
such that Fig.1.-2. commute.
\end{itemize}

\noindent\textbf{Step 3.} $\cE nd\cA$ is a symmetric 2-group.

(i)\ Every object of $\cE nd\cA$ is invertible, i.e. $\forall\
F\in\cE nd\cA, \exists F^* \in\cE nd\cA$, and natural isomorphism
$\eta_{F}:F^{*}+F\rightarrow 0$. In fact, let $F^*(a)\triangleq
(F(a))^*,\ (\eta_{F})_{a}\triangleq \eta_{Fa}$, where $(F(a))^*$ is
the inverse of $F(a)$, and $\eta_{Fa}:(F(a))^*+F(a)\rightarrow 0$ is
a natural isomorphism in $\cA$.

(ii)\ For any $F,\ G\in\cE nd\cA$, there is a natural isomorphism
$c_{F,G}:F+G\rightarrow G+F$, given by $(c_{F,G})_{a}\triangleq
c_{Fa,Ga}:Fa+Ga\rightarrow Ga+Fa$, and since $\cA$ is symmetric, so
$c_{Fa,Ga}$ is an isomorphism, such that $c_{Fa,Ga}\circ
c_{Ga,Fa}=id$, then $c_{F,G}\circ c_{G,F}=id$.

\noindent\textbf{Step 4.} $\cE nd\cA$ is a 2-ring.

We need to give another monoidal structure on $\cE nd\cA$, satisfy
some compatibilities.

(i)\ There are a bifunctor $\cdot:\cA\times\cA\rightarrow\cA,$ by
$F\cdot G\triangleq F\circ G,\ \tau\cdot\tau^{'}\triangleq
\tau^{'}\ast\tau$, $F\circ G$ is the composition of homomorphisms,
see Remark 2.(ii), and $\tau^{'}\ast\tau$ is the horizontal
composition of 2-morphisms in the 2-category (2-SGp).

Using the same method in proof of Theorem 1, we know that $F\cdot
G\in\cE nd\cA$, $\tau\cdot\tau^{'}$ is the morphism in $\cE nd\cA.$

(ii)\ The unit object of $\cE nd\cA$ under the monoidal structure
$\cdot$ is $(I,id,id)$, which is for any object $a$ and morphism f
in $\cA$, $I(a)\equiv a,\ I(f)\equiv f.$

(iii)\  There are natural isomorphisms:
\begin{align*}
&[F,G,H]:(F\cdot G)\cdot H\Rightarrow F\cdot(G\cdot H),\\
&\hspace{1.1cm}\lambda_{F}:I\cdot F\Rightarrow F,\\
&\hspace{1.1cm}\rho_{F}:F\cdot I\Rightarrow F,\\
&\hspace{0.5cm}[F^{G}_{G^{'}}>:F\cdot(G+G^{'})\Rightarrow F\cdot G+F\cdot G^{'},\\
&\hspace{0.3cm}<_{F^{'}}^{F}G]:(F+F^{'})\cdot G\Rightarrow F\cdot
G+F^{'}\cdot G
\end{align*}
given by:
\begin{align*}
&[F,G,H]_{a}=id:((F\cdot G)\cdot H)(a)=(F\cdot(G\cdot
H))(a),\\
&\hspace{0.8cm}(\lambda_{F})_{a}=id:(I\cdot F)(a)=(I\circ F)(a)=I(Fa)=Fa, \\
&\hspace{0.8cm}(\rho_{F})_{a}=id:(F\cdot I)(a)=(F\circ I)a=F(I(a))=Fa, \\
&({[F^{G}_{G^{'}}>})_{a}=id:F\cdot(G+G^{'})(a)=F((G+G^{'})(a))=
F(Ga+G^{'}a)\xrightarrow[]{F_{+}}(F\cdot G+F\cdot G^{'})(a),\\
&(< _{F^{'}}^{F}G])_{a}=id:((F+F^{'})\cdot
G)(a)=(F+F^{'})G(a)=F(Ga)+F^{'}(Ga)=(F\cdot G)(a)+(F^{'}\cdot G)(a).
\end{align*}

The above ingredients satisfy the following conditions:

1)\ ($I,\ [-,-,-],\ \lambda_-,\ \rho_-$) constitute a monoidal
structure (obviously, since all of them are identities.)

2)\ Fig.8-14. commute for all possible objects of $\cE nd\cA$, from
the definition of homomorphisms.
\end{proof}
\begin{Rek}
For an Ann-category $\cA$, the category $End(\cA)$ is also an
Ann-category(\cite{21}), then is a 2-ring(\cite{14}).
\end{Rek}

Similarly, as the 1-dimensional representation of rings\cite{15}, we
give
\begin{Def}
 A representation of a 2-ring $\cR$ is a 2-ring homomorphism
$F:\cR\rightarrow \cE nd\cA$ of 2-rings, where $\cA$ is a symmetric
2-group.
\end{Def}

\begin{Prop} Let $\cM$ be a symmetric 2-group.
Every representation $F:\cR\rightarrow \cE nd\cM$ equips $\cM$ with
the structure of a $\cR$-2-module. Conversely, every $\cR$-2-module
$\cM$ determines a representation $F:\cR\rightarrow \cE nd\cM$.
\end{Prop}

\begin{proof}
For a representation $F:\cR\rightarrow \cE nd\cM$, we will prove
$\cM$ is a $\cR$-2-module. \\
There are a bifunctor $\cdot:\cR\times \cM\rightarrow \cM$,
defined by $r\cdot m\triangleq F(r)(m).$ and natural isomorphisms:\\
\begin{align*}
&a_{m,n}^{r}:r\cdot(m+n)=F(r)(m+n)\xrightarrow[]{F(r)_+}
F(r)(m)+F(r)(n)=r\cdot m+r\cdot
n,\\
&b_{m}^{r,s}:(r+s)\cdot m=F(r+s)(m)\xrightarrow[]{F_+}
(F(r)+F(s))(m)=F(r)m+F(s)m=r\cdot m+s\cdot m,\\
&b_{r,s,m}:(rs)\cdot m=F(rs)(m)\xrightarrow[]{F_{\cdot}}
(F(r)F(s))(m)=F(r)(F(s)(m))=r\cdot(s\cdot m),\\
&i_{m}:I\cdot m=F(I)(m)\xrightarrow[]{F_1} I(m)=m,\\
&z_{r}:r\cdot0=F(r)(0)\xrightarrow[]{F(r)_0}0.
\end{align*}
It is easy to check Fig.18-31. commute, since F is homomorphism of
2-rings.

Conversely, if $\cM$ is a $\cR$-2-module, by Proposition 1, we know
that $\cE nd\cM$ is a 2-ring. Now we give a homomorphism
\begin{align*}
&F:\cR\longrightarrow \cE nd\cM\\
&\hspace{0.9cm}r\mapsto F(r)(m)\triangleq r\cdot m,\\
\end{align*}
where $r\cdot m$ is the operation of $\cR$ on $\cM$.\\
Next we will prove F is a homomorphism of 2-rings.
\begin{itemize}
\item $F:\cR\rightarrow\cE nd\cM$ is a functor.

(i) $F(r)\in\cE nd\cM$, i.e. F(r) is an endomorphism of
symmetric 2-group $\cM$.\\
In fact, $F(r)\triangleq r\cdot$ is a functor, for fixed $r\in\cR$
from $\cdot$ being a bifunctor. Moreover, there are natural
isomorphisms:
\begin{align*}
&F(r)_{+}(m_{1},m_{2})\triangleq
a^{r}_{m_{1},m_{2}}:F(r)(m_{1}+m_{2})=r\cdot(m_{1}+m_{2})\rightarrow
r\cdot m_{1}+r\cdot m_{2}\\
&\hspace{7cm}=F(r)(m_{1})+F(r)(m_{2}),\\
&F(r)_{0}\triangleq z_{r}:F(r)(0)=r\cdot 0\rightarrow 0.
\end{align*}
such that the following diagrams commute.
\begin{center}
\scalebox{0.9}[0.85]{\includegraphics{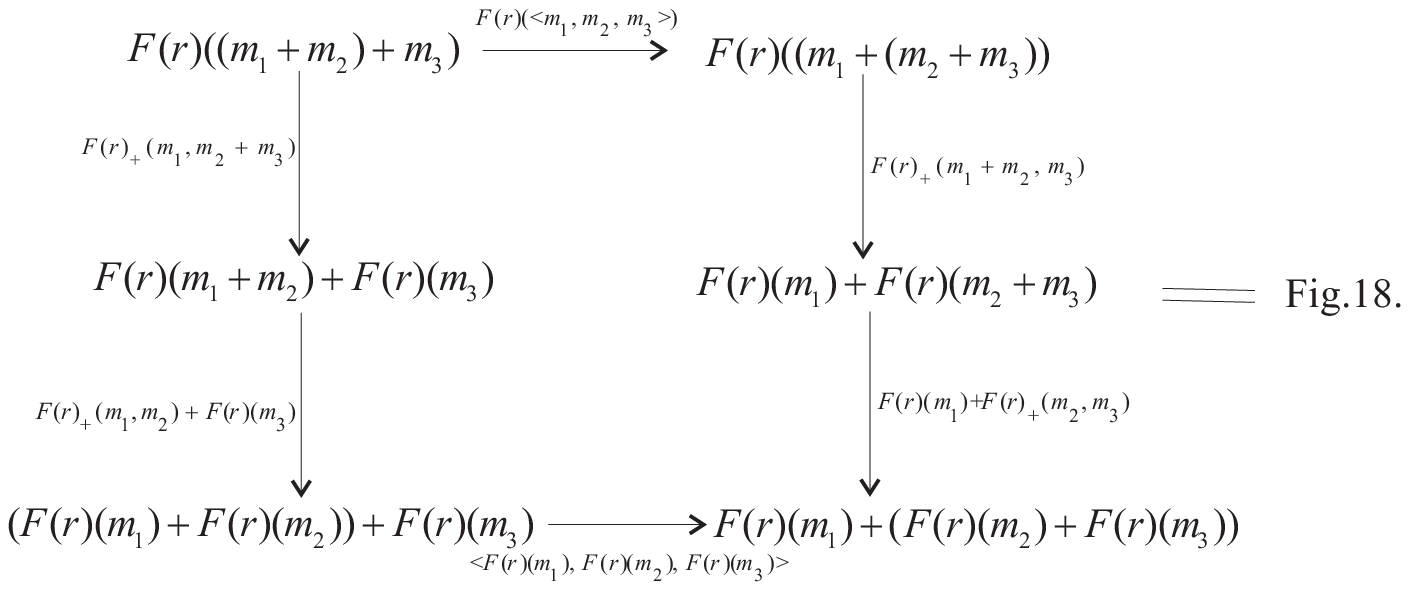}}
\scalebox{0.9}[0.85]{\includegraphics{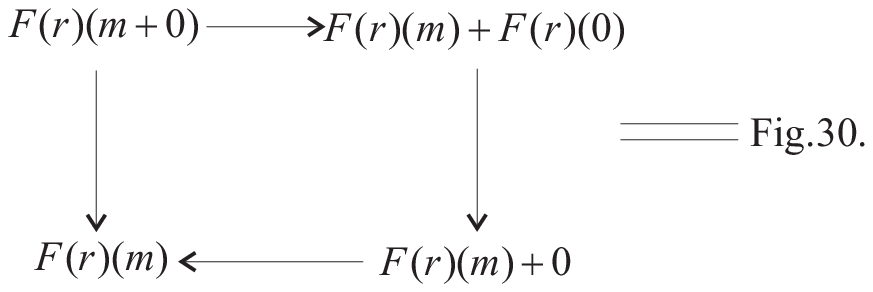}}
\scalebox{0.9}[0.85]{\includegraphics{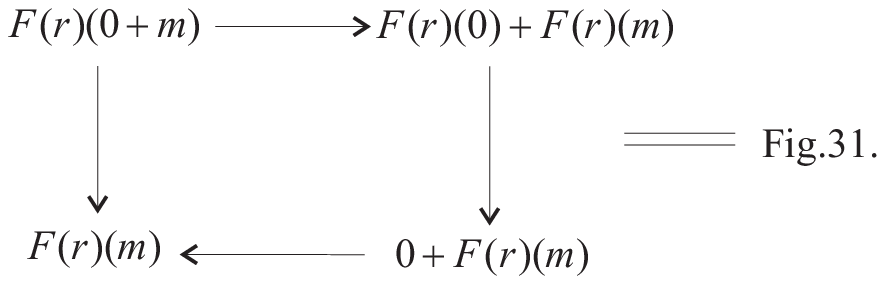}}
\end{center}
(ii) For any morphism $f:r_{1}\rightarrow r_{2}$ in $\cR$,
$F(f):F(r_{1})\rightarrow F(r_{2})$ is a morphism in $\cE nd\cM$.

(iii) For any morphisms $r_{1}\xrightarrow[]{f_{1}}
r_{2}\xrightarrow[]{f_{2}}r_{3}$, we have $F(f_{2}\circ
f_{1})=F(f_{2})\circ F(f_{1})$, and also for any object $r\in\cR$,
$F(1_{r})=1_{F(r)}$.

In fact, by the properties of bifunctor
$\cdot:\cR\times\cM\rightarrow \cM$, it is easy to prove (ii),
(iii).
\item There are natural isomorphisms:
\begin{align*}
&F_{+}(r_{0},r_{1}):F(r_{0}+r_{1})\rightarrow
F(r_{0})+F(r_{1}),\\
&F_{.}(r_{0},r_{1}):F(r_{0}r_{1})\rightarrow F(r_{0})F(r_{1}),\\
&F_{1}:F(1)\rightarrow 1
\end{align*}
given by the following ways: for any object $m\in\cM$,
\begin{align*}
&F_{+}(r_{0},r_{1})_{m}\triangleq
b_{r_{0},r_{1}}^{m}:F(r_{0}+r_{1})(m)=(r_{0}+r_{1})\cdot
m\rightarrow r_{0}\cdot m+r_{1}\cdot m\\
&\hspace{3.5cm}=(F(r_{0})+F(r_{1}))(m),\\
&F_{.}(r_{0},r_{1})_{m}\triangleq b_{r_{0},r_{1},m}:F(r_{0}r_{1})(m)
=(r_{0}r_{1})\cdot m\rightarrow r_{0}\cdot(r_{1}\cdot m)=F(r_{0})F(r_{1}),\\
&(F_{1})_{m}\triangleq i_{m}:F(1)(m)=1\cdot m\rightarrow m=1(m).
\end{align*}
\item The above ingredients satisfy the following conditions:

(i) $(F,F_{+},F_{0})$, $(F,F_{\cdot},F_{1})$ are monoidal functors
with respect to the monoidal structures on 2-rings $\cR$
and $\cE nd\cM$.\\

(ii)\ The following diagrams commute:
\begin{center}
\scalebox{0.9}[0.85]{\includegraphics{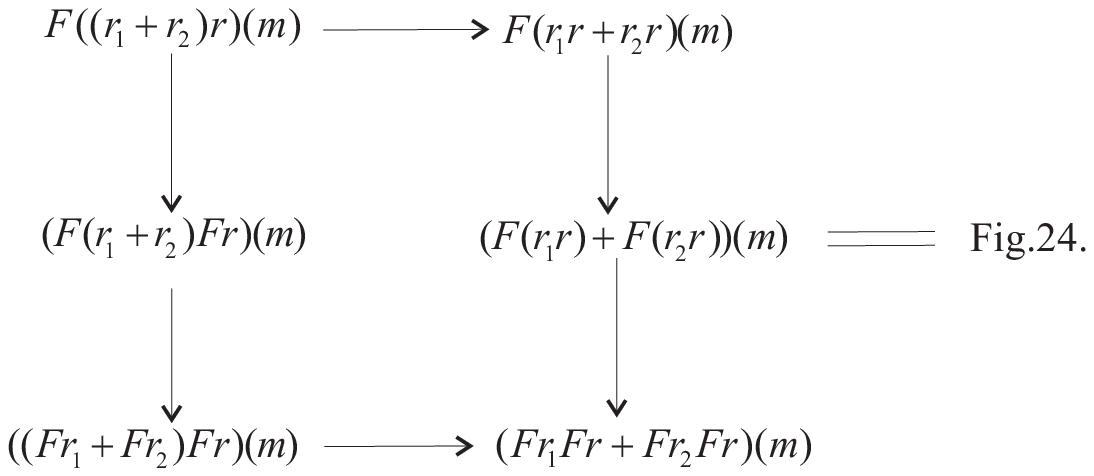}}
\end{center}
\begin{center}
\scalebox{0.9}[0.85]{\includegraphics{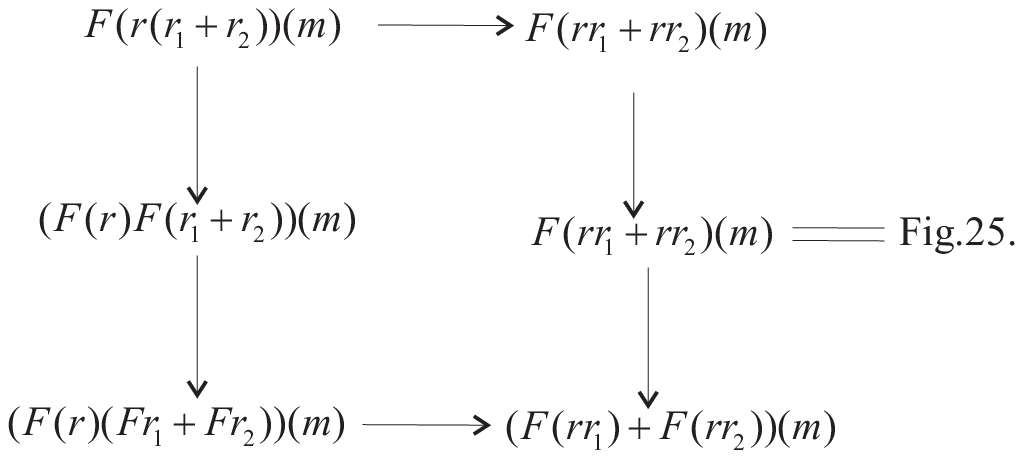}}
\end{center}
\end{itemize}
\end{proof}

\section{($\cR$-2-Mod) is a 2-Abelian $Gpd$-category }

In this section, we will show that $(\cR$-2-Mod) is a 2-abelian
$Gpd$-category under the the definition of 2-abelian $Gpd$-category
given in \cite {18}. First we will give several definitions similar
as \cite{5,6,7,11,18}.

\begin{Def}\cite {16}
Let $\cM,\ \cN$ be two $\cR$-2-modules, where $\cR$ is an 2-ring.
The functor $0:\cM\rightarrow\cN$ which sends each morphism to the
identity of the unit object of $\cN$, is a $\cR$-homomorphism,
called the zero-morphism.
\end{Def}

\begin{Def}
Let $F:\cA\rightarrow \cB$ be a 1-morphism in( $\cR$-2-Mod). The
kernel of F is a triple $(KerF,e_F,\varepsilon_F)$, where $KerF$ is
an $\cR$-2-module, $e_{F}:KerF\rightarrow \cA$ is a 1-morphism, and
$\varepsilon_{F}:F\circ e_{F}\Rightarrow 0$ is a 2-morphism,
satisfies the universal property in the following sense:\\
For given $\cK\in obj(\cR$-2-Mod), a 1-morphism $G:\cK\rightarrow
\cA$, and a 2-morphism $\varphi:F\circ G\Rightarrow 0$, there exist
a 1-morphism $G^{'}:\cK\rightarrow KerF$ and a 2-morphism
$\varphi^{'}:e_{F}\circ G^{'}\Rightarrow G$, such that $\varphi^{'}$
is compatible with $\varphi$ and $\varepsilon_{F}$, i.e. the
following diagram commutes:
\begin{center}
\scalebox{0.9}[0.85]{\includegraphics{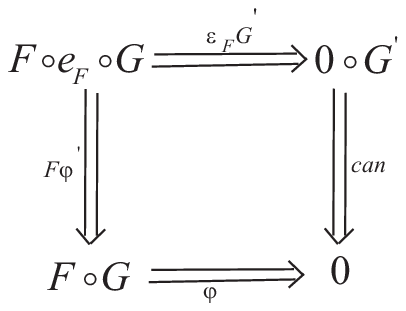}}{\footnotesize
Fig.38.}
\end{center}
Moreover, if $G^{''}$ and $\varphi^{''}$ satisfy the same conditions
as $G^{'}$ and $\varphi^{'}$, then there exists a unique 2-morphism
$\psi:G^{''}\Rightarrow G^{'}$, such that
\begin{center}
\scalebox{0.9}[0.85]{\includegraphics{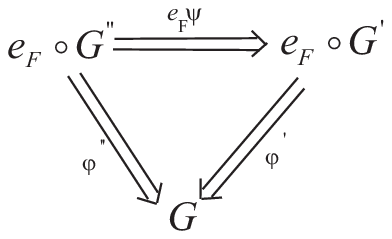}}{\footnotesize
Fig.39.}
\end{center}
\end{Def}

\begin{Thm}
For any 1-morphism $F:\cA\rightarrow\cB$ in ($\cR$-2-Mod), the
kernel of F exists.
\end{Thm}
\begin{proof}
We will construct the kernel of $F$ as follows:
\begin{itemize}
\item There is a category $KerF$ consists of:

$\cdot$  Object is a pair $(A,a)$, where $A\in obj(\cA)$,
$a:F(A)\rightarrow 0$ is a morphism in $\cB$.

$\cdot$ Morphism $f:(A,a)\rightarrow (A^{'},a^{'})$ is a morphism
$f:A\rightarrow A^{'}$ in $\cA$, such that $a^{'}\circ F(f)=a.$

$\cdot$ Composition of morphisms. Given morphisms
$(A,a)\xrightarrow[]{f}(A^{'},a^{'})\xrightarrow[]{f^{'}}(A^{''},a^{''})$,
their composition $f^{'}\circ f:A\rightarrow A^{''}$ is just the
composition of morphisms in $\cA$, such that $a^{''}\circ
F(f^{'}\circ f)=a^{''}\circ (F(f^{'})\circ F(f)=a^{'}\circ F(f)=a.$

The above ingredients satisfy the following axioms:

(1) For any $(A,a)\in obj(KerF)$, there exists
$1_{(A,a)}\triangleq1_{A}:(A,a)\rightarrow (A,a)$ such that
any$f:(A,a)\rightarrow (A^{'},a^{'})$, $f\circ 1_{(A,a)}=f$.

(2) Composition is associative: Given morphisms
$$
A_{1}\xrightarrow[]{f_{1}}A_{2}\xrightarrow[]{f_{2}}A_{3}\xrightarrow[]{f_{3}}A_{4}.
$$
We have $(f_{3}\circ f_{2})\circ f_{1}=f_{3}\circ (f_{2}\circ
f_{1})$ in $\cA$, so it is true in $KerF$.

\item $KerF$ is a groupoid with zero object $(0,F_{0})$, where 0 is unit object of $\cA,\\
F_{0}:F(0)\rightarrow 0$ is a morphism in $\cA$.\\
For any morphism $f:(A,a)\rightarrow (A_{1},a_{1})$ of $KerF$, is a
morphism $f:A\rightarrow A_{1}$ in $\cA$ such that $a_{1}\circ
Ff=a$. $\cA$ is an $\cR$-2-module, as a morphism in $\cA$,
$f:A\rightarrow A_{1}$ has inverse $g:A_{1}\rightarrow A$ such that
$g\circ f=1_{A}$. Then $g:(A_{1},a_{1}\rightarrow (A,a)$, such that
$g\circ f=1_{(A,a)}$ and $a\circ Fg=a_{1}$, so is the inverse of $f$
in $KerF$.
\item $KerF$ is a monoidal category.

There is a bifunctor
\begin{align*}
&\hspace{4cm}+:KerF\times KerF\longrightarrow KerF\\
&\hspace{4.4cm}((A_{1},a_{1}),A_{2},a_{2}))\mapsto (A,a)\triangleq(A_{1},a_{1})+(A_{2},a_{2})\\
&((A_{1},a_{1})\xrightarrow[]{f_{1}}(A_{1}^{'},a_{1}^{'}),(A_{2},a_{2})\xrightarrow[]{f_{2}}(A_{2}^{'},a_{2}^{'}))\mapsto
(A_{1},a_{1})+(A_{1}^{'},a_{1}^{'})\xrightarrow[]{f_{1}+f_{2}}
(A_{2},a_{2})+(A_{2}^{'},a_{2}^{'})
\end{align*}
where $A\triangleq A_{1}+A_{2},\ a$ is the composition
$$
F(A)=F(A_{1}+A_{2})\xrightarrow[]{F_{+}}
FA_{1}+FA_{2}\xrightarrow[]{a_{1}+a_{2}}0+0\xrightarrow[]{l_0}0,
$$
$f_{1}+f_{2}:A_{1}+A_{1}^{'}\rightarrow A_{2}+A_{2}^{'}$ is
 an addition of morphisms in $\cA$ under monoidal structure of $\cA$, such that the following diagram
 commutes:
\begin{center}
\scalebox{0.9}[0.85]{\includegraphics{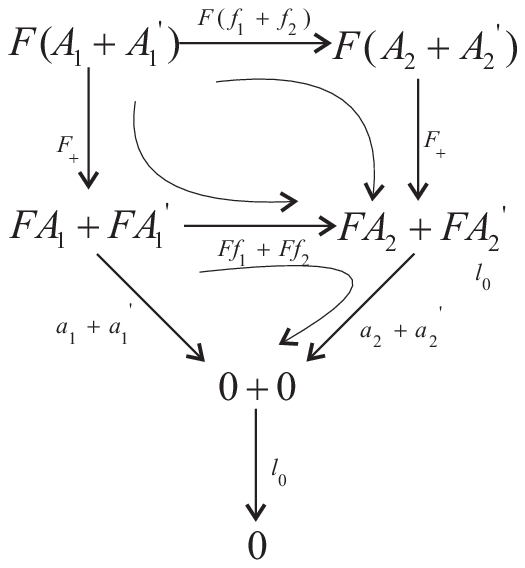}}
\end{center}
Then $f_{1}+f_{2}$ is a morphism in $KerF$.

Moreover, there are natural isomorphisms:

(i) $<(A_{1},a_{1}),(A_{2},a_{2}),(A_{3},a_{3})>
\triangleq<A_{1},A_{2},A_{3}>:((A_{1},a_{1})+(A_{2},a_{2}))+(A_{3},a_{3})
=((A_{1}+A_{2})+A_{3},a)\longrightarrow
(A_{1}+(A_{2}+A_{3}),a^{'})=(A_{1},a_{1})+((A_{2},a_{2})+(A_{3},a_{3})),$
such that the following diagram commutes:
\begin{center}
\scalebox{0.9}[0.85]{\includegraphics{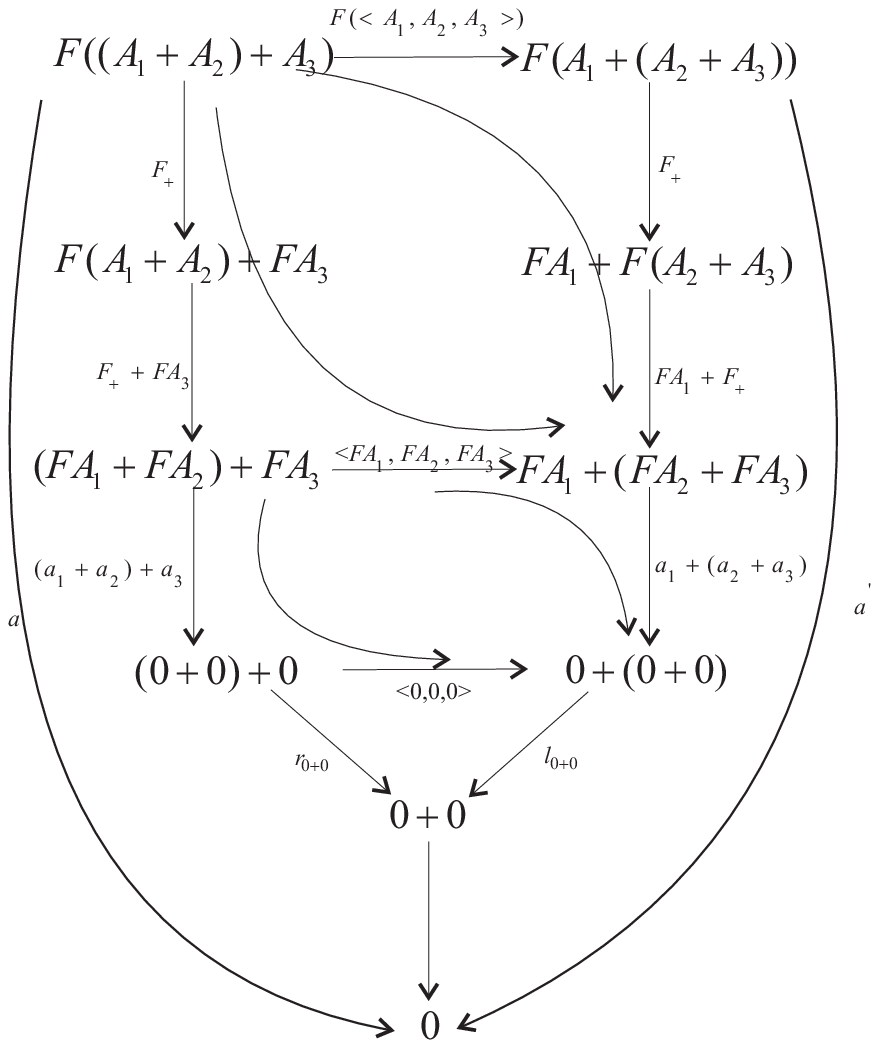}}
\end{center}
(ii) $l_{(A,a)}\triangleq
l_{A}:(0,F_{0})+(A,a)=(0+A,a^{'})\longrightarrow (A,a)$, such that
\begin{center}
\scalebox{0.9}[0.85]{\includegraphics{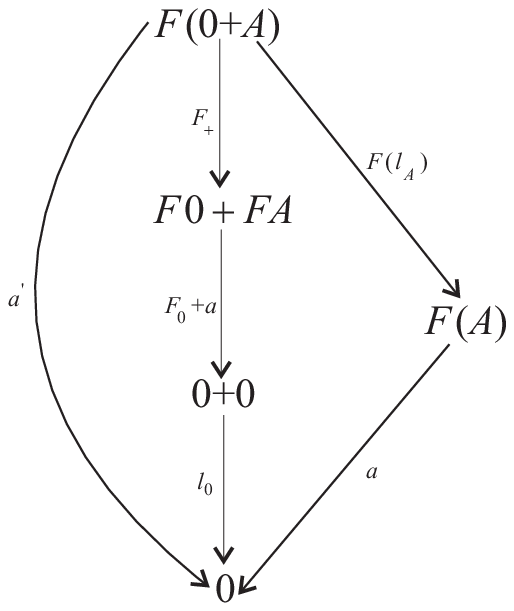}}
\end{center}
commutes.\\
(iii) $r_{(A,a)}\triangleq
r_{A}:(A,a)+(0,F_{0})=(A+0,a^{'})\longrightarrow (A,a)$, such that
\begin{center}
\scalebox{0.9}[0.85]{\includegraphics{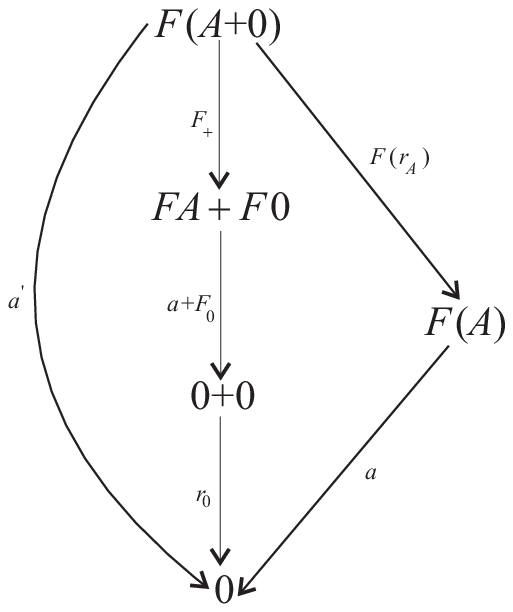}}
\end{center}
commutes.

From the definition of morphism in $KerF$ and $\cA$ being a monoidal
category, the above isomorphisms satisfy the Mac Lane coherence
conditions.

\item\ $KerF$ is a symmetric 2-group.

$\cdot$ For any $(A,a)\in obj(KerF)$, $A\in obj(\cA),\
a:F(A)\rightarrow 0$. Since $\cA$ is a symmetric 2-group, so $A$ has
an inverse $A^{*}$, and natural isomorphism
$\eta_{A}:A^{*}+A\rightarrow 0$. Let $a^{*}$ be a composition
$F(A^{*})\cong (F(A))^{*}\rightarrow 0^{*}\cong 0$. We have
$(A,a)^{*}\triangleq(A^{*},a^{*})$, and natural isomorphism
$\eta_{(A,a)}\triangleq\eta_{A}: (A,a)^{*}+(A,a)\rightarrow
(0,F_{0})$. It is easy to check $\eta_{(A,a)}$ is a morphism in
$KerF$.

$\cdot$ For any objects $(A_{1},a_{1}),(A_{2},a_{2})\in KerF$, there
is a natural isomorphism
$$
c_{(A_{1},a_{1}),(A_{2},a_{2})}:(A_{1},a_{1})+(A_{2},a_{2})\rightarrow
(A_{2},a_{2})+(A_{1},a_{1}),
$$
given by
$$
c_{A_{1},A_{2}}:A_{1}+A_{2}\rightarrow A_{2}+A_{1},
$$
and since $c_{A_{1},A_{2}}\circ c_{A_{2},A_{1}}=id$, so
$c_{A_{1},A_{2}}\circ c_{A_{2},A_{1}}=id$.
\item $KerF$ is an $\cR$-2-module.

There is a bifunctor
\begin{align*}
&\hspace{2.7cm}\cdot:\cR\times KerF\longrightarrow KerF\\
&\hspace{3.5cm}(r,(A,a))\mapsto r\cdot(A,a)\triangleq (r\cdot
A,a^{'})\\
&(r_{1}\xrightarrow[]{\varphi} r_{2},(A_{1},a_{1})\xrightarrow[]{f}
(A_{2},a_{2}))\mapsto r_{1}\cdot
(A_{1},a_{1})\xrightarrow[]{\varphi\cdot f}r_{2}\cdot(A_{2},a_{2})
\end{align*}
where $a^{'}:F(r\cdot A)\xrightarrow[]{F_{2}}r\cdot
F(A)\xrightarrow[]{r\cdot a}r\cdot 0\xrightarrow[]{z_{r}}0$,
$\varphi\cdot f$ is defined under the operation of $\cR$ on $\cA$,
and such that the following diagram commutes:
\begin{center}
\scalebox{0.9}[0.85]{\includegraphics{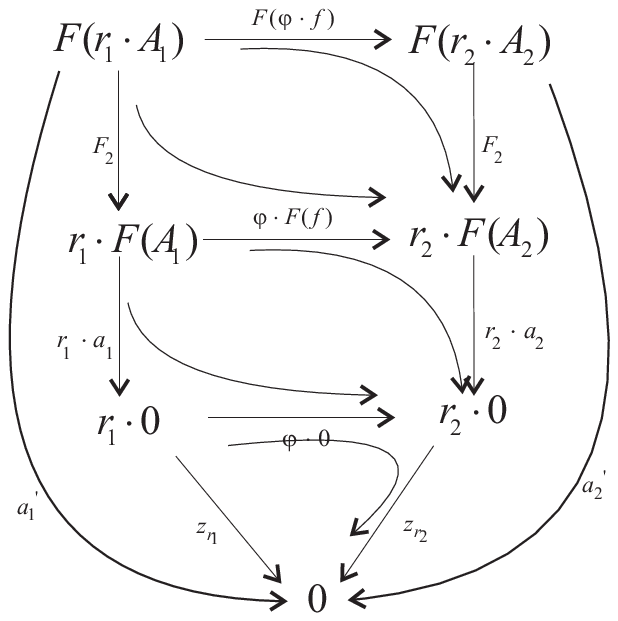}}
\end{center}
Moreover, there are natural isomorphisms:

(1) $a_{(A_{1},a_{1}),(A_{2},a_{2})}^{r}:r\cdot
((A_{1},a_{1})+(A_{2},a_{2}))\rightarrow r\cdot (A_{1},a_{1})+r\cdot
(A_{2},a_{2})$, as follows:\\
$r\cdot ((A_{1},a_{1})+(A_{2},a_{2}))\triangleq (r\cdot
(A_{1}+A_{2}),a^{'})$, $r\cdot (A_{1},a_{1})+r\cdot
(A_{2},a_{2})\triangleq (r\cdot A_{1}+r\cdot A_{2},a^{''})$,
$a_{(A_{1},a_{1}),(A_{2},a_{2})}^{r}\triangleq a_{A_{1},A_{2}}^{r}$,
where $a^{'}$ and $a^{''}$ are the compositions as in the following
diagrams, such that the diagram commutes:
\begin{center}
\scalebox{0.9}[0.85]{\includegraphics{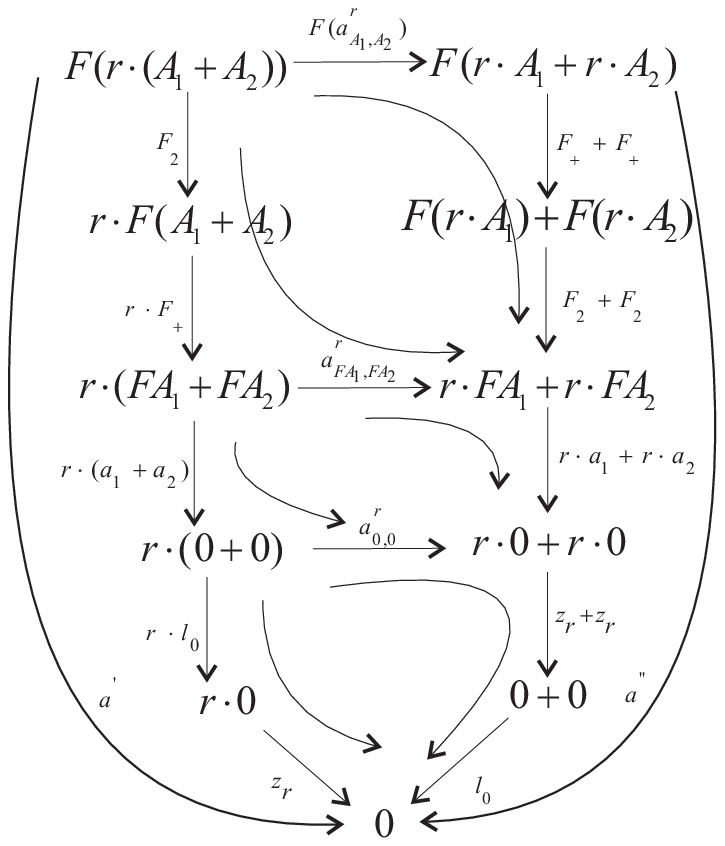}}
\end{center}
(2) $b_{(A,a)}^{r_{1},r_{2}}:(r_{1}+r_{2})\cdot (A,a)\rightarrow
r_{1}\cdot(A,a)+r_{2}\cdot(A,a),$ as follows:\\
$(r_{1}+r_{2})\cdot (A,a)\triangleq ((r_{1}+r_{2})\cdot A,a^{'})$,
$r_{1}\cdot(A,a)+r_{2}\cdot(A,a)\triangleq (r_{1}\cdot A+r_{2}\cdot
A,a^{''})$, $b_{(A,a)}^{r_{1},r_{2}}\triangleq b_{A}^{r_{1},r_{2}}$,
where $a^{'}$ and $a^{''}$ are the compositions as in the following
diagrams, such that the diagram commutes:
\begin{center}
\scalebox{0.9}[0.85]{\includegraphics{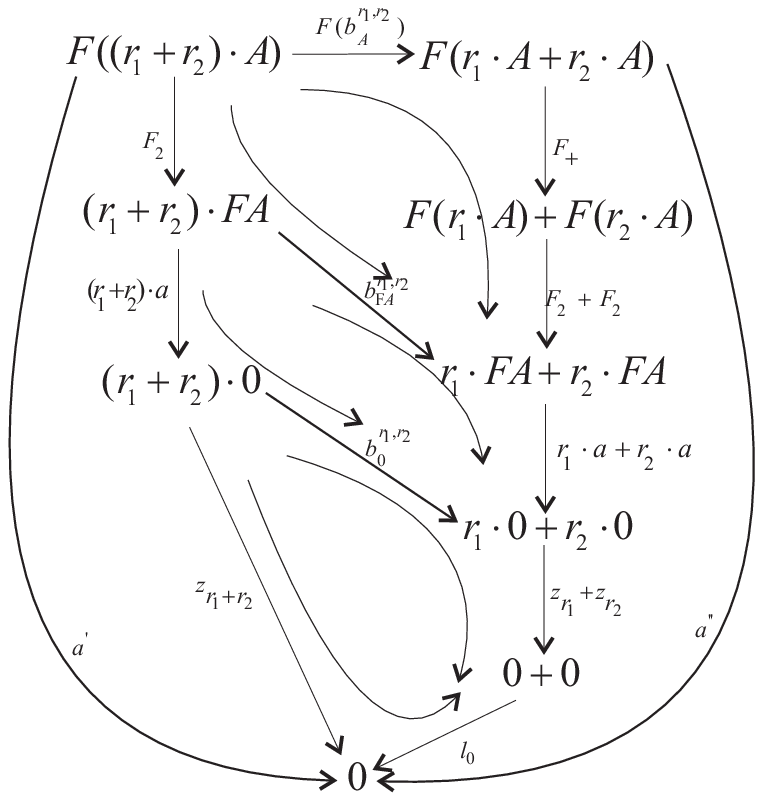}}
\end{center}

(3) $b_{r_{1},r_{2},(A,a)}:(r_{1}r_{2})\cdot (A,a)\rightarrow
r_{1}\cdot(r_{2}\cdot (A,a))$, as follows:\\
$(r_{1}r_{2})\cdot (A,a)\triangleq ((r_{1}r_{2})\cdot A,a^{'})$,
$(r_{1}\cdot(r_{2}\cdot (A,a))\triangleq (r_{1}\cdot(r_{2}\cdot
A),a^{''})$, $b_{r_{1},r_{2},(A,a)}\triangleq b_{r_{1},r_{2},A}$,
where $a^{'}$ and $a^{''}$ are the compositions as in the following
diagrams, such that the diagram commutes:
\begin{center}
\scalebox{0.9}[0.85]{\includegraphics{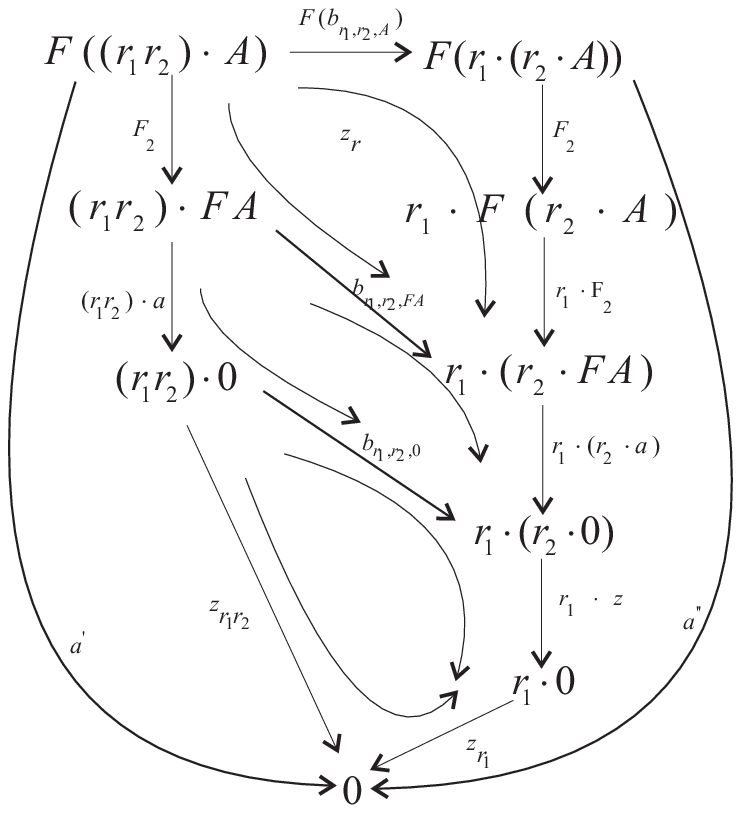}}
\end{center}

(4) $i_{(A,a)}:1\cdot(A,a)\rightarrow (A,a)$, as follows:\\
$1\cdot(A,a)\triangleq (1\cdot A,a^{'})$, $i_{(A,a)}\triangleq
1_{A}$, where $a^{'}$ is the composition as in the following
diagrams, such that the diagram commutes:
\begin{center}
\scalebox{0.9}[0.85]{\includegraphics{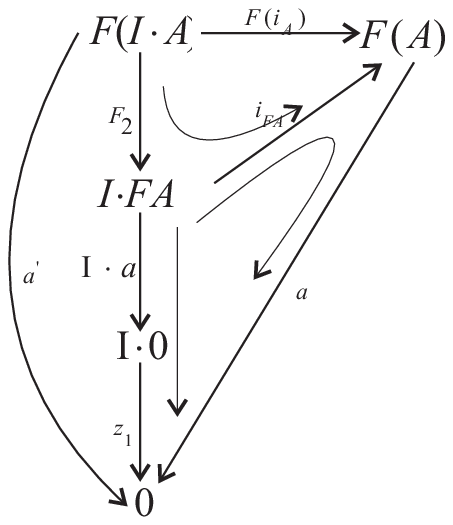}}
\end{center}

(5)\ $z_{r}:r\cdot(0,F_{0})\rightarrow (0,F_{0})$, as follows:\\
$r\cdot(0,F_{0})\triangleq(r\cdot 0,F_{0}^{'}),\ z_{r}$ is just the
natural isomorphism in $\cA$, such that the following diagram
commutes:
\begin{center}
\scalebox{0.9}[0.85]{\includegraphics{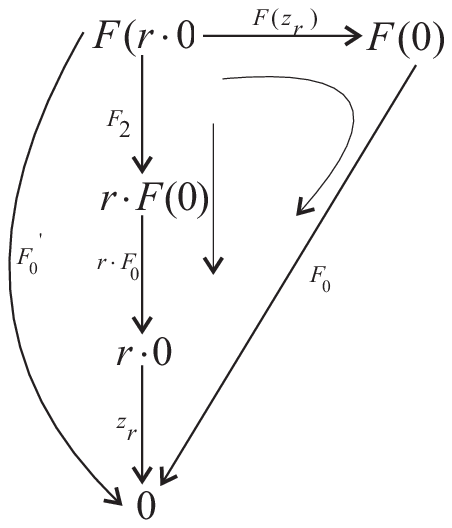}}
\end{center}

The above natural isomorphisms, together with the basic facts about
$\cR$-2-module $\cA$ satisfy Fig.18.-31.commute.

\item $KerF$ is a kernel of $F:\cA\rightarrow\cB$.

(i) There is a functor
\begin{align*}
&\hspace{1.3cm}e_{F}:KerF\longrightarrow \cA\\
&\hspace{2.2cm}(A,a)\mapsto A\\
&(A_{1},a_{1})\xrightarrow[]{f}(A_{2},a_{2})\mapsto
f:A_{1}\rightarrow A_{2}
\end{align*}
such that, for $(A_{1},a_{1})\xrightarrow[]{f_{1}}
(A_{2},a_{2})\xrightarrow[]{f_{2}}(A_{3},a_{3})$, we have
$e_{F}(f_{2}\circ f_{1})=f_{2}\circ f_{1}=e_{F}(f_{2})\circ
e_{F}(f_{1})$.\\
And also there are natural isomorphisms:
\begin{align*}
&(e_{F})_{+}=id:e_{F}((A_{1},a_{1})+(A_{2},a_{2})=A_{1}+A_{2}=
e_{F}((A_{1},a_{1})+e_{F}((A_{2},a_{2})),\\
&(e_{F})_{0}=id:e_{F}(0,F_{0})=0,\\
&(e_{F})_{2}=id:e_{F}(r\cdot(A,a))=r\cdot A=r\cdot e_{F}(A,a).
\end{align*}
Obviously, Fig.32.-36. commute.

(ii)\ There is a 2-morphism:
$$\varepsilon_{F}:F\circ e_{F}\Rightarrow 0$$
given by, $(\varepsilon_{F})_{(A,a)}\triangleq a:(F\circ
e_{F})(A,a)=F(A)\rightarrow 0$, for any $(A,a)\in obj(KerF)$.

(iii)\ $(KerF,e_{F},\varepsilon_{F})$ satisfies the universal
property.\\
For any $\cK\in obj(\cR$-2-Mod), 1-morphism $G:\cK\rightarrow \cA$,
and 2-morphism $\varphi:F\circ G\Rightarrow 0$. There is a
1-morphism:
\begin{align*}
&\hspace{0.5cm}G^{'}:\cK\longrightarrow KerF\\
&\hspace{1.2cm}A\mapsto (G(A),\varphi_{A})\\
&A_{1}\xrightarrow[]{f}A_{2}\mapsto
(G(A_{1}),\varphi_{A_{1}})\xrightarrow[]{G(f)}(G(A_{2}),\varphi_{A_{2}})
\end{align*}
It is easy to see $G^{'}$ is well-defined and $G^{'}$ is an $\cR$-homomorphism.\\
Define a 2-morphism
$$
\varphi^{'}:e_{F}\circ G^{'}\Rightarrow G
$$
by
$$
\varphi^{'}_{A}\triangleq1_{GA}:(e_{F}\circ
G^{'})(A)=e_{F}(GA,\varphi_{A})=GA\rightarrow GA ,
$$
for any $A\in
obj(\cK)$.

$\varphi^{'}$  is compatible with $\varphi$ and $\varepsilon_{F}$,
i.e. Fig.38. commutes.

If $G^{''},\varphi^{''}$ satisfy the same conditions as
$G^{'},\varphi^{'}$.\\
To define $\psi:G^{''}\Rightarrow G^{'}$, we need to define, for any
$A\in obj(\cK)$, a morphism $\psi_{A}:G^{''}(A)\rightarrow
G^{'}(A)$, which is defined by $\psi_{A}\triangleq\varphi^{''}_{A}$
from the commutative diagram Fig.39..

Obviously, $\psi^{''}$ is unique for $\varphi^{''}$ is.
\end{itemize}
\end{proof}

\begin{Def}
Let $F:\cA\rightarrow \cB$ be a 1-morphism in($\cR$-2-Mod). The
cokernel of F is the triple $(CokerF,p_F,\pi_F)$, where $CokerF$ is
an $\cR$-2-module, $p_{F}:\cB\rightarrow CokerF$ is a 1-morphism,
and $\pi_{F}:p_{F}\circ F\Rightarrow 0$ is a 2-morphism, satisfies
the universal
property in the following sense:\\
Given an $\cR$-2-module $\cK$, a 1-morphism $G:\cB\rightarrow \cK$,
and a 2-morphism $\varphi:G\circ F\Rightarrow 0$, there exist a
1-morphism $G^{'}:CokerF\rightarrow \cK$ and a 2-morphism
$\varphi^{'}:G^{'}\circ p_{F} \Rightarrow G$, such that
$\varphi^{'}$ is compatible with $\varphi$ and $\pi_{F}$, i.e. the
following diagram commutes.
\begin{center}
\scalebox{0.9}[0.85]{\includegraphics{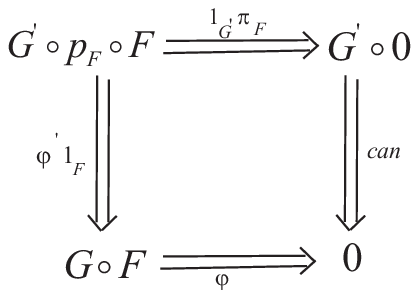}}{\footnotesize
Fig.40.}
\end{center}
Moreover, if $G^{''}$ and $\varphi^{''}$ satisfy the same conditions
as $G^{'}$ and $\varphi^{'}$, then there exists a unique 2-morphism
$\psi:G^{''}\Rightarrow G^{'}$, such that
\begin{center}
\scalebox{0.9}[0.85]{\includegraphics{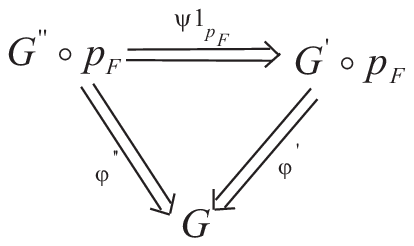}}{\footnotesize
Fig.41.}
\end{center}
commutes.
\end{Def}

\begin{Thm}
For any 1-morphism $F:\cA\rightarrow\cB$ in ($\cR$-2-Mod), the
cokernel of F exists.
\end{Thm}

\begin{proof}
We need to construct an the cokernel of $F$ as follows:
\begin{itemize}
\item Let $CokerF$ be a category consisting of:

$\cdot$ Objects are those of $\cB$.

$\cdot$ Morphism of $B_{1}\rightarrow B_{2}$ is the equivalence
class of $(f,A)$, denoted by $[f,A]$, where $A\in obj(\cA)$,
$f:B_{1}\rightarrow B_{2}+F(A)$, and for two morphisms $(f,A),\
(f^{'},A^{'}):B_{1}\rightarrow B_{2}$ are equivalent if and only if
there exists an isomorphism $\alpha:A\rightarrow A^{'}$ in $\cA$,
such that $(1_{B_{2}}+F(\alpha))\circ f=f^{'}.$

$\cdot$\  Composition of morphisms in $CokerF$.\\
Let
$B_{1}\xrightarrow[]{[f_{1},A_{1}]}B_{2}\xrightarrow[]{[f_{2},A_{2}]}B_{3}$
be morphisms in $CokerF$. Then
$[f_{2},A_{2}]\circ[f_{1},A_{1}]\triangleq
[f,A]$, where $A=A_{1}+A_{2},\ f$ is the composition \\
\begin{align*}
& B_{1}\xrightarrow[]{f_{1}}
B_{2}+FA_{1}\xrightarrow[]{f_{1}+1_{FA_{1}}}
(B_{3}+FA_{2})+FA_{1}\xrightarrow[]{<B_{3},FA_{2},FA_{1}>}
B_{3}+(FA_{2}+FA_{1})\\
&\ \ \ \ \ \ \  \xrightarrow[]{1_{B_{3}}+c_{FA_{2},FA_{1}}}
B_{3}+(FA_{1}+FA_{2})\xrightarrow[]{1_{B_{3}}+F_{+}^{-1}}
B_{3}+F(A_{1}+A_{2}).
\end{align*}
The above composition of morphisms in $CokerF$ is well-defined,
since, if $(f_{1},A_{1})$ and $(f_{1}^{'},A_{1}^{'})$ are
equivalent, i.e. $\exists\ \alpha:A_{1}\rightarrow A_{1}^{'}$, such
that $(1_{B_{2}}+F(\alpha))\circ f_{1}=f_{1}^{'}$. Then there exists
$\beta\triangleq\alpha+1_{A_{2}}:A_{1}+A_{2}\rightarrow
A_{1}^{'}+A_{2}$, such that the following diagram commutes:
\begin{center}
\scalebox{0.9}[0.85]{\includegraphics{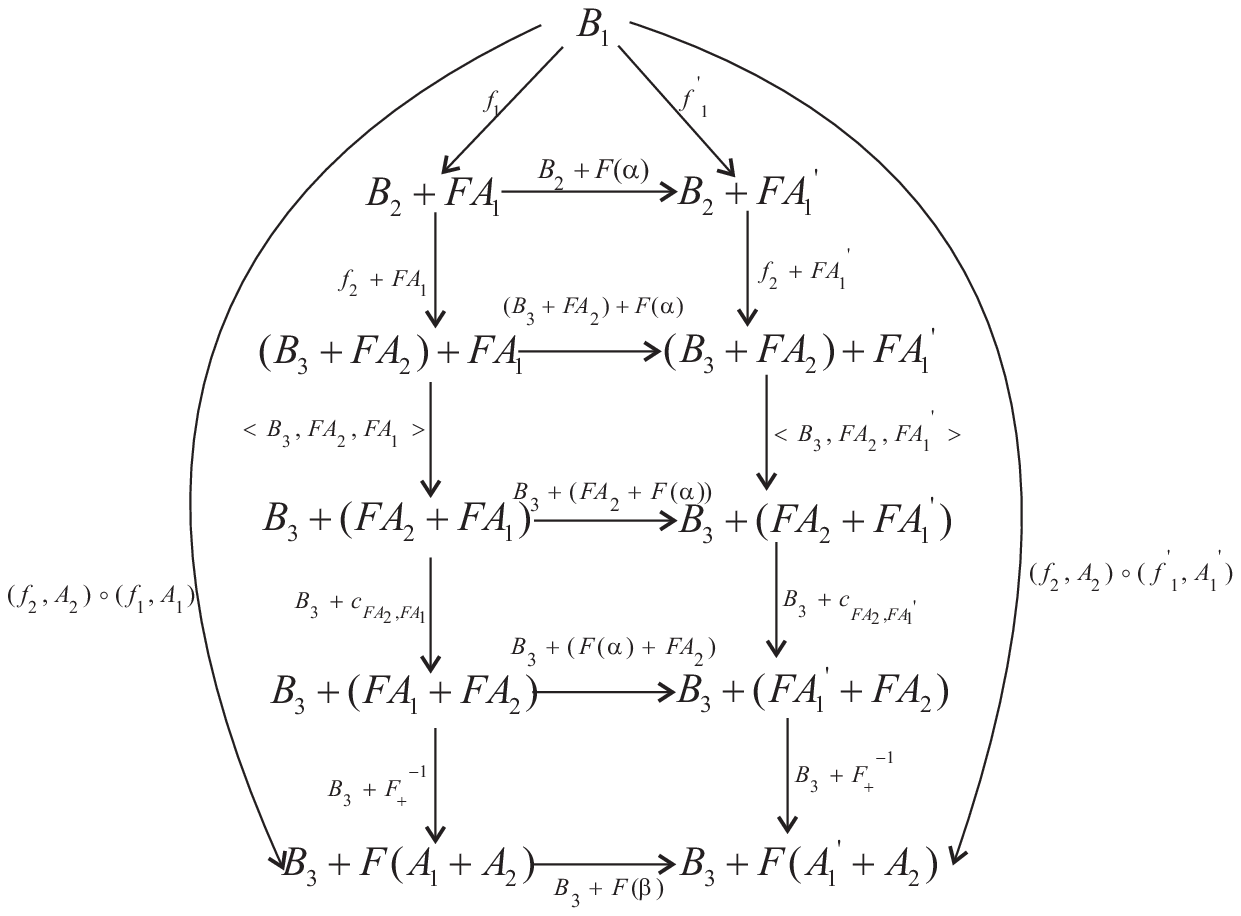}}
\end{center}
The above objects and morphisms, together with composition of
morphisms satisfy the following axioms:

(1) For any $B\in\cB$, there exists identity morphism
$[\widetilde{1_{B}},0]:B\rightarrow B$ in $CokerF$, where 0 is unit
object of $\cA$, $\widetilde{1_{B}}:B\xrightarrow[]{r_{B}^{-1}}
B+0\xrightarrow[]{1_{B}+F_{0}^{-1}}B+F0$, such that for any morphism
$[f,A]:B\rightarrow B_{1}$, $[f,A]\circ[\widetilde{1_{B}},0]=[f,A]$.

In fact, $[f,A]\circ[\widetilde{1_{B}},0]=[f^{'},A^{'}]$, where
$A^{'}=A+0$, $f^{'}$ is the composition
$B\xrightarrow[]{r_{B}^{-1}}B+0\xrightarrow[]{1_{B}+F_{0}^{-1}}B+F0\xrightarrow[]{f+1_{F0}}(B_{1}+FA)+F0
\xrightarrow[]{<B_{1},FA,F0>}B_{1}+(FA+F0)\xrightarrow[]{1_{B_{1}}+c_{FA,F0}}B_{1}+(F0+FA)
\xrightarrow[]{1_{B_{1}}+F_{+}^{-1}}B_{1}+F(0+A)$. There exists an
isomorphism $l_{A}:0+A\rightarrow A$, together with the properties
of F, such that $(f^{'},A^{'})$ is equivalent to $(f,A)$, so they
are same morphism in $CokerF$.

(ii)\ Associativity of compositions. Given morphisms
$$
B_{1}\xrightarrow[]{[f_{1},A_{1}]}B_{2}\xrightarrow[]{[f_{2},A_{2}]}B_{3}\xrightarrow[]{[f_{3},A_{3}]}B_{4}.
$$
There exists an isomorphism
$<A_{1},A_{2},A_{3}>:(A_{1}+A_{2})+A_{3}\rightarrow
A_{1}+(A_{2}+A_{2})$, such that
$$
[f_{3},A_{3}]\circ([f_{2},A_{2}]\circ[f_{1},A_{1}])=([f_{3},A_{3}]\circ[f_{2},A_{2}])\circ[f_{1},A_{1}]
$$
up to equality in $CokerF$.

\item $CokerF$ is a groupoid with zero object 0, where 0 is zero
object of $\cB$.

For any morphism $[f,A]:B_{1}\rightarrow B_{2}$ in $CokerF$, since
$f:B_{1}\rightarrow B_{2}+FA$ is a morphism in $\cB$, and $\cB$ is a
groupoid, there exists $f^{'}:B_{2}+FA\rightarrow B_{1}$, such that
$f^{'}\circ f=id$. For any $A\in obj(\cA)$, since $\cA$ is 2-group,
there exists $A^{*}\in obj(\cA)$, and a natural isomorphism
$\eta_{A}:A+A^{*}\rightarrow 0$.

Let $f^{*}$ be the composition
\begin{align*}
&B_{2}\rightarrow B_{2}+0\rightarrow B_{2}+F0\rightarrow
B_{2}+F(A+A^{*})\rightarrow\\
&B_{2}+(FA+F(A^{*})\rightarrow(B_{2}+FA)+F(A^{*})\rightarrow
B_{1}+F(A^{*}).
\end{align*}
From $\eta_{A}$, the properties of 2-groups and $F$, we get
$$
[(f^{*},A^{*})\circ(f,A)]=[1_{B_{1}},0].
$$

\item $CokerF$ is a monoidal category.

There is a bifunctor
\begin{align*}
&\hspace{1cm}+:CokerF\times CokerF\longrightarrow CokerF\\
&\hspace{3.5cm}(B_{1},B_{2})\mapsto B_{1}+B_{2},\\
&(B_{1}\xrightarrow[]{[f_{1},A_{1}]}B_{2},B_{3}\xrightarrow[]{[f_{2},A_{2}]}B_{4})\mapsto
B_{1}+B_{3}\xrightarrow[]{[f,A]}B_{2}+B_{4}
\end{align*}
where $B_{1}+B_{2}$ is the addition under the monoidal structure of
$\cB$, $A\triangleq A_{1}+A_{2}$, f is the following composition
\begin{align*}
&B_{1}+B_{3}\xrightarrow[]{f_{1}+f_{2}}(B_{2}+FA_{1})+(B_{4}+FA_{2})\xrightarrow[]{\langle^{B_{2}\
FA_{1}}_{B_{4}\
FA_{2}}\rangle}(B_{1}+B_{4})+(FA_{1}+FA_{2})\\
&\hspace{1.3cm}\xrightarrow[]{1+F_{+}^{-1}}(B_{1}+B_{4})+F(A_{1}+A_{2}).
\end{align*}
From the definition of $+$, the natural isomorphisms are those of
them in $\cB$, then the Mac Lane coherence conditions hold.

\item $CokerF$ is a symmetric 2-group.

For any object $B\in CokerF$, $B\in\cB$, since $\cB$ is a 2-group,
there exist $B^{*}\in\cB$, and $\eta_{B}:B^{*}+B\rightarrow 0$. Let
$\eta_{B}^{'}:B^{*}+B\rightarrow B^{*}+B+0\rightarrow B^{*}+B+F0$,
then $(\eta_{B}^{'},0)$ is an isomorphism in $CokerF$.\\
For any two objects $B_{1},B_{2}\in CokerF$, $B_{1},B_{2}\in\cB$,
$\cB$ is symmetric monoidal category, there exists
$c_{B_{1},B_{2}}:B_{1}+B_{2}\rightarrow B_{2}+B_{1}$, such that
$c_{B_{1},B_{2}}\circ c_{B_{2},B_{1}}=1_{B_{1}+B_{2}}$. Let
$c_{B_{1},B_{2}}^{'}$ be the composition
$$
B_{1}+B_{2}\xrightarrow[]{c_{B_{1},B_{2}}}B_{2}+B_{1}\xrightarrow[]{r_{B_{2}+B_{1}}^{-1}}
(B_{2}+B_{1})+0\xrightarrow[]{1_{B_{2}+B_{1}}+F_{0}}
(B_{2}+B_{1})+F0.
$$
Obviously, $c_{B_{1},B_{2}}^{'}\circ
c_{B_{2},B_{1}}^{'}=1_{B_{1}+B_{2}}$.

\item $CokerF$ is an $\cR$-2-module.

There is a bifunctor\\
\begin{align*}
&\hspace{1cm}\cdot:\cR\times CokerF\longrightarrow CokerF\\
&\hspace{2.7cm} (r,B)\mapsto r\cdot B,\\
&(r_{1}\xrightarrow[]{\varphi}r_{2},B_{1}\xrightarrow[]{[f,A]}B_{2}\mapsto
r_{1}\cdot B_{1}\xrightarrow[]{[f^{'},A^{'}]}r_{2}\cdot B_{2}
\end{align*}
where $r\cdot B\in\cB$ under the operation of $\cR$ on $\cB$,
$A^{'}\triangleq r_{2}\cdot A$, $f^{'}$ is the composition
$$
r_{1}\cdot B_{1}\rightarrow r_{2}\cdot B_{1}\rightarrow
r_{2}\cdot(B_{2}+FA)\rightarrow r_{2}\cdot B_{2}+r_{2}\cdot
FA\rightarrow r_{2}\cdot B_{2}+ F(r_{2}\cdot A).
$$
The natural isomorphisms are all induced by the the natural
isomorphisms in $\cB$, so they satisfy Fig.18--31. commute.

\item There is a functor:
\begin{align*}
&\hspace{0.4cm}p_{F}:\cB\longrightarrow CokerF\\
&\hspace{1.2cm}B\mapsto B,\\
&B_{1}\xrightarrow[]{f} B_{2}\mapsto
B_{1}\xrightarrow[]{[f^{'},0]}B_{2}
\end{align*}
where $f^{'}$ is the composition
$$
f^{'}:B_{1}\xrightarrow[]{f}B_{2}\xrightarrow[]{r_{B_{2}}}B_{2}+0\xrightarrow[]{1_{B_{2}}+F_{0}^{-1}}B_{2}+F0.
$$
$\mathbf{\cdot}$ $p_{F}$ is a functor.\\
Let $B_{1}\xrightarrow[]{f_{1}}B_{2}\xrightarrow[]{f_{2}}B_{3}$ be
morphisms in $\cB$, we need to prove $p_{F}(f_{2}\circ
f_{1})=p_{F}(f_{2})\circ p_{F}(f_{1})$. In fact,
$$
p_{F}(f_{1})=[f_{1}^{'},0]:B_{1}\rightarrow B_{2}\rightarrow
B_{2}+0\rightarrow B_{2}+F0,
$$
$$
p_{F}(f_{2})=[f_{2}^{'},0]:B_{2}\rightarrow B_{3}\rightarrow
B_{3}+0\rightarrow B_{3}+F0,
$$
$$
p_{F}(f_{1}\circ f_{2})=[f,0]:B_{1}\rightarrow B_{3}\rightarrow
B_{3}+0\rightarrow B_{3}+F0,
$$
$$
p_{F}(f_{2})\circ p_{F}(f_{1}):B_{1}\rightarrow
B_{2}+F0\rightarrow(B_{3}+F0)+F0\rightarrow B_{3}+(F0+F0)\rightarrow
B_{3}+F(0+0)\rightarrow B_{3}+F0.
$$
For any identity morphism $1_{B}:B\rightarrow B$ in $\cB$, we have
$p_{F}(1_{B})=(\widetilde{1_{B}},0):B\rightarrow B$.\\
$\mathbf{\cdot}$ There are natural isomorphisms:
\begin{align*}
&(p_{F})_{+}=id:p_{F}(B_{1}+B_{2})=B_{1}+B_{2}=p_{F}(B_{1})+p_{F}(B_{2}),\\
&(p_{F})_{0}=id:p_{F}(0)=0,\\
&(p_{F})_{2}=id:p_{F}(r\cdot B)=r\cdot B=r\cdot p_{F}(B).
\end{align*}
Obviously, $p_{F}$ is an $\cR$-homomorphism.

Define a 2-morphism
$$
\pi_{F}:p_{F}\circ F\Rightarrow 0
$$
given by, for any $A\in obj(\cA)$,
$$
(\pi_{F})_{A}=[(\pi_{F})_{A}^{'},A]:(p_{F}\circ
F)(A)=p_{F}(FA)=FA\rightarrow 0
$$
where $(\pi_{F})_{A}^{'}:FA\rightarrow FA+0\rightarrow FA+F0.$\\
Next, we will show $p_{F},\pi_{F}$ satisfy the universal property.
For any $\cK\in obj(\cR$-2-Mod), 1-morphism $G:\cB\rightarrow\cK$,
and $\varphi:G\circ F\Rightarrow 0$ in ($\cR$-2-Mod), define a
homomorphism
\begin{align*}
&G^{'}:CokerF\longrightarrow \cK\\
&\hspace{1.6cm}B\mapsto G(B),\\
&B_{1}\xrightarrow[]{[f,A]}B_{2}\mapsto
G(B_{1})\xrightarrow[]{G^{'}([f,A])}G(B_{2})
\end{align*}
Since $G$ is an $\cR$-homomorphism, so $G^{'}$ is. \\
Define a 2-morphism $\varphi^{'}:G^{'}\circ p_{F}\Rightarrow G $, by
$$
\varphi^{'}_{B}\triangleq 1_{GB}:(G^{'}\circ
p_{F})(B)=G^{'}(B)=G(B),
$$
for any $B\in obj(\cB)$, such that Fig.40. commutes.

If $G^{''}$ and $\varphi^{''}$ satisfy the same conditions as
$G^{'}$ and $\varphi^{'}$, there is a 2-morphism
$\psi:G^{''}\Rightarrow G^{'}$ defined by, $\psi_{B}\triangleq
\varphi^{''}_{B}:G^{''}(B)\rightarrow G^{'}(B)$,for each $B\in
obj(\cB)$, such that the diagram Fig.41. commutes .

Obviously, $\psi^{''}$ is the unique 2-morphism for $\varphi^{''}$
is.
\end{itemize}
\end{proof}

A groupoid enriched category (for short, a $Gpd$-category) is in
fact a 2-category, satisfies some special properties, which plays an
important role in 2-abelian category(more details see \cite{18}).
Next we will give some results about ($\cR$-2-Mod).

\begin{Lemma}The 2-category ($\cR$-2-Mod) is a $Gpd$-category.
\end{Lemma}
\begin{proof}
($\cR$-2-Mod) contains the following ingredients:

$\mathbf{(1).}$ For any $\cA,\cB\in obj$($\cR$-2-Mod),
$Hom(\cA,\cB)$ is a groupoid, with $\cR$-homomorphisms from $\cA$ to
$\cB$ as its objects and  the morphisms of two $\cR$-homomorphisms
as its morphisms.

Composition of morphisms. Let $\tau:F\Rightarrow G,\
\sigma:G\Rightarrow H$ be morphisms in $Hom(\cA,\cB)$.
$\sigma\circ\tau:F\Rightarrow H$ is given by
$$
(\sigma\circ\tau)_{A}\triangleq\sigma_{A}\circ\tau_{A}:FA\xrightarrow[]{\tau_{A}}GA\xrightarrow[]{\sigma_{A}}HA.
$$
It is easy to check $\sigma\circ\tau$ is a morphism of
$\cR$-homomorphisms from $F$ to $H$(see Theorem 1).

The above objects and morphisms satisfy the following axioms:

(i) For any $F\in Hom(\cA,\cB)$, $\exists\ 1_{F}:F\Rightarrow F$,
defined by $(1_{F})_{A}\triangleq 1_{FA},\ \forall\ A\in\cA$, such
that for any $\tau:F\Rightarrow G$, and $\sigma:H\Rightarrow F$, we
have $\tau\circ1_{F}=\tau,\ 1_{F}\circ\sigma=\sigma$, since
$(\tau\circ1_{F})(A)=\tau_{A}\circ(1_{F}){A}=\tau_{A}\circ1_{FA}=\tau_{A}$,
$(1_{F}\circ\sigma)(A)=(1_{F})_{A}\circ\sigma_{A}=1_{FA}\circ\sigma_{A}=\sigma_{A}$,
for any $A\in\cA$.

(ii) Associativity of the composition. Given morphisms
$$
F_{1}\xrightarrow[]{\tau_{1}}F_{2}\xrightarrow[]{\tau_{2}}F_{3}\xrightarrow[]{\tau_{3}}F_{4}
$$
in $Hom(\cA,\cB)$. Then
$(\tau_{3}\circ\tau_{2})\circ\tau_{1}=\tau_{3}\circ(\tau_{2}\circ\tau_{1})$,
since
$$
((\tau_{3}\circ\tau_{2})\circ\tau_{1})_{A}=((\tau_{3})_{A}\circ(\tau_{2})_{A})\circ(\tau_{1})_{A}
=(\tau_{3})_{A}\circ((\tau_{2})_{A}\circ(\tau_{1})_{A})=(\tau_{3}\circ(\tau_{2}\circ\tau_{1}))_{A},
$$
for $\forall\ A\in\ obj(\cA)$.

(iii) For any morphism $\tau:F\Rightarrow G:\cA\rightarrow \cB$,
$\exists\ \tau^{*}:G\Rightarrow F$, such that $\tau^{*}\circ
\tau\backsimeq 1_{F}$. In fact, $\forall\ A\in\cA,\
\tau_{A}:FA\rightarrow GA$ is a morphism in $\cB$, and $\cB$ is a
groupoid, so $\exists\ (\tau_{A})^{*}:GA\rightarrow FA$, such that $
(\tau_{A})^{*}\circ\tau_{A}\backsimeq 1_{FA}$. Define
$\tau^{*}:G\Rightarrow F$ by $(\tau^{*})_{A}\triangleq
(\tau_{A})^{*}$, such that $(\tau^{*}\circ
\tau)_{A}=(\tau^{*})_{A}\circ\tau_{A}\backsimeq 1_{FA}=(1_{F})_{A}$.

$\mathbf{(2).}$ For any $\cA\in obj(\cR$-2-Mod), there is an
$\cR$-homomorphism $1_{\cA}:\cA\rightarrow \cA$, defined by
$1_{\cA}(A)=A,\ \forall A\in\cA$.

$\mathbf{(3).}$ For any $\cA,\ \cB,\ \cC\in obj(\cR$-2-Mod), there
is a functor composition
\begin{align*}
&comp:Hom(\cA,\cB)\times Hom(\cB,\cC)\longrightarrow Hom(\cA,\cC)\\
&\hspace{4.3cm}(F\ ,\ G)\mapsto comp(F,G)\triangleq G\circ F,\\
&\hspace{1.2cm}(\alpha:F\Rightarrow F^{'},\beta:G\Rightarrow
G^{'})\mapsto comp(\alpha,\beta)\triangleq\beta\ast\alpha
\end{align*}
where $\beta\ast\alpha$ is the horizontal composition in Theorem 1.

From Theorem 1, we also have

$\mathbf{(4)}.$ For all $\cA,\ \cB,\ \cC,\ \cD\in obj(\cR$-2-Mod),
and $F\in Hom(\cA,\cB),\ G\in Hom(\cB,\cC),\ H\in Hom(\cC,\cD)$,
there is a natural transformation
$$
\alpha_{H,G,F}=id:(H\circ G)\circ F\Rightarrow H\circ(G\circ F).
$$

$\mathbf{(5)}.$ For all $\cA,\ \cB\in obj(\cR$-2-Mod), $F\in
Hom(\cA,\cB)$, there are natural transformations
$$
\rho_{F}:F\circ1_{\cA}\Rightarrow F,
$$
$$
\lambda_{F}:1_{\cB}\circ F\Rightarrow F
$$
given by, for any $A\in obj(\cA),\ f:A_{1}\rightarrow A_{2}$ in
$\cA$, $(\rho_{F})_{A}=id,\ (\lambda_{F})_{A}=id,\
(\rho_{F})_{f}=F(f),\ (\lambda_{F})_{f}=F(f)$. Obviously,
$\rho_{F},\ \lambda_{F}$ are morphisms in $Hom(\cA,\cB)$.

Given $\cR$-2-module homomorphisms
$\cA\xrightarrow[]{F}\cB\xrightarrow[]{G}\cC\xrightarrow[]{H}\cD\xrightarrow[]{K}\cE$.
Since $\alpha_{-,-,-},\rho_{-},\lambda_{-}$ are identities, so they
satisfy the following diagrams commute:
\begin{center}
\scalebox{0.9}[0.85]{\includegraphics{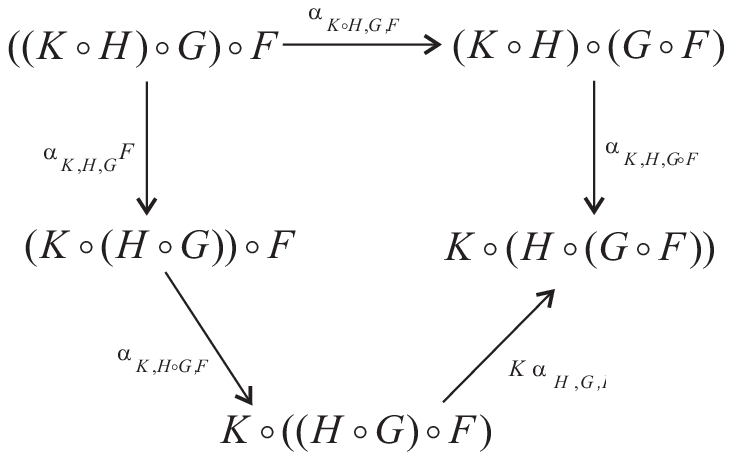}}
\scalebox{0.9}[0.85]{\includegraphics{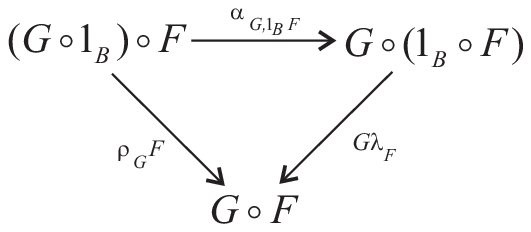}}
\end{center}
\end{proof}

\begin{cor}
($\cR$-2-Mod) is a $Gpd^{*}$-category, where $Gpd^*$ means the
pointed groupoid appearing in \cite{18}.
\end{cor}
\begin{proof}
For any $\cA,\ \cB\in obj(\cR$-2-Mod), $Hom(\cA,\cB)$ is a pointed
groupiod where the point is just the zero $\cR$-homomorphism from
$\cA$ to $\cB$. Using the similar methods in the proof of Lemma 1,
it is easy to prove it.
\end{proof}

\begin{Lemma}
For any $\cA,\ \cB\in obj(\cR$-2-Mod), $Hom(\cA,\cB)$ is a symmetric
2-group.
\end{Lemma}
\begin{proof}
From Lemma 1, $Hom(\cA,\cB)$ is a groupoid, so we need to give the
monoidal structure on it, and prove it is symmetric under this
monoidal structure.

There is a bifunctor:
\begin{align*}
&+:Hom(\cA,\cB)\times Hom(\cA,\cB)\longrightarrow Hom(\cA,\cB)\\
&\hspace{4.2cm}(F,G)\mapsto F+G,\\
&\hspace{0.8cm}(\tau:F\Rightarrow F^{'},\sigma:G\Rightarrow
G^{'})\mapsto \tau+\sigma:F+G\Rightarrow F^{'}+G^{'}
\end{align*}
given by, $(F+G)(A)\triangleq FA+GA$ , $(\tau+\sigma)_{A}\triangleq
\tau_{A}+\sigma_{A}$ under the monoidal addition in
$\cB$, for any $A\in\cA$. \\
The unit object $0\in Hom(\cA,\cB)$ is just the zero
$\cR$-homomorphism of $\cR$-2-modules.

Moreover, there are natural isomorphisms:
\begin{align*}
&<F,G,H>:(F+G)+H\Rightarrow F+(G+H),\\
&\hspace{1.9cm}l_{F}:0+F\Rightarrow F,\\
&\hspace{1.9cm}r_{F}:F+0\Rightarrow F
\end{align*}
defined by, $<F,G,H>_{A}\triangleq <FA,GA,HA>,\ (l_{F})_{A}=l_{FA},\
(r_{F})_{A}=r_{FA}$, $\forall\ A\in\cA$. Since $\cB$ is a monoidal
category, the Mac Lane coherence conditions hold, i.e. Fig.1-2.
commute.

For any $F,G\in obj(Hom(\cA,\cB))$, $(c_{F,G})_{A}\triangleq
c_{FA,GA},\exists\  c_{F,G}:F+G\Rightarrow G+F$, such that
$c_{FA,GA}\circ c_{GA,FA}\backsimeq 1_{GA+FA}$. Then $c_{F,G}\circ
c_{G,F}= 1_{G+F}$.

For any $F\in obj(Hom(\cA,\cB))$. Define $F^{*}:\cB\rightarrow\cA$
by $F^{*}(A)=(FA)^{*},\ \forall A\in\cA$, where $(FA)^{*}$ is the
inverse of FA in $\cB$, with natural isomorphism
$\eta_{FA}:(FA)^{*}+FA\rightarrow 0$. So there is a natural
isomorphism $\eta_{F}:F^{*}+F\Rightarrow 0$, given by
$(\eta_{F})_{A}\triangleq \eta_{FA}$.
\end{proof}

\begin{Lemma}
($\cR$-2-Mod) is a presemiadditive $Gpd$-category refer to the
Definition 218 in \cite {18}.
\end{Lemma}
\begin{proof}
For any $\cA,\ \cB\in obj(\cR$-2-Mod), $Hom(\cA,\cB)$ is a symmetric
monoid groupoid (Lemma 2) with transformations natural in each
variables:
\begin{align*}
&\varphi_{G_{1},G_{2}}^{H}:H\circ(G_{1}+G_{2})\Rightarrow H\circ
G_{1}+H\circ G_{2},\\
&\psi_{G}^{H_{1},H_{2}}:(H_{1}+H_{2})\circ G\Rightarrow H_{1}\circ
G+H_{2}\circ G,\\
&\varphi_{0}^{H}:0\Rightarrow H\circ 0,\\
&\psi_{G}^{0}:0\Rightarrow 0\circ G,
\end{align*}
defined by, $\forall A\in\cA$,
\begin{align*}
&(\varphi_{G_{1},G_{2}}^{H})_{A}\triangleq
(H_{+})_{G_{1}A,G_{2}A}:(H\circ(G_{1}+G_{2}))(A)=H(G_{1}A+G_{2}A)\rightarrow
H(G_{1}A)+H(G_{2}A),\\
&(\psi_{G}^{H_{1},H_{2}})_{A}\triangleq id:((H_{1}+H_{2})\circ
G)(A)=H_{1}(GA)+H_{2}(GA),\\
&(\varphi_{0}^{H})_{A}\triangleq H_{0}:0(A)=0\rightarrow
H(0A)=H(0),\\
&(\psi_{G}^{0})_{A}\triangleq id:0(A)=0\rightarrow 0(GA)=0.
\end{align*}
The above natural transformations satisfy the following conditions:

\noindent\textbf{(a)}\ For $\cA\in obj(\cR$-2-Mod), $H\in
Hom(\cB,\cC)$, the functor $H\circ-:Hom(\cA,\cB)\rightarrow
Hom(\cA,\cC)$, with $\varphi_{G_{1},G_{2}}^{H}$ and
$\varphi_{0}^{H}$, is symmetric monoidal.

In fact, $(G_{1}+G_{2})(A)\triangleq G_{1}A+G_{2}A,\ \forall\
A\in\cA$. Moreover, F satisfies the following commutative diagrams:
\begin{center}
\scalebox{0.9}[0.85]{\includegraphics{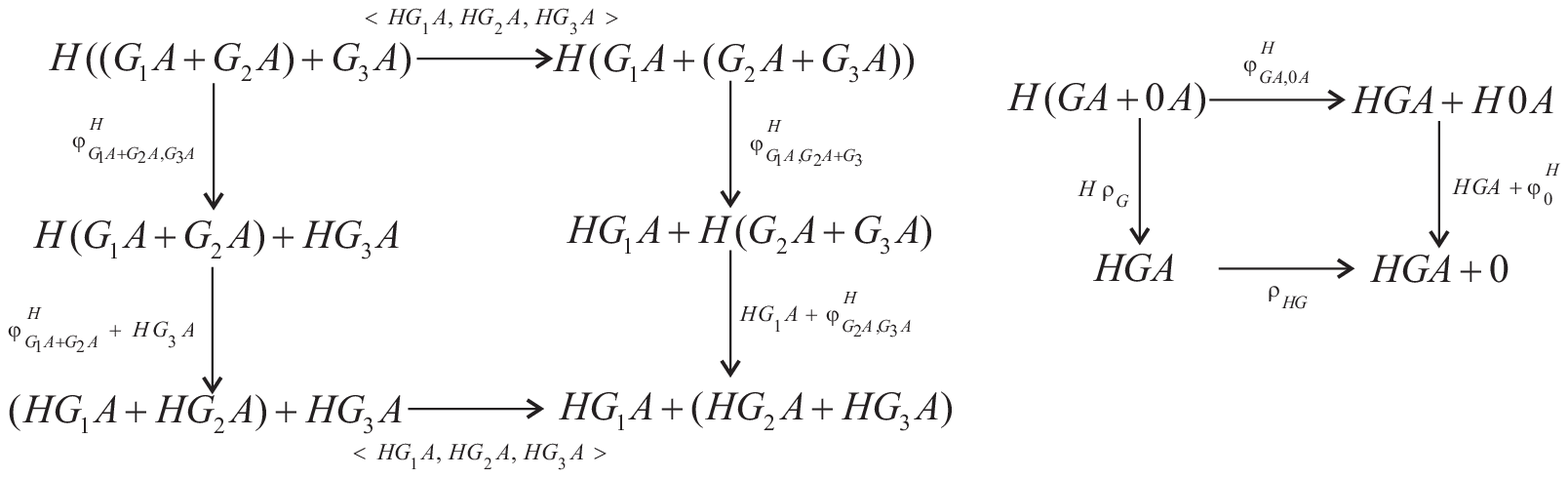}}
\scalebox{0.9}[0.85]{\includegraphics{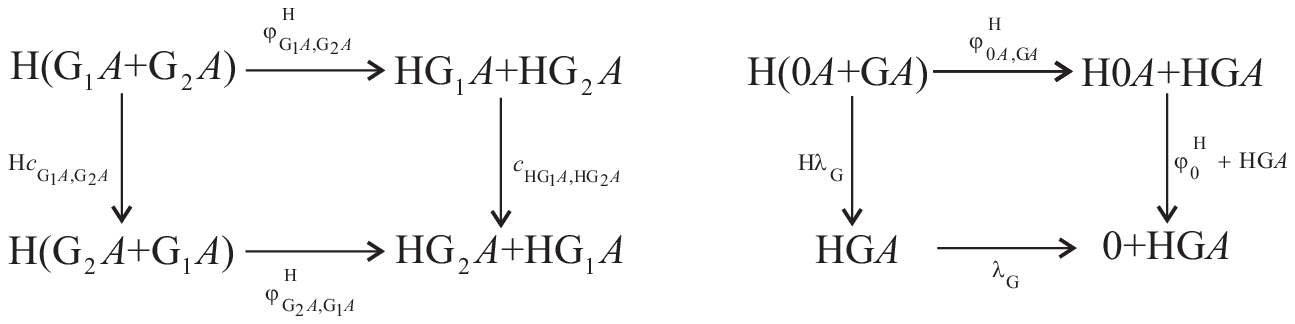}}
\end{center}
So $H\circ-$ is a symmetric monoidal functor.

\noindent\textbf{(b)} For $G\in Hom(\cA,\cB)$, $\cC\in
obj(\cR$-2-Mod), the functor $-\circ G:Hom(\cB,\cC)\rightarrow
Hom(\cA,\cC)$, with $\psi_{G}^{H_{1},H_{2}}$ and $\psi_{G}^{0}$
symmetric monoidal in the same methods as in \noindent\textbf{(a)}.

\noindent\textbf{(c)} For all $H_{1},H_{2}\in Hom(\cB,\cC)$, the
transformations $\psi^{H_{1},H_{2}}_{-}$ and $\psi^{0}_{-}$ are
monoidal functors, Since $H_{1},H_{2}$ are $\cR$-homomorphisms.

\noindent\textbf{(d)} For $G:\cB\rightarrow \cC,\
H:\cC\rightarrow\cD$, and $F,F^{'}\in Hom(\cA,\cB)$,
$\alpha_{H,G,-}$ is a monoidal natural identity, i.e. the following
diagrams commute:
\begin{center}
\scalebox{0.9}[0.85]{\includegraphics{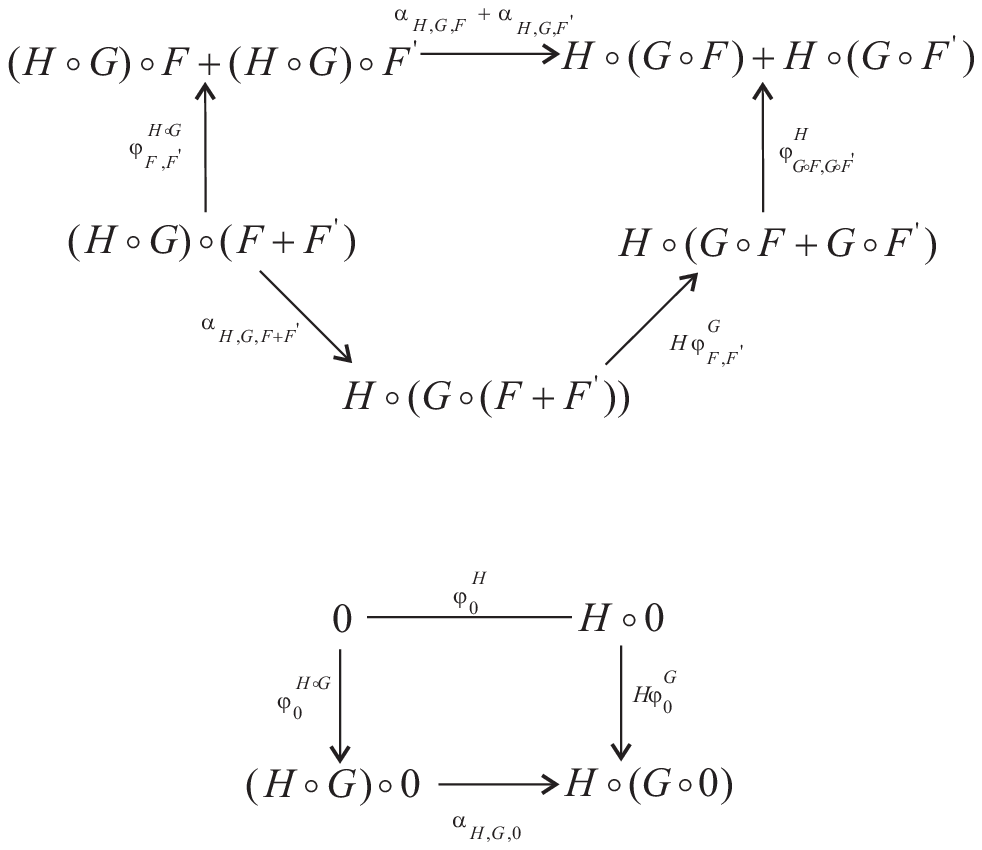}}
\end{center}
Similarly, the following conditions hold.

\noindent\textbf{(d)} For $F\in Hom(\cA,\cB),\ H\in Hom(\cC,\cD)$,
$\alpha_{H,-,F}$ is a monoidal natural transformation.

\noindent\textbf{(e)} For $F\in Hom(\cA,\cB),\ G\in Hom(\cB,\cC)$,
$\alpha_{-,G,F}$ is a monoidal natural transformation.

\noindent\textbf{(f)}  For all $\cA,\ \cB\in obj(\cR$-2-Mod), since
$\lambda_{-}$ and $\rho_{-}$ are identities, so the following unit
natural transformations are monoidal:
\begin{center}
\scalebox{0.9}[0.85]{\includegraphics{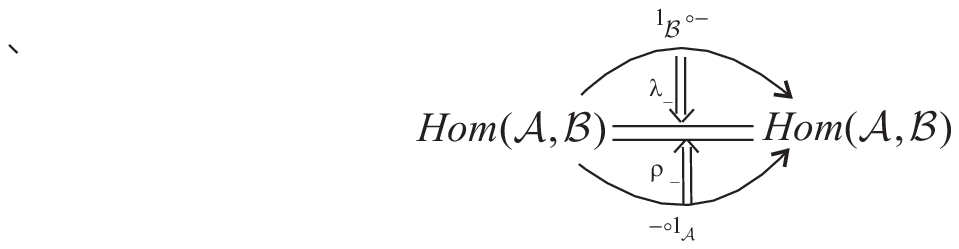}}
\end{center}
i.e. the following diagrams commute:
\begin{center}
\scalebox{0.9}[0.85]{\includegraphics{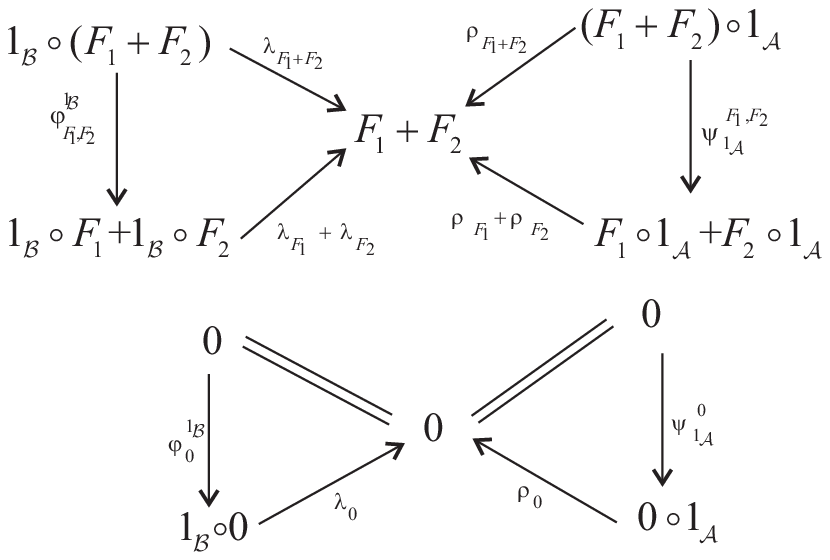}}
\end{center}
\end{proof}

\begin{Prop}
($\cR$-2-Mod) has all finite biproducts.
\end{Prop}
\begin{proof}
For any two objects $\cA,\cB$ in ($\cR$-2-Mod), there is a new
object $\cA\times\cB$ in ($\cR$-2-Mod).
\begin{itemize}
\item $\cA\times\cB$ is a category consisting of the following data:

$\cdot$\ Objects are pairs $(A,B)$, where $A\in obj(\cA),\ B\in
obj(\cB)$.

$\cdot$\ A morphism between $(A_{1},B_{1})$ and $(A_{2},B_{2})$ is a
pair $(f,g)$, where $f:A_{1}\rightarrow A_{2}$ is a morphism in
$\cA$, $g:B_{1}\rightarrow B_{2}$ is a morphism in $\cB$.

$\cdot$\ Composition of morphisms. Given morphisms
$(A_{1},B_{1})\xrightarrow[]{(f_{1},g_{1})}(A_{2},B_{2})\xrightarrow[]{(f_{2},g_{2})}(A_{3},B_{3})$,
$(f_{2},g_{2})\circ(f_{1},g_{1})\triangleq(f_{2}\circ
f_{1},g_{2}\circ g_{1})$ is also a morphism in $\cA\times\cB$.

The above ingredients satisfy the following axioms:

(1) For any $(A,B)\in obj(\cA\times\cB)$, there exists a morphism
$1_{(A,B)}\triangleq(1_{A},1_{B}):(A,B)\rightarrow (A,B)$, such that
for any morphism $(f,g):(A,B)\rightarrow (A_{1},B_{1})$,
$(f,g)\circ1_{(A,B)}=(f,g)$.

(2) Associativity of compositions. Given morphisms
$$
(A_{1},B_{1})\xrightarrow[]{(f_{1},g_{1})}(A_{2},B_{2})\xrightarrow[]{(f_{2},g_{2})}(A_{3},B_{3})
\xrightarrow[]{(f_{3},g_{3})}(A_{4},B_{4}),
$$
for $(f_{3}\circ f_{2})\circ f_{1}=f_{3}\circ(f_{2}\circ f_{1})$ and
$(g_{3}\circ g_{2})\circ g_{1}=g_{3}\circ(g_{2}\circ g_{1})$, we
have
$$
((f_{3},g_{3})\circ (f_{2},g_{2}))\circ
(f_{1},g_{1})=(f_{3},g_{3})\circ((f_{2},g_{2})\circ (f_{1},g_{1})).
$$

\item $\cA\times\cB$ is a symmetric 2-group.

For any morphism $(f,g):(A_{1},B_{1})\rightarrow (A_{2},B_{2})$ in
$\cA\times\cB$, and $\cA,\ \cB$ are groupoids, there exist
$f^{*}:A_{2}\rightarrow A_{1},\ g^{*}:B_{2}\rightarrow B_{1}$, such
that $f^{*}\circ f=1_{A_{1}},\ g^{*}\circ g=1_{B_{1}}$. So there
exists $(f,g)^{*}\triangleq(f^{*},g^{*}):(A_{2},B_{2})\rightarrow
(A_{1},B_{1})$, such that $(f,g)^{*}\circ (f,g)=1_{(A_{1},B_{1})}$.

There is an unit object $0=(0_{\cA},0_{\cB})$ in $\cA\times\cB$,
where $0_{\cA}$ and $0_{\cB}$ are unit objects of $\cA$ and $\cB$,
respectively.

There are a bifunctor
\begin{align*}
&+:(\cA\times\cB)\times
(\cA\times\cB)\longrightarrow(\cA\times\cB)\\
&\hspace{0.6cm}((A_{1},B_{1}),(A_{2},B_{2}))\mapsto
(A_{1},B_{1})+(A_{2},B_{2})\triangleq(A_{1}+A_{2}),B_{1}+B_{2}),\\
&\hspace{1cm}((f_{1},g_{1}),(f_{2},g_{2}))
\mapsto(f_{1},g_{1})+(f_{2}+g_{2})\triangleq(f_{1}+f_{2},g_{1}+g_{2})
\end{align*}
and natural isomorphisms:
\begin{align*}
&<(A_{1},B_{1}),(A_{2},B_{2}),(A_{3},B_{3})>:((A_{1},B_{1})+(A_{2},B_{2}))+(A_{3},B_{3})\\
&\hspace{5cm}\rightarrow
(A_{1},B_{1})+((A_{2},B_{2})+(A_{3},B_{3})),\\
&\hspace{4cm}l_{(A,B)}:0+(A,B)\rightarrow (A,B),\\
&\hspace{4cm}r_{(A,B)}:(A,B)+0\rightarrow (A,B)
\end{align*}
given by $<(A_{1},B_{1}),(A_{2},B_{2}),(A_{3},B_{3})>\triangleq
 (<A_{1},A_{},A_{3}>,<B_{1},B_{2},B_{3}>),$
 $l_{(A,B)}\triangleq (l_{A},l_{B}),\
 r_{(A,B)}\triangleq (r_{A},r_{B})$. Obviously, the Mac Lane coherence conditions hold.

For any $(A,B)\in obj(\cA\times\cB)$, $A\in\cA,\ B\in\cB$, there
exist $A^{*}\in\cA, B^*\in\cB$, and natural isomorphisms
$\eta_{A}:A^{*}+A\rightarrow 0,\ \eta_{B}:B^{*}+B\rightarrow 0.$
Then there exists $(A,B)^{*}\triangleq(A^{*},B^{*})$, and natural
isomorphism $\eta_{(A,B)}\triangleq
(\eta_{A},\eta_{B}):(A,B)^{*}+(A,B)\rightarrow 0.$

For any two objects $(A_{1},B_{1}),(A_{2},B_{2})\in(\cA\times\cB)$,
since there are natural isomorphisms
$c_{A_{1},A_{2}}:A_{1}+A_{2}\rightarrow A_{2}+A_{1},\
c_{B_{1},B_{2}}:B_{1}+B_{2}\rightarrow B_{2}+B_{1}$, with
$c_{A_{1},A_{2}}\circ c_{A_{2},A_{1}}=id,\ c_{B_{1},B_{2}}\circ
c_{B_{2},B_{1}}=id$. Then we get a natural isomorphism
$c_{(A_{1},B_{1}),(A_{2},B_{2})}\triangleq
(c_{A_{1},A_{2}},c_{B_{1},B_{2}}):(A_{1},B_{1})+(A_{2},B_{2})\rightarrow(A_{2},B_{2})+(A_{1},B_{1})
$, with $c_{(A_{1},B_{1}),(A_{2},B_{2})}\circ
c_{(A_{2},B_{2}),(A_{1},B_{1})}=id.$

\item $\cA\times\cB$ is an $\cR$-2-module.

There is a bifunctor
\begin{align*}
&\hspace{2.8cm}\cdot:\cR\times(\cA\times\cB)\longrightarrow(\cA\times\cB)\\
&\hspace{3.8cm}(r,(A,B))\mapsto r\cdot(A,B)\triangleq(r\cdot A,r\cdot B),\\
&(r_{1}\xrightarrow[]{\varphi}r_{2},(A_{1},B_{1}))\xrightarrow[]{(f,g)}(A_{2},B_{2})\mapsto
(r_{1}\cdot A_{1},r_{2}\cdot A_{2})\xrightarrow[]{(\varphi\cdot
f,\varphi\cdot g)}(r_{1}\cdot B_{1},r_{2}\cdot B_{2})
\end{align*}
since $\cA,\ \cB$ are $\cR$-2-modules .

Also, there are natural isomorphisms:
\begin{align*}
&a^{r}_{(A_{1},B_{1}),(A_{2},B_{2})}\triangleq(a^{r}_{A_{1},A_{2}},a^{r}_{B_{1},B_{2}})
:r\cdot((A_{1},B_{1})+(A_{2},B_{2}))\rightarrow
(r\cdot(A_{1}+A_{2}),r\cdot(B_{1}+B_{2})),\\
&b^{r_{1},r_{2}}_{(A,B)}\triangleq(b^{r_{1},r_{2}}_{A},b^{r_{1},r_{2}}_{B}):
(r_{1}+r_{2})\cdot(A,B)\rightarrow
r_{1}\cdot(A,B)+r_{2}\cdot(A,B),\\
&b_{r_{1},r_{2},(A,B)}\triangleq
(b_{r_{1},r_{2},A},b_{r_{1},r_{2},B}):(r_{1}r_{2})\cdot(A,B)\rightarrow
r_{1}\cdot(r_{2}\cdot(A,B)),\\
&i_{(A,B)}\triangleq(i_{A},i_{B}):I\cdot(A,B)\rightarrow (I\cdot
A,I\cdot B),\\
&z_{r}\triangleq (z_{r},z_{r}):r\cdot 0=r\cdot(0,0)\rightarrow 0.
\end{align*}
Since $\cA,\ \cB$ are $\cR$-2-modules, their natural isomorphisms
make Fig.18.-Fig.31. commute, so do
$a^{r}_{(A_{1},B_{1}),(A_{2},B_{2})},\ b^{r_{2},r_{1}}_{(A,B)},\
b_{r_{1},r_{2},(A,B)},\ i_{(A,B)}$ and $z_{r}$.

\item $\cA\times\cB$ is the biproduct of $\cA$ and $\cB$.

We need to prove $\cA\times\cB$ is not only the product but also the
coproduct of $\cA$ and $\cB$.

There are $\cR$-homomorphisms:
$$
\cA\xleftarrow[]{p_{1}}\cA\times\cB\xrightarrow[]{p_{2}}\cB
$$
$$
A\leftarrow(A,B)\rightarrow B,
$$
$$
A_{1}\xleftarrow[]{f}A_{2}\leftarrow((A_{1},B_{1})\xrightarrow[]{(f,g)}(A_{2},B_{2}))\rightarrow
B_{1}\xrightarrow[]{g} B_{2}
$$
and the faithful $\cR$-homomorphisms:
$$
\cA\xrightarrow[]{i_{1}}\cA\times\cB\xleftarrow[]{i_{2}}\cB
$$
$$
A\rightarrow(A,0),\ \ \ \ \ \
$$
$$
\ \ \ \ \ \ \ \ \ \ \ (0,B)\leftarrow B
$$
Next, we will show that ($\cA\times\cB,i_{1},i_{2})$ satisfies the
universal property of coproduct\cite{5,18}.\\
For any $\cR$-2-module $\cK$ and $\cR$-homomorphisms
$F_{1}:\cA\rightarrow\cK,\ F_{2}:\cB\rightarrow\cK$, there is an
$\cR$-homomorphism $G:\cA\times\cB\rightarrow\cK$, with
$$G(A,B)=F_{1}A+F_{2}B,$$
and isomorphisms $l_{1}:G\circ i_{1}\Rightarrow F_{1},\ l_{2}:G\circ
i_{2}\Rightarrow F_{2}$, given by
$$(l_{1})_{A}:(G\circ
i_{1})(A)=G(A,0)=F_{1}A+F_{2}0\rightarrow F_{1}A+0\rightarrow
F_{1}A,$$
$$(l_{2})_{B}:(G\circ
i_{2})(B)=G(0,B)=F_{1}0+F_{2}B\rightarrow 0+F_{2}B\rightarrow
F_{2}B.
$$
So $\cA\times\cB$ is a coproduct.

Using  the Proposition\ 225 in \cite{18},
($\cA\times\cB,p_{1},p_{2})$ is a product, such that
$$
p_{2}\circ i_{1}=0,\ p_{1}\circ i_{1}=1_{\cA},\ p_{1}\circ i_{2}=0,\
p_{2}\circ i_{2}=1_{\cB}.
$$
So, $\cA\times\cB$ is a biproduct.
\end{itemize}
\end{proof}
\begin{Def}(\cite{18} Definition\ 242)
Let $\cC$ be a $Gpd$-category.

1. We say that $\cC$ is semiadditive if it is presemiadditive and
has all finite biproducts.

2. We say that $\cC$ is additive if it is preadditive and has all
finite biproducts.
\end{Def}

\begin{cor}
($\cR$-2-Mod) is an additive $Gpd$-category.
\end{cor}

Next, we use the definitions of pips(copips) and roots(coroots)
given in \cite{8,18} to give the next definitions in ($\cR$-2-Mod).
\begin{Def}
Let $F:\cA\rightarrow\cB$ be $\cR$-homomorphism.\\
\begin{itemize}
\item The pip of F is given by an $\cR$-2-module $PipF$, two zero
$\cR$-2-module homomorphisms $0:PipF\rightarrow \cA$, and morphism
$\sigma:0\Rightarrow 0$ of $\cR$-homomorphisms as in the following
diagram:
\begin{center}
\scalebox{0.9}[0.85]{\includegraphics{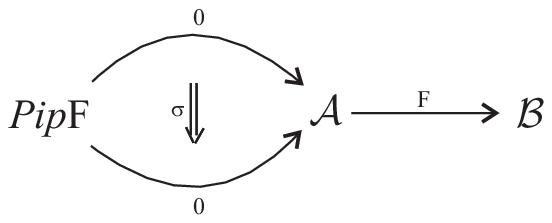}}
\end{center}
such that $F\ast\sigma=1_{0}:0\Rightarrow 0:PipF\rightarrow \cB$,
and for any other $\cD$ as in the following diagram
\begin{center}
\scalebox{0.9}[0.85]{\includegraphics{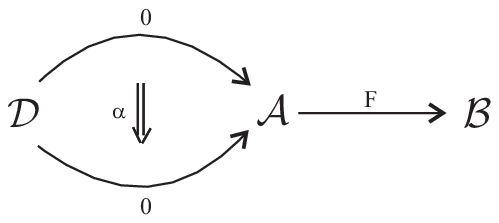}}
\end{center}
with $F\ast\alpha=1_{0}$. There is an $\cR$-homomorphism
$G:\cD\rightarrow PipF$, such that $\sigma\ast G=\alpha$. $G$ is
unique up to an invertible morphism of $\cR$-homomorphisms, i.e. if
there is a $G^{'}:\cD\rightarrow PipF$, with $\sigma\ast
G^{'}=\alpha$, there exists a unique isomorphism
$\tau:G^{'}\Rightarrow G$, such that
\begin{center}
\scalebox{0.9}[0.85]{\includegraphics{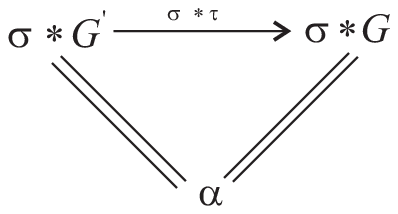}}
\end{center}
commutes.

\item The copip of F is given by an $\cR$-2-module $CopipF$, two zero
$\cR$-homomorphisms $0:\cB\rightarrow CopipF$, and morphism of
$\cR$-homomorphisms $\sigma:0\Rightarrow 0:\cB\rightarrow CopipF$,
as in in the following diagram:
\begin{center}
\scalebox{0.9}[0.85]{\includegraphics{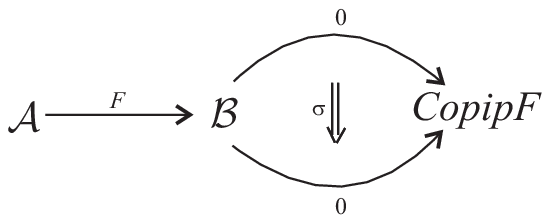}}
\end{center}
such that $\sigma\ast F=1_{0}:0\Rightarrow 0:\cA\rightarrow CopipF$,
and for any other
\begin{center}
\scalebox{0.9}[0.85]{\includegraphics{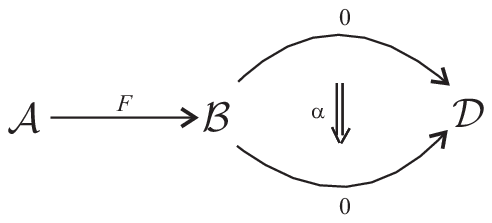}}
\end{center}
with $\alpha\ast F=1_{0}$. There is an $\cR$-homomorphism
$G:CopipF\rightarrow\cD$, such that $G\ast\sigma=\alpha$. $G$ is
unique up to an invertible morphism of $\cR$-homomorphisms, i.e. if
there is a $G^{'}:CopipF\rightarrow \cD$, with $
G^{'}\ast\sigma=\alpha$, there exists a unique isomorphism
$\tau:G^{'}\Rightarrow G$.
\end{itemize}
\end{Def}

\begin{Def}
Let $\alpha:0\Rightarrow 0:\cA\rightarrow\cB$ be 2-morphism in
($\cR$-2-Mod).
\begin{itemize}
\item The root of $\alpha$ is an $\cR$-2-module $Root \alpha$ and an
$\cR$-homomorphism $F:Root\alpha\rightarrow\cA$ as in the following
diagram
\begin{center}
\scalebox{0.9}[0.85]{\includegraphics{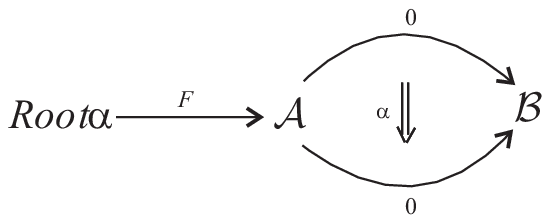}}
\end{center}
such that $\alpha\ast F=1_{0}:0\Rightarrow 0:\cD\rightarrow\cB$, for
any other as in
\begin{center}
\scalebox{0.9}[0.85]{\includegraphics{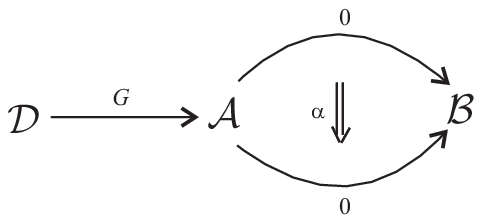}}
\end{center}
with $\alpha\ast G=1_{0}$, there exist an $\cR$-homomorphism
$G^{'}:\cD\rightarrow Root \alpha$ and an invertible morphism
$\varphi:F\circ G^{'}\Rightarrow G$. The pair $(G^{'},\varphi)$ is
unique up to an invertible morphism, i.e. If $(G^{''},\varphi^{'})$
satisfies the same conditions as $(G^{'},\varphi)$, there exists a
unique $\tau:G^{''}\Rightarrow G^{'}$, such that
\begin{center}
\scalebox{0.9}[0.85]{\includegraphics{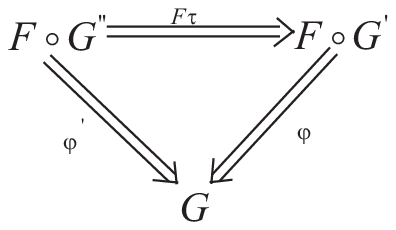}}
\end{center}
commutes.
\item The coroot of $\alpha$ is an $\cR$-2-module $Coroot \alpha$,
and an $\cR$-homomorphism $F:\cB\rightarrow Coroot\alpha$ with
$F\ast\alpha=1_{0}:0\Rightarrow 0:\cA\rightarrow Coroot\alpha$. For
any other $G:\cB\rightarrow\cD$, with $G\ast\alpha=1_{0}$, there
exist an $\cR$-homomorphism $G^{'}:Coroot\alpha\rightarrow\cD$, and
an invertible 2-morphism $\varphi:G^{'}\circ F\Rightarrow G$ as in
the following diagram
\begin{center}
\scalebox{0.9}[0.85]{\includegraphics{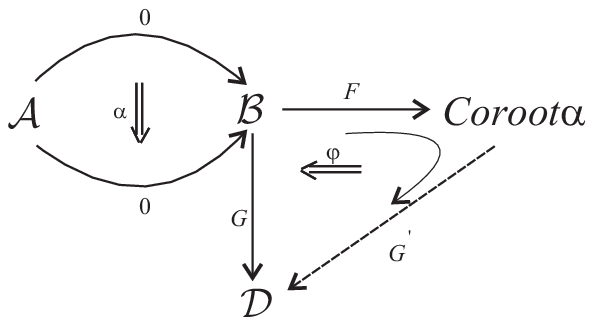}}
\end{center}
The pair $(G^{'},\varphi)$ is unique up to an invertible 2-morphism,
i.e. if $(G^{''},\varphi^{'})$ satisfies the same conditions as
$(G^{'},\varphi)$, there exists a unique $\tau:G^{''}\Rightarrow
G^{'}$, such that
\begin{center}
\scalebox{0.9}[0.85]{\includegraphics{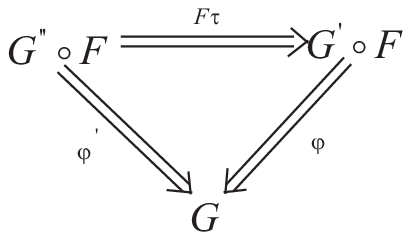}}
\end{center}
commutes.
\end{itemize}
\end{Def}

Next, we will use definitions 16-17 to give the existence of (Co)Pip
and (Co)Root, based on the results of symmetric 2-groups in \cite
{8,18}.

\begin{itemize}
\item Existence of pip of F.\\
(i)\ There is a category $PipF$ consisting the following data:

$\cdot$\ Object is a morphism $a:0\rightarrow 0$ in $\cA$, such that
$Fa=1_{F0}$, where 0 is the unit object of $\cA$.

$\cdot$\ Morphism of $\varphi:a\Rightarrow b:0\rightarrow 0$ is just
the identity if and only if $a=b$.

$\cdot$\ Composition of morphisms is a composition of identity
morphisms in $\cA$ which is also identity.

Obviously, the above data satisfy the necessary conditions in the
definition of category, and $PipF$ is also a groupiod.

(ii)\ $PipF$ is a symmetric 2-group with monoidal structure as
follows:

$PipF$ has a unit object $1:0\rightarrow 0$, which is the identity
of 0.

There is a bifunctor
\begin{align*}
&\hspace{0.3cm}+:PipF\times PipF\longrightarrow PipF\\
&(a:0\rightarrow0,b:0\rightarrow0)\mapsto a+b\triangleq a\circ
b:0\rightarrow 0\\
\end{align*}
Also, $+$ maps identity morphism to identity morphism in $PipF$.
Moreover, there are natural identities:
\begin{align*}
&(a+b)+c=(a\circ b)\circ c=a\circ(b\circ c)=a+(b+c),\\
&1+a=1\circ a=a,\\
&a+1=a\circ 1=a
\end{align*}
satisfy the Mac Lane coherence conditions(Fig.1-2.).

For any object $a:0\rightarrow 0$ in $PipF$, which is a morphism in
$\cA$,
and since $\cA$ is a groupoid, $a$ is invertible.\\
For any two objects $a,b:0\rightarrow0$ in $PipF$, there exists
identity $a+b=b+a$, since $a,\ b$ are endomorphisms of 0.

(iii) $PipF$ is an $\cR$-2-module.

There is a bifunctor:
\begin{align*}
&\hspace{0.5cm}\star:\cR\times PipF\longrightarrow PipF\\
&\hspace{1.2cm}(r,0\xrightarrow[]{a}0)\mapsto r\star a\\
&(r_{1}\xrightarrow[]{\varphi}r_{2},a\xrightarrow[]{id}a)\mapsto
0\xrightarrow[]{\varphi\star id }0
\end{align*}
where $r\star a=z_{r}\circ r\cdot a\cdot z_{r}^{-1}$ is the
composition
$$0\xrightarrow[]{z_{r}^{-1}}r\cdot0\xrightarrow[]{r\cdot
a}r\cdot0\xrightarrow[]{z_{r}}0,$$ such that $F(r\star a)=1_{F0}.$
$\varphi\star id=z_{r_{2}}\ast(\varphi\cdot id)\ast
z_{r_{1}}^{-1}=id.$

In fact, using Fig.37. and $Fa=1_{F0},$\\
$F(r\star a)=F(z_{r}\circ r\cdot a\cdot z_{r}^{-1})=F(z_{r})\circ
F(r\cdot a)\circ F(z_{r}^{-1})=F(z_{r})\circ F_{2}^{-1}\circ r\cdot
F(a)\circ F_{2}\circ F(z_{r}^{-1})=1_{F0}.$

Moreover, there are natural isomorphisms:

(1)\ $(r_{1}r_{2})\star a =r_{1}\star(r_{2}\star a)$.

In fact,
\begin{align*}
&r\star(a_{1}+a_{2})=r\star(a_{1}\circ a_{2})\triangleq
z_{r}\circ(r\cdot (a_{1}\circ a_{2}))\circ
z_{r}^{-1}=z_{r}\circ(r\cdot a_{1}\circ r\cdot a_{2})\circ
z_{r}^{-1}\\
&\hspace{2cm}=z_{r}\circ r\cdot a_{1}\circ z_{r}^{-1}\circ z_{r}
r\cdot a_{2}\circ z_{r}^{-1}=r\star a_{1}\circ r\star a_{1}=r\star
a_{1}+ r\star a_{1},\\
&(r_{1}r_{2})\star a\triangleq z_{r_{1}r_{2}}\circ
((r_{1}r_{2})\cdot a)\circ z_{r_{1}r_{2}}^{-1}\\
&\hspace{2cm}=z_{r_{1}}\circ (r_{1}\cdot z_{r_{2}})\circ
b_{r_{1},r_{2},0}\circ ((r_{1}r_{2})\cdot a)\circ
b_{r_{1},r_{2},0}^{-1}\circ (r_{1}\cdot z_{r_{2}}^{-1})\circ
z_{r{1}}^{-1},\\
&r_{2}\star a\triangleq z_{r_{2}}\circ (r_{2}\circ a)\circ
z_{r_{2}}^{-1},\\
&r_{1}\star (r_{2}\star a)\triangleq z_{r_{1}}\circ r_{1}\cdot
(r_{2}\star a)\circ z_{r_{1}}^{-1}=z_{r_{1}}\circ r_{1}\cdot
z_{r_{2}}\circ r_{1}\cdot(r_{2}\cdot a)\circ r_{1}\cdot
z_{r_{2}}^{-1}\circ z_{r_{1}}^{-1}.
\end{align*}
Next, we will check $b_{r_{1},r_{2},0}\circ ((r_{1}r_{2})\cdot
a)\circ b_{r_{1},r_{2},0}^{-1}= r_{1}\cdot(r_{2}\cdot a)$.\\
Let $F,G:\cR\times\cR\times \cA\rightarrow \cA$ be functors, given
by,
\begin{align*}
&F(r_{1},r_{2},0)\triangleq (r_{1}r_{2})\cdot 0,\\
&G(r_{1},r_{2},0)\triangleq r_{1}\cdot(r_{2}\cdot 0).
\end{align*}
for $r_{1},r_{2}\in obj(\cR)$. $b_{-,-,-}$ is a natural
transformation from F to G, so for $a:0\rightarrow 0$ in $\cA$, we
have the following commutative diagram:
\begin{center}
\scalebox{0.9}[0.85]{\includegraphics{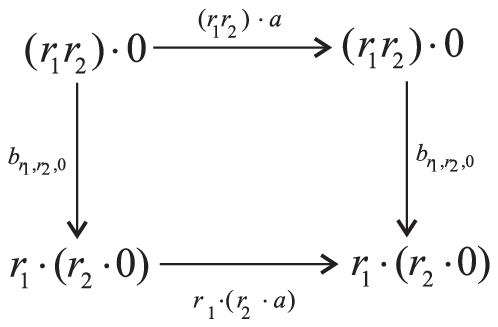}}
\end{center}

(2)\ $(r_{1}+r_{2})\star a=r_{1}\star a+r_{2}\star a.$

In fact,
\begin{align*}
&(r_{1}+r_{2})\star a\triangleq z_{r_{1}+r_{2}}\circ
(r_{1}+r_{2})\cdot a\circ
z_{r_{1}+r_{2}}^{-1}=\\
&\hspace{2cm}(z_{r_{1}}+z_{r_{2}})\circ b_{0}^{r_{1},r_{2}}\circ
(r_{1}+r_{2})\cdot a\circ (b^{r_{1},r_{2}}_{0})^{-1}\circ
(z_{r_{1}}^{-1}+z_{r_{2}}^{-1}).
\end{align*}
Next, we will check $b_{0}^{r_{1},r_{2}}\circ (r_{1}+r_{2})\cdot
a\circ (b^{r_{1},r_{2}}_{0})^{-1}=r_{1}\cdot a+r_{2}\cdot a$.\\
Let $F,G:\cR\times\cR\times \cA\rightarrow \cA$ be functors, given
by, for $r_{1},r_{2}\in obj(\cR)$,
\begin{align*}
&F(r_{1},r_{2},0)\triangleq (r_{1}+r_{2})\cdot 0,\\
&G(r_{1},r_{2},0)\triangleq r_{1}\cdot0+r_{2}\cdot 0.
\end{align*}
$b^{-,-}_{-}$ is a natural transformation from F to G, so for
$a:0\rightarrow 0$ in $\cA$, we have the following commutative
diagram:
\begin{center}
\scalebox{0.9}[0.85]{\includegraphics{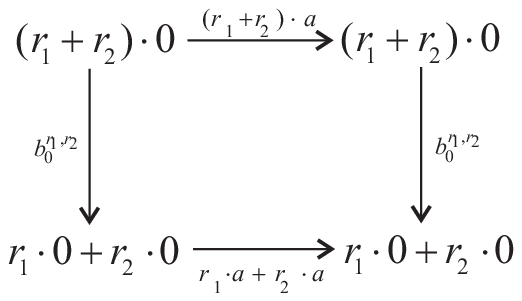}}
\end{center}

Similarly, we have (3)\  $1\star a=a,$ (4)\ $ r\star 1=1$.

The above natural isomorphisms are identities, so they make
Fig.18.-31. commute, then $PipF\in obj(\cR$-2-Mod).

(iv) $PipF$ is the pip of F.

There are zero homomorphisms:
\begin{align*}
&0:PipF\longrightarrow \cA\\
&\hspace{0.4cm}0\xrightarrow[]{a}0\mapsto 0.
\end{align*}
\begin{Rek}
For one object $a\in PipF$, it is a morphism in $\cA$, so we can
consider the above  zero morphisms as one maps to the source of $a$,
another maps to the target of $a$.
\end{Rek}
A 2-morphism
\begin{align*}
&\sigma:0\Rightarrow 0\\
&(\sigma)_{a}\triangleq a
\end{align*}
such that, $(1_{F}\ast\sigma)_{a}=(1_{F})_{0(a)}\circ
F((\sigma)_{a})=1_{F0}$, then $F\ast \sigma=1_{0}$.

If $\cD\in obj(\cR$-2-Mod), and $\alpha:0\Rightarrow
0:\cD\rightarrow \cA$, with $F\ast\alpha=1_{0}$. There exists an
$\cR$-homomorphism $G:\cD\rightarrow PipF$, by $G(d)\triangleq
\alpha_{d}:0\rightarrow 0$, and since $F\ast\alpha=1_{0}$, so
$F(\alpha_{d})=1_{F0}$. Also for any $d\in\cD,\ (\sigma\ast
G)_{d}={\sigma}_{Gd}=G(d)=\alpha_{d}$, i.e. $\sigma\ast G=\alpha$.

From the given $\cR$-homomorphism G, it is easy to see that G is
unique up to an invertible 2-morphism.

\item Existence of the copip of F.

(i) There is a category consisting of the following data:

$\cdot$\ A unique object is denoted by $\divideontimes$.

$\cdot$\ Morphism from $\divideontimes$ to  $\divideontimes$ is the
object $B\in\cB$. Two morphisms
$B_{1},B_{2}:\divideontimes\rightarrow\divideontimes$ are equal, if
there exist $A\in obj(\cA)$, and $b:B_{1}\rightarrow FA+B_{2}$.
Denote the equivalence class of morphisms by [B].

$\cdot$\ Composition of morphisms:

Let
$\divideontimes\xrightarrow[]{[B_{1}]}\divideontimes\xrightarrow[]{[B_{2}]}\divideontimes$
be morphisms in $CopipF$. We have
$[B_{2}]\circ[B_{1}]\triangleq[B_{1}+B_{2}]$, which is well-defined.
In fact, if $B_{1},B_{1}^{'}$ are equal, i.e. $\exists\ A_{1}\in
obj(\cA)$, and $b_{1}:B_{1}\rightarrow FA_{1}+B_{1}^{'}$. There
exist $A_{1}\in obj(\cA)$, and $b:B_{1}+B_{2}\rightarrow
(FA_{1}+B_{1}^{'})+B_{2}\rightarrow FA_{1}+(B_{1}^{'}+B_{2})$, so
$B_{1}+B_{2},\ B_{1}^{'}+B_{2}$ are equal.

The above data satisfy the following axioms:

(1) For the unique object $\divideontimes$, there exists an identity
morphism $1:\divideontimes\rightarrow\divideontimes$, which in fact
is the unit object $0$ of $\cB$, such that for any morphism [B],
there are $[B]\circ 1=[0+B]=[B], 1\circ [B]=[B+0]=[B]$. We will to
show the first equality. In fact, there exist $0\in obj(\cA)$, and
$b:0+B\xrightarrow[]{F_{0}^{-1}+1_{B}}F0+B$.

(2) Associativity of composition.

Given morphisms
$\divideontimes\xrightarrow[]{[B_{1}]}\divideontimes\xrightarrow[]{[B_{2}]}
\divideontimes\xrightarrow[]{[B_{3}]}\divideontimes$,
$[B_{3}]\circ([B_{2}]\circ [B_{1}])\triangleq(B_{1}+B_{2})+B_{3}$ is
equal to $([B_{3}]\circ [B_{2}])\circ [B_{1}]\triangleq
B_{1}+(B_{2}+B_{3})$, since there exist $0\in obj(\cA)$, and
morphism
$b:(B_{1}+B_{2})+B_{3}\xrightarrow[]{<B_{1},B_{2},B_{3}>}B_{1}+(B_{2}+B_{3})
\xrightarrow[]{1_{B_{1}+(B_{2}+B_{3})}}0+B_{1}+(B_{2}+B_{3})
\xrightarrow[]{F_{0}^{-1}+1_{B_{1}+(B_{2}+B_{3})}}F0+B_{1}+(B_{2}+B_{3})$.

For any morphism [B] in $CopipF$, $B\in obj(\cB)$, and $\cB$ is a
2-group, so B is invertible, then [B] is invertible, i.e. $CopipF$
is a groupoid.

(ii) $CopipF$ is a symmetric 2-group.

There is a bifunctor
\begin{align*}
&+:CopipF\times CopipF\longrightarrow CopiF\\
&\hspace{2.6cm}(\divideontimes,\divideontimes)\mapsto\divideontimes\\
&(\divideontimes\xrightarrow[]{[B_{1}]}\divideontimes,\divideontimes\xrightarrow[]{[B_{2}]}\divideontimes)\mapsto
\divideontimes\xrightarrow[]{[B_{1}+B_{2}]}\divideontimes
\end{align*}
Moreover, the natural isomorphisms are identities, so they satisfy
the Mac Lane coherence conditions. The inverse of an object is just
itself.

(iii) $CopipF$ is an $\cR$-2-module.

We can give the trivial bifunctor
\begin{align*}
&\star:\cR\times CopipF\longrightarrow CopipF\\
&\hspace{1.7cm}(r,\divideontimes)\mapsto \divideontimes
\end{align*}

(iv) $CopipF$ is the copip of F.

There are two zero $\cR$-morphisms
\begin{align*}
&0:\cB\rightarrow CopipF\\
&\hspace{0.5cm}B\mapsto \divideontimes
\end{align*}
\begin{Rek}
For one object $B\in\cB$, it is a morphism in $CopipF$, so we can
consider the above  zero morphisms as one maps to the source of B,
another maps to the target of B.
\end{Rek}
There is a morphism between the above two zero $\cR$-homomorphisms
\begin{align*}
&\sigma:0\Rightarrow 0\\
&\sigma_{B}\triangleq B:\divideontimes\rightarrow\divideontimes
\end{align*}
such that $\sigma\ast F=1_{0}$.

If $\cD\in obj(\cR$-2-Mod), and $\alpha:0\Rightarrow
0:\cB\rightarrow\cD$, with $\alpha\ast F=1_{0}$. There exists an
$\cR$-homomorphism
\begin{align*}
&G:CopipF\longrightarrow\cD\\
&\hspace{1.7cm}\divideontimes\mapsto G(\divideontimes)\equiv0,\\
&[B]:\divideontimes\rightarrow\divideontimes\mapsto
G([B])\triangleq\alpha_{B}:0\rightarrow 0
\end{align*}
where 0 is the unit object of $\cD$. For any $B\in (obj\cB),\
(G\ast\sigma)_{B}=(1_{G})_{0(B)}\circ
G(\sigma_{B})=(1_{G})_{\divideontimes}\circ G(B)=G(B)=\alpha_{B}$,
i.e. $G\ast\sigma =\alpha$.\\
From the definition of G, G is unique up to an invertible
2-morphism.
\end{itemize}

\begin{Def}(\cite {18}, Proposition\  179.)
Let $\cC$ be a $Gpd^{*}$-category with zero object and all the
kernels and cokernels. we say that $\cC$ is 2-Puppe-exact if the
following property holds.

For every morphism $f:A\rightarrow B$ in $\cC$,
$\overline{\omega}_{f}$ and $\omega_{f}$ are equivalent in the
following diagrams:
\begin{center}
\scalebox{0.9}[0.85]{\includegraphics{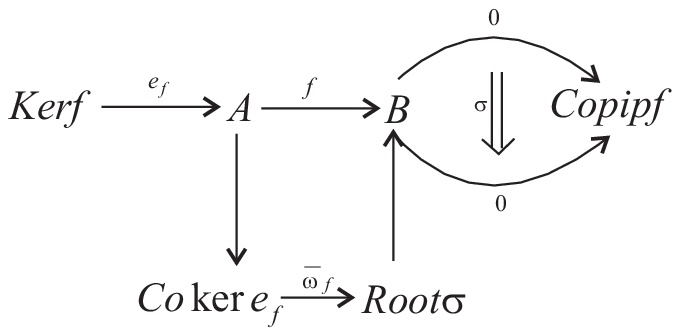}}
\scalebox{0.9}[0.85]{\includegraphics{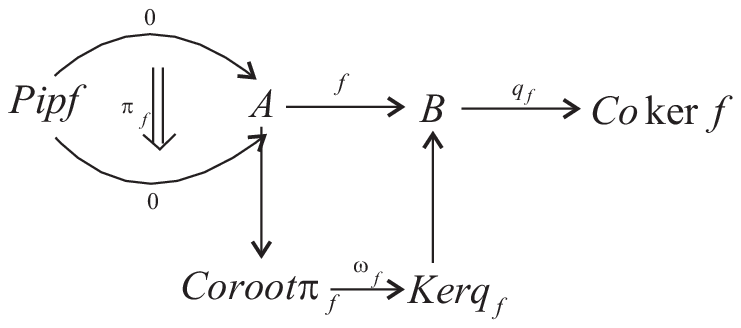}}
\end{center}
\end{Def}

\begin{Def}(\cite {18},Definition 183.) A 2-abelian $Gpd$-category is a 2-Puppe-exact
$Gpd^{*}$-category which has all the finite products and coproducts.
\end{Def}
The following four Propositions are the factorization systems in
2-category ($\cR$-2Mod) in the sense of \cite{7,11,18}.
\begin{Prop}
Every 1-morphism $F:\cA\rightarrow\cB$ in ($\cR$-2-Mod) factors as
the following composite, where $\widehat{E}_{F}$ is surjective,
$\widehat{\Omega}_{F}$ is an equivalence, and $\widehat{M}_{F}$ is
full and faithful.
$$
\cA\xrightarrow[]{\widehat{E}_{F}}Im^{1}_{pl}F\xrightarrow[]{\widehat{\Omega}_{F}}Im^{2}_{pl}F
\xrightarrow[]{\widehat{M}_{F}}\cB.
$$
\end{Prop}
\begin{proof}
\noindent\textbf{Step 1.} The $\cR$-2-module $Im^{1}_{pl}F$ is
described in the following way:
\begin{itemize}
\item  Category $Im^{1}_{pl}F$ consists of:

$\cdot$\ Objects are those of $\cA$.

$\cdot$\ Morphism $f:A\rightarrow A^{'}$ in $Im^{1}_{pl}F$ is a
morphism $F(f):FA\rightarrow FA^{'}$ in $\cB$. The composition of
morphisms are those of $\cB$, i.e. morphisms
$A_{1}\xrightarrow[]{f_{1}}A_{2}\xrightarrow[]{f_{2}}A_{3}$ in
$Im^{1}_{pl}F$ are
$FA_{1}\xrightarrow[]{F(f_{1})}FA_{2}\xrightarrow[]{F(f_{2})}FA_{3}$
in $\cB$, so the composition $f_{2}\circ f_{1}:A_{1}\rightarrow
A_{3}$ is $F(f_{2}\circ f_{1})=F(f_{2})\circ F(f_{1})$ in $\cB$.

The above ingredients satisfying the following axioms:

(1) For any object $A\in Im^{1}_{pl}F$, there exist an identity
morphism $1_{A}:A\rightarrow A$ given by $F(1_{A})\triangleq
1_{FA}:FA\rightarrow FA$ in $\cB$, such that for any morphism
$f:A_{1}\rightarrow A_{2}$, we have $f\circ 1_{A_{1}}=f,\
1_{A_{2}}\circ f=f$, since $F(f\circ1_{A_{1}})=F(f)\circ
F(1_{A_{1}})=F(f)\circ 1_{FA_{1}}=F(f),\ F(1_{A_{2}}\circ
f)=F(1_{A_{2}})\circ F(f)=1_{FA_{2}}\circ F(f)=F(f)$.

(2) Given morphisms
$$
A_{1}\xrightarrow[]{f_{1}}A_{2}\xrightarrow[]{f_{2}}A_{3}\xrightarrow[]{f_{3}}A_{4},
$$
for $F((f_{3}\circ f_{2})\circ f_{1})=F(f_{3})\circ F(f_{2})\circ
F(f_{1})=F(f_{3}\circ (f_{2}\circ f_{1}))$, so the associativity of
composition is true.

\item $Im^{1}_{pl}F$ is a symmetric monoidal groupoid.

The unit object of $Im^{1}_{pl}F$ is just the unit object of $\cA$.

There is a bifunctor
\begin{align*}
&+:Im^{1}_{pl}F\times Im^{1}_{pl}F\longrightarrow
Im^{1}_{pl}F \\
&\hspace{2.3cm}(A,A^{'})\mapsto A+A^{'},\\
&(A_{1}\xrightarrow[]{f}A_{2},A_{1}^{'}\xrightarrow[]{f^{'}}A_{2}^{'})
\mapsto A_{1}+A_{1}^{'}\xrightarrow[]{f+f^{'}}A_{2}+A_{2}^{'}
\end{align*}
where $A+A^{'}$ is an addition of objects of $\cA$, $f+f^{'}$ is a
composition morphism
$F(f+f^{'}):F(A_{1}+A_{1}^{'})\xrightarrow[]{F_{+}}FA_{1}+FA_{1}^{'}
\xrightarrow[]{Ff+Ff^{'}}FA_{2}+FA_{2}^{'}\xrightarrow[]{F_{+}^{-1}}F(A_{2}+A_{2}^{'}).$

Moreover, there are natural isomorphisms:
\begin{align*}
&<A_{1},A_{2},A_{3}>:(A_{1}+A_{2})+A_{3}\rightarrow
A_{1}+(A_{2}+A_{3}),\\
&\hspace{2.2cm}l_{A}:0+A\rightarrow A,\\
&\hspace{2.2cm}r_{A}:A+0\rightarrow A,\\
&\hspace{1.6cm}c_{A_{1},A_{2}}:A_{1}+A_{2}\rightarrow A_{2}+A_{1}
\end{align*}
defined by the images of the natural isomorphisms in $\cB$ under F.
As $\cB$ is a symmetric monoidal groupoid, so they satisfy the Mac
Lane coherence conditions, and $Im^{1}_{pl}F$ is a symmetric
monoidal groupoid.

\item $Im^{1}_{pl}F$ is a symmetric 2-group.

We need to show every object of $Im^{1}_{pl}F$ is invertible. In
fact, for any object $A\in Im^{1}_{pl}F$, $A\in obj(\cA)$, and $\cA$
is 2-group, there exist $A^{*}\in obj(\cA)$ and natural isomorphism
$\eta_{A}:A^{*}+A\rightarrow 0$, so there are $A^{*}\in
obj(Im^{1}_{pl}F)$, and composition isomorphism
$\eta_{FA}:FA^{*}+FA\xrightarrow[]{F_{+}^{-1}}F(A^{*}+A)
\xrightarrow[]{F(\eta_{A})}F0\xrightarrow[]{F_{0}}0$ in $\cB$, which
gives a natural isomorphism in $Im^{1}_{pl}F$.

\item $Im^{1}_{pl}F$ is an $\cR$-2-module.\\

There is a bifunctor
\begin{align*}
&\ \ \ \cdot \hspace{0.6cm} \cR\times Im^{1}_{pl}F\longrightarrow Im^{1}_{pl}F\\
&\hspace{2.3cm}(r,A)\mapsto r\cdot A,\\
&(r_{1}\xrightarrow[]{\varphi}r_{2},A_{1}\xrightarrow[]{f}A_{2})\mapsto
r_{1}\cdot A_{1}\xrightarrow[]{\varphi\cdot f}r_{2}\cdot A_{2}
\end{align*}
where $r\cdot A$ is the operation of $\cR$ on $\cA$, $\varphi\cdot
f:r_{1}\cdot A_{1}\rightarrow r_{2}\cdot A_{2}$ is a composition
$F(\varphi\cdot f):F(r_{1}\cdot
A_{1})\xrightarrow[]{F_{2}}r_{1}\cdot
FA_{1}\xrightarrow[]{\varphi\cdot F(f)}r_{2}\cdot
FA_{2}\xrightarrow[]{F_{2}^{-1}}F(r_{2}\cdot A_{2}) $ in $\cB$.

Moreover, there are natural isomorphisms:
\begin{align*}
&F(a_{A_{1},A_{2}}^{r}):r\cdot(A_{1}+A_{2})\rightarrow r\cdot
A_{1}+r\cdot A_{2},\\
&\hspace{0.2cm}F(b_{A}^{r_{1},r_{2}}):(r_{1}+r_{2})\cdot
A\rightarrow r_{1}\cdot
A+r_{2}\cdot A,\\
&F(b_{r_{1},r_{2},A}):(r_{1}r_{2})\cdot A\rightarrow
r_{1}\cdot(r_{2}\cdot A),\\
&\hspace{0.8cm}F(i_{A}):I\cdot A\rightarrow A,\\
&\hspace{0.8cm}F(z_{r}):r\cdot 0\rightarrow 0
\end{align*}
defined by the images of natural isomorphisms in $\cA$ under F.
Since $\cA$ is an $\cR$-2-module and F is a functor, so Fig.18.-31.
commute.
\end{itemize}

\noindent\textbf{Step 2.} The $\cR$-2-module $Im^{2}_{pl}F$ is
described in the following way:
\begin{itemize}
\item  Category $Im^{2}_{pl}F$ consists of the following data:

$\cdot$ Objects are the triples $(A,\varphi, B)$, where $A\in
obj(\cA),\ B\in obj(\cB),\ \varphi:FA\rightarrow B$.

$\cdot$ Morphism of $(A_{1},\varphi_{1},B_{1})\rightarrow
(A_{2},\varphi_{2},B_{2})$ is the morphism $g:B_{1}\rightarrow
B_{2}$ in $\cB$. Composition of morphisms and identity morphism are
those of $\cB$, and associativity of composition naturally holds.

\item $Im^{2}_{pl}F$ is a symmetric monoidal groupiod.

The unit object is $(0,F_{0},0)$, where the first 0 is the unit
object of $\cA$, the last one is the unit object of $\cB$,
$F_{0}:F0\rightarrow 0$.

There is a bifunctor
\begin{align*}
&\hspace{1cm}+:Im^{2}_{pl}F\times Im^{2}_{pl}F\longrightarrow Im^{2}_{pl}F\\
&((A_{1},\varphi_{1},B_{1}),(A_{2},\varphi_{2},B_{2}))\mapsto
(A_{1},\varphi_{1},B_{1})+(A_{2},\varphi_{2},B_{2})\triangleq(A,\varphi,B),\\
&\hspace{2.8cm}(g\hspace{0.4cm},\hspace{0.3cm}g^{'})\mapsto g+g^{'}
\end{align*}
where $A\triangleq A_{1}+A_{2},\ B\triangleq B_{1}+B_{2},\ \varphi$
is a composition
$F(A_{1}+A_{2})\xrightarrow[]{F_{+}}FA_{1}+FA_{1}\xrightarrow[]{\varphi_{1}+\varphi_{2}}B_{1}+B_{2}$,
$g+g^{'}$ is a morphism  of $g,\ g^{'}$ under the monoidal structure
of $\cB$.

Moreover, there are natural isomorphisms:
\begin{align*}
&<(A_{1},\varphi_{1},B_{1}),(A_{2},\varphi_{2},B_{2}),(A_{3},\varphi_{3},B_{3})>\triangleq
<B_{1},B_{2},B_{3}>:(
(A_{1},\varphi_{1},B_{1})+\\
&\hspace{2cm}(A_{2},\varphi_{2},B_{2}))+
(A_{3},\varphi_{3},B_{3})\rightarrow(A_{1},\varphi_{1},B_{1})+((A_{2},\varphi_{2},B_{2})+(A_{3},\varphi_{3},B_{3})),\\
& l_{(A,\varphi,B)}\triangleq
l_{B}:(0,F_{0},0)+(A,\varphi,B)\rightarrow (A,\varphi,B),\\
&r_{A,\varphi, B}\triangleq
r_{B}:(A,\varphi,B)+(0,F_{0},0)\rightarrow (A,\varphi,B),\\
&c_{(A_{1},\varphi_{1},B_{1}),(A_{2},\varphi_{2},B_{2})}\triangleq
c_{B_{1},B_{2}}:(A_{1},\varphi_{1},B_{1})+(A_{2},\varphi_{2},B_{2})\rightarrow
(A_{2},\varphi_{2},B_{2})+(A_{1},\varphi_{1},B_{1})
\end{align*}
 where $<B_{1},B_{2},B_{3}>,\ l_{B},\ r_{B},\ c_{B_{1},B_{2}}$ are
 natural isomorphisms in $\cB$, satisfy the Mac Lane coherence conditions, and
 also $F$ is a functor, then $Im^{2}_{pl}F$ is a symmetric monoidal
 category.

 For any morphism $g:(A_{1},\varphi_{1},B_{1})\rightarrow
 (A_{2},\varphi_{2},B_{2})$ in $Im^{2}_{pl}F$ is a morphism
 $g:B_{1}\rightarrow B_{2}$ in $\cB$ and $\cB$ is a groupoid, there
 exists $g^{*}:B_{2}\rightarrow B_{1}$, such that $g^{*}\circ
 g=1_{B_{1}}$, so there is $g^{*}:(A_{2},\varphi_{2},B_{2})\rightarrow
 (A_{1},\varphi_{1},B_{1})$, such that $g^{*}\circ
 g=1_{(A_{1},\varphi_{1},B_{1})}$ in $Im^{2}_{pl}F$.

\item $Im^{2}_{pl}F$ is a symmetric 2-group.

For any object $(A,\varphi,B)$ in $Im^{2}_{pl}F$, $A\in obj(\cA),\
B\in obj(\cB),\ \varphi:FA\rightarrow B$, for $\cA,\ \cB$ are
2-groups, there exist $A^{*}\in obj(\cA),\ B^{*}\in obj(\cB)$, and
natural isomorphisms $\eta_{A}:A^{*}+A\rightarrow 0,\
\eta_{B}:B^{*}+B\rightarrow 0$. So there exist
$(A,\varphi,B)^{*}\triangleq(A^{*},\varphi^{*},B^{*})$, and natural
isomorphism $\eta_{(A,\varphi,B)}\triangleq\eta_{B}$.

\item $Im^{2}_{pl}F$ is an $\cR$-2-module.

There is a bifunctor
\begin{align*}
&\star:\cR\times Im^{2}_{pl}F\longrightarrow Im^{2}_{pl}F\\
&\hspace{0.3cm}(r,(A,\varphi,B))\mapsto
r\star(A,\varphi,B)\triangleq (r\cdot
A,\widetilde{r\cdot\varphi},r\cdot B)\\
&(r_{1}\xrightarrow[]{\pi}r_{2},(A_{1},\varphi_{1},B_{1})\xrightarrow[]{g}(A_{2},\varphi_{2},B_{2}))
\mapsto \pi\cdot g
\end{align*}
where $r\cdot A,\ r\cdot B$ are operations of $\cR$ on $\cA,\ \cB$,
respectively, $\widetilde{r\cdot\varphi}$ is a composition $F(r\cdot
A)\xrightarrow[]{F_{2}}r\cdot FA\xrightarrow[]{r\cdot\varphi}r\cdot
B,\ \pi\cdot g$ is under the operation of $\cR$ on $\cB$.

Moreover, there are natural isomorphisms given by the natural
isomorphisms in $\cB$ from its $\cR$-2-module structure.
\end{itemize}

\noindent\textbf{Step 3.} $Im^{1}_{pl}F$, $Im^{2}_{pl}F$ are
equivalent $\cR$-2-modules.

Define a functor
\begin{align*}
&\widehat{\Omega}_{F}: Im^{1}_{pl}F\longrightarrow Im^{2}_{pl}F\\
&\hspace{1.6cm}A\mapsto (A,1_{FA},FA),\\
&\hspace{0.4cm}A_{1}\xrightarrow[]{g} A_{2}\mapsto
FA_{1}\xrightarrow[]{g}FA_{2}
\end{align*}
From the definition of $\widehat{\Omega}_{F}$, we see that
$\widehat{\Omega}_{F}$ restricts on morphisms of $Im^{1}_{pl}F$ to
be identity, so $\widehat{\Omega}_{F}$ is a functor.

Also, there are natural morphisms:
\begin{align*}
&(\widehat{\Omega}_{F})_{+}\triangleq
F_{+}:\widehat{\Omega_{F}}(A_{1}+A_{2})=(A_{1}+A_{2},1_{F(A_{1}+A_{2})},F(A_{1}+A_{2}))\rightarrow
\widehat{\Omega}_{F}(A_{1})+\widehat{\Omega}_{F}(A_{2})\\
&\hspace{2cm}=(A_{1},1_{FA_{1}},FA_{1})+(A_{2},1_{FA_{2}},FA_{2})=(A_{1}+A_{2},-,FA_{1}+FA_{2}),\\
&(\widehat{\Omega}_{F})_{0}\triangleq
F_{0}:\widehat{\Omega}_{F}(0)=(0,1_{F0},F0)\rightarrow (0,F_{0},0),\\
&(\widehat{\Omega}_{F})_{2}\triangleq
F_{2}:\widehat{\Omega}_{F}(r\cdot A)=(r\cdot A,1_{F(r\cdot
A)},F(r\cdot A))\rightarrow r\star\widehat{\Omega}_{F}(A)=r\star(A,1_{FA},FA)\\
&\hspace{3.5cm }=(r\cdot A,-,r\cdot FA)
\end{align*}
such that
$(\widehat{\Omega}_{F},(\widehat{\Omega}_{F})_{+},(\widehat{\Omega}_{F})_{0},(\widehat{\Omega}_{F})_{2})$
 is an $\cR$-homomorphism, since F is.

Define a functor
\begin{align*}
&\hspace{2.3cm}\widehat{\Omega}_{F}^{-1}:Im^{1}_{pl}F\longrightarrow
Im^{2}_{pl}F\\
&\hspace{3cm}(A,\varphi,B)\mapsto A,\\
&(A_{1},\varphi_{1},B_{1})\xrightarrow[]{g}(A_{2},\varphi_{2},B_{2})\mapsto
\widehat{\Omega}_{F}^{-1}(g)
\end{align*}
where $\widehat{\Omega}_{F}^{-1}(g)=\varphi_{2}^{-1}\circ
g\circ\varphi_{1}$ is the composition
$FA_{1}\xrightarrow[]{\varphi_{1}}B_{1}\xrightarrow[]{g}B_{2}\xrightarrow[]{\varphi_{2}^{-1}}FA_{2}$.

For any identity morphism $1_{B}:(A,\varphi,B)\rightarrow
(A,\varphi, B)$ in $Im^{2}_{pl}F$, we have
$\widehat{\Omega}_{F}^{-1}(1_{B})=\varphi^{-1}\circ
1_{B}\circ\varphi=1_{FA}$. Given morphisms
$(A_{1},\varphi_{1},B_{1})\xrightarrow[]{g_{1}}(A_{2},\varphi_{2},B_{2})\xrightarrow[]{g_{2}}(A_{3},\varphi_{3},B_{3})$,
we have $\widehat{\Omega}_{F}^{-1}(g_2\circ
g_{1})=\varphi_{3}^{-1}\circ(g_{2}\circ g_{1})\circ\varphi_{1},\
\widehat{\Omega}_{F}^{-1}(g_{2})=\varphi_{3}^{-1}\circ
g_{2}\circ\varphi_{2},\
\widehat{\Omega}_{F}^{-1}(g_{1})=\varphi_{2}^{-1}\circ
g_{1}\circ\varphi_{1}$, so $\widehat{\Omega}_{F}^{-1}(g_{2}\circ
g_{1})=\widehat{\Omega}_{F}^{-1}(g_{2})\circ
\widehat{\Omega}_{F}^{-1}(g_{1})$. Then $\widehat{\Omega}_{F}^{-1}$
is a functor.

Also, there are natural morphisms:
\begin{align*}
&(\widehat{\Omega}_{F}^{-1})_{+}=id:\widehat{\Omega}_{F}^{-1}((A_{1},\varphi_{1},B_{1})+(A_{2},\varphi_{2},B_{2}))
=\widehat{\Omega}_{F}^{-1}(A_{1}+A_{2},-,B_{1}+B_{2})=A_{1}+A_{2}\\
&\hspace{2.2cm}\rightarrow
\widehat{\Omega}_{F}^{-1}(A_{1},\varphi_{1},B_{1})+\widehat{\Omega}_{F}^{-1}(A_{2},\varphi_{2},B_{2})=A_{1}+A_{2},\\
&(\widehat{\Omega}_{F}^{-1})_{0}=id:\widehat{\Omega}_{F}^{-1}(0,1_{F0},0)=0\rightarrow
0,\\
&(\widehat{\Omega}_{F}^{-1})_{2}=id:\widehat{\Omega}_{F}^{-1}(r\star(A,\varphi,B))=\widehat{\Omega}_{F}^{-1}(r\cdot
A,\widetilde{(r\cdot\varphi)},r\cdot B)=r\cdot A\rightarrow
r\cdot\widehat{\Omega}_{F}^{-1}(A,\varphi,B)=r\cdot A.
\end{align*}
Obviously,
$(\widehat{\Omega}_{F}^{-1},(\widehat{\Omega}_{F}^{-1})_{+},
(\widehat{\Omega}_{F}^{-1})_{0},(\widehat{\Omega}_{F}^{-1})_{2})$ is
an $\cR$-homomorphism.

Next, we will check $\widehat{\Omega}_{F}^{-1}\circ
\widehat{\Omega}_{F}=1,\
\widehat{\Omega}_{F}\circ\widehat{\Omega}_{F}^{-1}\Rightarrow 1$.
\begin{align*}
&(\widehat{\Omega}_{F}^{-1}\circ\widehat{\Omega}_{F})(A)=\widehat{\Omega}_{F}^{-1}(A,1_{FA},FA)=A,\  \forall\ A\in obj(Im^{1}_{pl}F),\\
&(\widehat{\Omega}_{F}^{-1}\circ
\widehat{\Omega}_{F})(g)=(\widehat{\Omega}_{F}^{-1})(g)=g,\ \forall\
g\in Mor(Im^{1}_{pl}F).
\end{align*}
There is a morphism of $\cR$-homomorphisms
$$
\tau:\widehat{\Omega}_{F}\circ\widehat{\Omega}_{F}^{-1}\Rightarrow
1:Im^{2}_{pl}F\longrightarrow Im^{2}_{pl}F,$$ given by
$$\tau_{(A,\varphi,B)}\triangleq
\varphi:(\widehat{\Omega}_{F}\circ\widehat{\Omega}_{F}^{-1})(A,\varphi,B)=\widehat{\Omega}_{F}(A)=(A,1_{FA},FA)\rightarrow
(A,\varphi,B).$$
For a morphism
$g:(A_{1},\varphi_{1},B_{1})\rightarrow (A_{2},\varphi_{2},B_{2}),\
(\widehat{\Omega}_{F}\circ\widehat{\Omega}_{F}^{-1})(g)=\widehat{\Omega}_{F}(\varphi_{2}^{-1}\circ
g\circ\varphi_{1})=\varphi_{2}^{-1}\circ g\circ\varphi_{1}$, i.e.
$\tau$ is a natural transformation.\\
\begin{align*}
&(\widehat{\Omega}_{F}\circ\widehat{\Omega}_{F}^{-1})((A_{1},\varphi_{1},B_{1})+(A_{2},\varphi_{2},B_{2}))
=(A_{1}+A_{2},1_{FA_{1}+FA_{2}},F(A_{1}+A_{2}),\\
&(\widehat{\Omega}_{F}\circ\widehat{\Omega}_{F}^{-1})(A_{1},\varphi_{1},B_{1})
+(\widehat{\Omega}_{F}\circ\widehat{\Omega}_{F}^{-1})(A_{2},\varphi_{2},B_{2})
=(A_{1}+A_{2},-,FA_{1}+FA_{2}),\\
&\tau_{(A_{1},\varphi_{1},B_{1})+(A_{2},\varphi_{2},B_{2})}=(\varphi_{1}+\varphi_{2})\circ
F_{+},\\
&\tau_{(A_{i},\varphi_{i},B_{i})}=\varphi_{i}, \mbox{ for i=1,\ 2.}
\end{align*}
Thus Fig.7 commutes. We can also get commutative diagram Fig.37. in
the similar ways.

Since $\varphi:FA\rightarrow B$ is a morphism in $\cB$, and $\cB$ is
a groupoid, so $\tau_{(A,\varphi,B)}$ is an isomorphism.

\noindent\textbf{Step 4.} There is a surjective $\cR$-homomorphism.

Define a functor
\begin{align*}
&\hspace{0.9cm}\widehat{E}_{F}:\cA\longrightarrow Im^{1}_{pl}F\\
&\hspace{1.8cm}A\mapsto A,\\
&f:A_{1}\rightarrow A_{2}\mapsto Ff:FA_{1}\rightarrow FA_{2}
\end{align*}
For any $1_{A}:A\rightarrow A$ in $\cA$,
$$
\widehat{E}_{F}(1_{A})=F(1_{A})=1_{FA}.
$$
Given morphisms
$A_{1}\xrightarrow[]{f_{1}}A_{2}\xrightarrow[]{f_{2}}A_{3}$ in
$\cA$,
$$
\widehat{E}_{F}(f_{2}\circ f_{1})=F(f_{2}\circ f_{1})=Ff_{2}\circ
Ff_{1}=\widehat{E}_{F}(f_{2})\circ\widehat{E}_{F}(f_{1}).
$$
There are natural isomorphisms:
\begin{align*}
&(\widehat{E}_{F})_{+}=id:(\widehat{E}_{F})(A_{1}+A_{2})=A_{1}+A_{2}=\widehat{E}_{F}(A_{1})+\widehat{E}_{F}(A_{2}),\\
&(\widehat{E}_{F})_{0}=id:\widehat{E}_{F}(0)=0,\\
&(\widehat{E}_{F})_{2}=id:\widehat{E}_{F}(r\cdot A)=r\cdot
A=r\cdot\widehat{E}_{F}(A).
\end{align*}
Obviously,
$(\widehat{E}_{F},(\widehat{E}_{F})_{+},(\widehat{E}_{F})_{0},(\widehat{E}_{F})_{2})$
is an $\cR$-homomorphism.

For every $\cC\in obj(\cR$-2-Mod), there is a functor
\begin{align*}
&\Phi:Hom(Im^{1}_{pl}F,\cC)\longrightarrow
Hom(\cA,\cC)\\
&\hspace{3cm} G\mapsto G\circ \widehat{E}_{F},\\
&\hspace{1.2cm}\alpha:G_{1}\Rightarrow G_{2}\mapsto
\alpha\ast\widehat{E}_{F}\triangleq\alpha\ast1_{\widehat{E}_{F}}
\end{align*}
 where $\alpha\ast1_{\widehat{E}_{F}}$ is the horizontal composition of
 2-morphisms. It is easy to check $\Phi$ is an $\cR$-homomorphism.

If $\alpha,\beta:G_{1}\Rightarrow G_{2}:Im^{1}_{pl}F\rightarrow\cC$,
such that $\alpha\ast \widehat{E}_{F}=\beta\ast\widehat{E}_{F}$,
i.e. $\forall\ A\in obj(\cA)$,
$$
(\alpha\ast\widehat{E}_{F})(A)=G_{2}((1_{\widehat{E}_{F}})_{A})\circ\alpha_{\widehat{E}_{F}(A)}
=1_{G_{2}A}\circ\alpha_{A}=\alpha_{A},
$$
$$
(\beta\ast\widehat{E}_{F})(A)=G_{2}((1_{\widehat{E}_{F}})_{A})\circ\beta_{\widehat{E}_{F}(A)}
=1_{G_{2}A}\circ\beta_{A}=\beta_{A}.
$$
Then $\alpha=\beta$. So we proved $\Phi$ is a faithful functor.

From above results, we get a surjective $\cR$-homomorphism
$\widehat{E}_{F}:\cA\rightarrow Im^{1}_{pl}F$.

\noindent\textbf{Step 5.} There is a fully faithful
$\cR$-homomorphism.

Define a functor
\begin{align*}
&\hspace{2.9cm}\widehat{M}_{F}:Im^{2}_{pl}F\longrightarrow\cB \\
&\hspace{3.5cm}(A,\varphi,B)\mapsto B,\\
&g:(A_{1},\varphi_{1},B_{1})\rightarrow
(A_{2},\varphi_{2},B_{2})\mapsto g:B_{1}\rightarrow B_{2}.
\end{align*}
For any identity morphism
$1_{(A,\varphi,B)}=1_{B}:(A,\varphi,B)\rightarrow (A,\varphi,B)$,
$$
\widehat{M}_{F}(1_{(A,\varphi,B)})=\widehat{M}_{F}(1_{B})=1_{B}=1_{(A,\varphi,B)}.
$$
For morphisms
$(A_{1},\varphi_{1},B_{1})\xrightarrow[]{g_{1}}(A_{2},\varphi_{2},B_{2})\xrightarrow[]{g_{2}}(A_{3},\varphi_{3},B_{3})$,
$$
\widehat{M}_{F}(g_{2}\circ g_{1})=g_{2}\circ
g_{1}=\widehat{M}_{F}(g_{2})\circ\widehat{M}_{F}(g_{1}).
$$
There are natural morphisms:
\begin{align*}
&(\widehat{M}_{F})_{+}=id:\widehat{M}_{F}((A_{1},\varphi_{1},B_{1})+(A_{2},\varphi_{2},B_{2}))
=\widehat{M}_{F}(A_{1}+A_{2},-,B_{1}+B_{2})=B_{1}+B_{2}\\
&\hspace{4cm}\rightarrow
\widehat{M}_{F}(A_{1},\varphi_{1},B_{1})+\widehat{M}_{F}(A_{2},\varphi_{2},B_{2}))=B_{1}+B_{2},\\
&(\widehat{M}_{F})_{0}=id:\widehat{M}_{F}(0,F_{0},0)=0,\\
&(\widehat{M}_{F})_{2}=id:\widehat{M}_{F}(r\star(A,\varphi,B))=\widehat{M}_{F}(r\cdot
A,-,r\cdot B)=r\cdot B=r\cdot\widehat{M}_{F}(A,\varphi,B).
\end{align*}
$(\widehat{M}_{F},(\widehat{M}_{F})_{+},(\widehat{M}_{F})_{0},(\widehat{M}_{F})_{2})$
is an $\cR$-homomorphism.

For any pairs of objects $(A_{1},\varphi_{1},B_{1}),\
(A_{2},\varphi_{2},B_{2})$ in $Im^{2}_{pl}F$ with
$g:B_{1}\rightarrow B_{2}$ in $\cB$, there exists
$g:(A_{1},\varphi_{1},B_{1})\rightarrow (A_{2},\varphi_{2},B_{2})$
in $Im^{2}_{pl}F$, such that $\widehat{M}_{F}(g)=g$. Thus
$\widehat{M}_{F}$ is full.

For given morphisms
$g_{1},g_{2}:(A_{1},\varphi_{1},B_{1})\rightarrow
(A_{2},\varphi_{2},B_{2})$ in $Im^{2}_{pl}F$ such that
$\widehat{M}_{F}(g_{1})=\widehat{M}_{F}(g_{2}):B_{1}\rightarrow
B_{2}$, since $\widehat{M}_{F}(g_{i})=g_{i},\ i=1,\ 2.$ So
$g_{1}=g_{2}$. Thus $\widehat{M}_{F}$ is faithful.

\noindent\textbf{Step 6.}
$F=\widehat{M}_{F}\circ\widehat{\Omega}_{F}\circ\widehat{E}_{F}.$

For any $A\in obj(\cA),$
$$
(\widehat{M}_{F}\circ\widehat{\Omega}_{F}\circ\widehat{E}_{F})(A)
=(\widehat{M}_{F}\circ\widehat{\Omega}_{F})(A)
=\widehat{M}_{F}(A,1_{FA},FA)=FA.
$$
For any morphism $f:A_{1}\rightarrow A_{2}$ in $\cA$,
$$
(\widehat{M}_{F}\circ\widehat{\Omega}_{F}\circ\widehat{E}_{F})(f)
=(\widehat{M}_{F}\circ\widehat{\Omega}_{F})(Ff)=\widehat{M}_{F}(Ff)=Ff.
$$
\end{proof}
\begin{Prop}
For each $\cR$-homomorphism $F:\cA\rightarrow\cB$ in ($\cR$-2-Mod),
$\widehat{E}_{F}:\cA\rightarrow Im^{1}_{pl}F$ is the cokernel of the
kernel of F.
\end{Prop}
\begin{proof}
We know that the kernel of F is $(KerF,\ e_{F},\ \varepsilon_{F})$,
where $KerF$ is a category  with objects are pairs $(A,a)$, where
$A\in obj(\cA),\ a:FA\rightarrow 0$, $e_{F}(A,a)=A,\
(\varepsilon_{F})_{(A,a)}=a$.
\begin{itemize}
\item Let us describe the $Cokere_{F}$ in the following ways:

$\cdot$\ Objects are the objects of $\cA$.

$\cdot$\ Morphism from $A_{1}$ to $A_{2}$ is the equivalence class
of the triple $(N,a,f)$, denote by $[N,a,f]$, where $N\in obj(\cA),\
a:FN\rightarrow 0$ is a morphism in $\cB$, $f:A_{1}\rightarrow
N+A_{2}$ in $\cA$, and for two morphisms
$(N_{1},a_{1},f_{1}),(N_{2},a_{2},f_{2}):A_{1}\rightarrow A_{2}$ are
equal if there exists a morphism $n:N_{1}\rightarrow N_{2}$ in
$\cA$, such that the following diagrams commute:
\begin{center}
\scalebox{0.9}[0.85]{\includegraphics{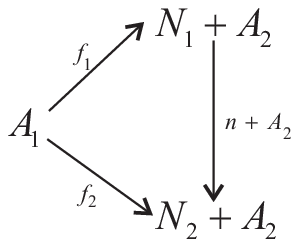}}{\footnotesize
Fig.42.}
\scalebox{0.9}[0.85]{\includegraphics{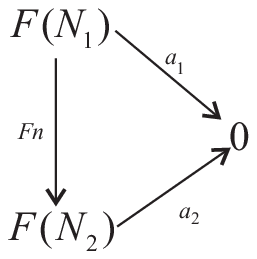}}{\footnotesize
Fig.43.}
\end{center}

$\cdot$\ Composition of morphisms. Let
$A_{1}\xrightarrow[]{[N_{1},a_{1},f_{1}]}A_{2}\xrightarrow[]{[N_{2},a_{2},f_{2}]}A_{3}$
be morphisms in $Cokere_{F}$,
$[N,a,f]\triangleq[N_{2},a_{2},f_{2}]\circ[N_{1},a_{1},f_{1}]$,
where $N\triangleq N_{1}+N_{2},\ a$ is the composition
$F(N_{1}+N_{2})\xrightarrow[]{F_{+}}FN_{1}+FN_{2}\xrightarrow[]{a_{1}+a_{2}}0+0=0,$
$f$ is the composition
$A_{1}\xrightarrow[]{f_{1}}N_{1}+A_{2}\xrightarrow[]{1_{
N_{1}}+f_{2}}N_{1}+(N_{2}+A_{3})\xrightarrow[]{<N_{1},N_{2},A_{3}>^{-1}}(N_{1}+N_{2})+A_{3}=N+A_{3}$.

This composition is well-defined, since
 $(N_{1},a_{1},f_{1}),(N_{1}^{'},a_{1}^{'},f_{1}^{'})$ are equal,
 i.e. $\exists\  n_{1}:N_{1}\rightarrow N_{1}^{'}$, such that
 $(n_{1}+1_{A_{2}})\circ f_{1}=f_{1}^{'},\ a_{1}^{'}\circ Fn_{1}=a_{1}.$
There exists $n\triangleq n_{1}+1_{N_{2}}:N_{1}+N_{2}\rightarrow
N_{1}^{'}+N_{2},$ such that Fig.42.-43. commute, then
$(N_{2},a_{2},f_{2})\circ(N_{1},a_{1},f_{1})$ is equal to
$(N_{2},a_{2},f_{2})\circ(N_{1}^{'},a_{1}^{'},f_{1}^{'})$.

\item There is a morphism of $\cR$-homomorphisms
$$
\delta:\widehat{E}_{F}\circ e_{F}\Rightarrow 0
$$
given by
$$
\delta_{(A,a)}\triangleq a:(\widehat{E}_{F}\circ
e_{F})(A,a)=(\widehat{E}_{F})(A)=A\rightarrow 0(A,a)=0,
$$
for any $(A,a)\in obj(KerF)$, which is a morphism in $Im^{1}_{pl}F$.

For any morphism $f:(A_{1},a_{1})\rightarrow (A_{2},a_{2})$, we have
$$
(\widehat{E}_{F}\circ e_{F})(f)=f:(\widehat{E}_{F}\circ
e_{F})(A_{1},a_{1})=A_{1}\rightarrow (\widehat{E}_{F}\circ
e_{F})(A_{2},a_{2})=A_{2}
$$
in $Im^{1}_{pl}F$ is the morphism $Ff:FA_{1}\rightarrow FA_{2}$ in
$\cB$, $\delta_{(A_{i},a_{i})}=a_{i},\ i=1,2.$ From $a_{2}\circ
Ff=a_{1}$, $\delta$ is a natural transformation.\\
Using the properties of F, it is easy to check $\delta$ is a
morphism of $\cR$-homomorphisms.

By the universal property of the cokernel, there is an
$\cR$-homomorphism
\begin{align*}
&\hspace{1cm}\Phi:Cokere_{F}\longrightarrow Im^{1}_{pl}F\\
&\hspace{3cm}A\mapsto A,\\
&[N,a,f]:A_{1}\rightarrow A_{2}\mapsto \Phi(N,a,f):A_{1}\rightarrow
A_{2}
\end{align*}
where $(N,a,f)$ is the representation element of equivalence class
of $[N,a,f]$, $\Phi(N,a,f)$ is the following composition of
morphisms in $\cB$,
$$
FA_{1}\xrightarrow[]{Ff}F(N+A_{2})\xrightarrow[]{F_{+}}FN+FA_{2}
\xrightarrow[]{a+1_{FA_{2}}}0+FA_{2}\xrightarrow[]{l_{FA_{2}}}FA_{2}.
$$
$\Phi$ is well-defined, if $(N_{1},a_{1},f_{1}),
 (N_{2},a_{2},f_{2})$ are equal, i.e. $\exists\ n:N_{1}\rightarrow N_{2},$ such that
$$
(n+1_{A_{2}})\circ f_{1}=f_{2},\ a_{2}\circ Fn=a_{1}.
$$
We get $\Phi(N_{1},a_{1},f_{1})=\Phi(N_{2},a_{2},f_{2})$ from the
following commutative diagrams
\begin{center}
\scalebox{0.9}[0.85]{\includegraphics{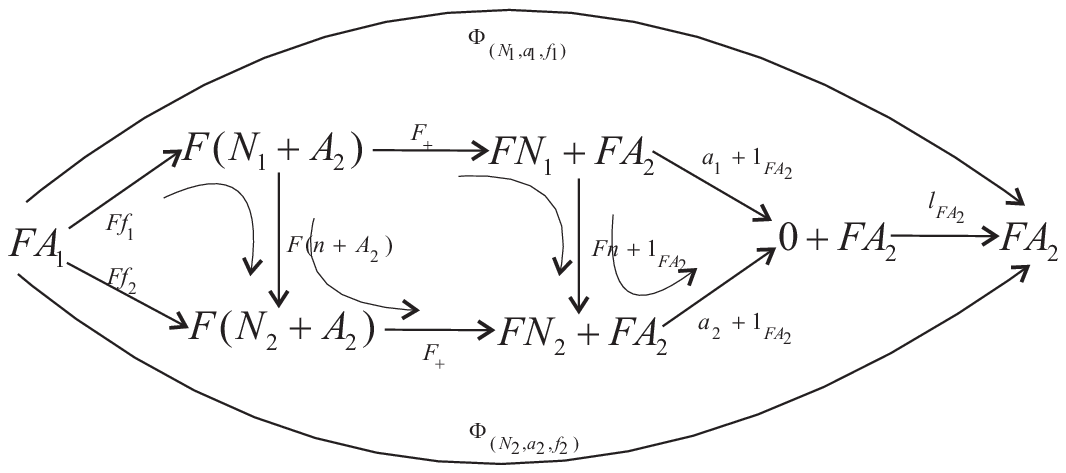}}
\end{center}
For any identity morphism $[0,F_{0},l_{a}^{-1}]:A\rightarrow A$ in
$Cokere_{F}$, we have
$$
\Phi(0,F_{0},l_{A}^{-1})=l_{FA}\circ(F+1_{FA})\circ F_{+}\circ
F(l_{A}^{-1}))=1_{FA}.
$$
For given morphisms
$A_{1}\xrightarrow[]{[N_{1},a_{1},f_{1}]}A_{2}\xrightarrow[]{[N_{2},a_{2},f_{2}]}A_{3}$
in $Cokere_{F}$, from the basic properties of morphisms in
$Cokere_{F}$ and $F$, we have
$$
\Phi((N_{2},a_{2},f_{2})\circ(N_{1},a_{1},f_{1}))=\Phi(N_{2},a_{2},f_{2})\circ\Phi(N_{1},a_{1},f_{1}).
$$

Moreover, there are natural morphisms:
\begin{align*}
&\Phi_{+}=id:\Phi(A_{1}+A_{2})=A_{1}+A_{2}=\Phi A_{1}+\Phi A_{2},\\
&\Phi_{0}=id:\Phi(0)=0,\\
&\Phi_{2}=id:\Phi(r\cdot A)=r\cdot A=r\cdot \Phi A.
\end{align*}
$\Phi$ is an $\cR$-homomorphism.

\item $\Phi$ is an equivalent $\cR$-homomorphism.

For every $\cC\in obj(\cR$-2-Mod), there is a functor
\begin{align*}
&\Psi:Hom(Im^{1}_{pl}F,\cC)\longrightarrow
Hom(Cokere_{F},\cC)\\
&\hspace{3cm} G\mapsto G\circ \Phi,\\
&\hspace{1.2cm}\alpha:G_{1}\Rightarrow G_{2}\mapsto \alpha\ast\Phi
\end{align*}
 where $\alpha\ast1_{\Phi}$ is the horizontal composition of
 2-morphisms.

 If $\alpha,\beta:G_{1}\Rightarrow G_{2}:Im^{1}_{pl}F\rightarrow\cC$,
such that $\alpha\ast \Phi=\beta\ast\Phi$, i.e. $\forall\ A\in
obj(Cokere_{F})$,
$$
(\alpha\ast\Phi)(A)=G_{2}((1_{\Phi})_{A})\circ\alpha_{\Phi(A)}
=1_{G_{2}A}\circ\alpha_{A}=\alpha_{A},
$$
$$
(\beta\ast\Phi)(A)=G_{2}((1_{\Phi})_{A})\circ\beta_{\Phi(A)}
=1_{G_{2}A}\circ\beta_{A}=\beta_{A}.
$$
Then $\alpha=\beta$. So $\Psi$ is a faithful functor, and $\Phi$ is
surjective.

For any two objects $A_{1},\ A_{2}$ in $Cokere_{F}$, and morphism
$g:FA_{1}\rightarrow FA_{2}$ in $Im^{1}_{pl}F$. Set
$N=A_{1}+A_{2}^{*}$, where $A_{2}^{*}$ is an inverse of $A_{2}$,
together with natural isomorphism
$\eta_{A_{2}}:A_{2}^{*}+A_{2}\rightarrow 0$, denote by
$\eta_{A_{2}}^{'}=\eta_{A_{2}}\circ
c_{A_{2},A_{2}^{*}}:A_{2}+A_{2}^{*}\rightarrow 0$, $a$ is the
composition
\begin{align*}
&F(A_{1}+A_{2}^{*})\xrightarrow[]{F_{+}}FA_{1}+FA_{2}^{*}
\xrightarrow[]{g+1_{F(A_{2}^{*})}}FA_{2}+F(A_{2}^{*})\xrightarrow[]{F_{+}^{-1}}F(A_{2}+A_{2}^{*})
\xrightarrow[]{F\eta_{A_{2}}^{'}}F0\xrightarrow[]{F_{0}}0
\end{align*}
and $f$ is the composition
$$
A_{1}\xrightarrow[]{r_{A_{1}}^{-1}}A_{1}+0\xrightarrow[]{1_{A_{1}}+\eta_{A_{2}}^{-1}}A_{1}+(A_{2}^{*}+A_{2})
\xrightarrow[]{<A_{1},A_{2}^{*},A_{2}>^{-1}}(A_{1}+A_{2}^{*})+A_{2}=N+A_{2}.
$$
After calculations, we have $\Phi(N,b,f)=g$, then $\Phi$ is full.

For two morphisms
$[N_{1},a_{1},f_{1}],[N_{2},a_{2},f_{2}]:A_{1}\rightarrow A_{2}$ in
$Cokere_{F}$, such that
$\Phi(N_{1},a_{1},f_{1})=\Phi(N_{2},a_{2},f_{2})$.

 Let n be the composition
\begin{align*}
&N_{1}\xrightarrow[]{r_{N_{1}}}N_{1}+0\xrightarrow[]{1+\eta_{A_{2}}^{'}}N_{1}+(A_{2}+A_{2}^{*})
\xrightarrow[]{<N_{1},A_{2},A_{2}^{*}>^{-1}}(N_{1}+A_{2})+A_{2}^{*}
\xrightarrow[]{f_{1}^{-1}+1_{A_{2}^{*}}}A_{1}+A_{2}^{*}\\
&\hspace{1cm}\xrightarrow[]{f_{2}+1_{A_{2}^{*}}}(N_{2}+A_{2})+A_{2}^{*}
\xrightarrow[]{<N_{2},A_{2},A_{2}^{*}>}N_{2}+(A_{2}+A_{2}^{*})
\xrightarrow[]{1_{N_{2}}+\eta_{A_{2}}^{'}}N_{2}+0\xrightarrow[]{l_{N_{2}}}N_{2},
\end{align*}
then $[N_{1},a_{1},f_{1}],\ [N_{2},a_{2},f_{2}]:A_{1}\rightarrow
A_{2}$ are equal in $Cokere_{F}$, thus $\Phi$ is faithful.

$\Phi:Cokere_{F}\rightarrow Im^{1}_{pl}F$ is an equivalent
$\cR$-homomorphism.
\end{itemize}
\end{proof}

\begin{Prop}
For each $\cR$-homomorphism $F:\cA\rightarrow\cB$ in ($\cR$-2-Mod),
$\widehat{M}_{F}:Im^{2}_{pl}F\rightarrow\cB$ is the root of the
copip of F.
\end{Prop}
\begin{proof}
The copip of F, is the $\cR$-2-module $CopipF$, together with
2-morphism $\sigma:0\Rightarrow 0:\cB\rightarrow CopipF$. As a
category, $CopipF$ has the unique object $\divideontimes$ and
morphisms are objects of $\cB$, and
$\sigma(B)=B:\divideontimes\rightarrow\divideontimes$.

\noindent\textbf{Step 1.} Let us describe the root of
$\sigma:0\Rightarrow 0:\cB\rightarrow CopipF$.
\begin{itemize}
\item\ The category $Root\sigma$ consists of the following data:

$\cdot$\ Objects are $B\in obj(\cB)$, such that $\sigma_{B}=B=0$ in
$CopipF$, i.e. there exist $A\in obj(\cA)$ and $g:B\rightarrow
FA+0$.

$\cdot$\ Morphisms are morphisms in $\cB$.

$\cdot$\ Composition of morphisms and the unit object are just the
composition of morphisms in $\cB$, and the unit object 0 in $\cB$,
respectively.

From the $\cR$-2-module structure of $\cB$, we can give $Root\sigma$
a $\cR$-2-module structure.

\item The $\cR$-2-module $Root\sigma$ is the root of $\sigma$.

There is an $\cR$-homomorphism
\begin{align*}
&\hspace{0.2cm}R:Root\sigma\longrightarrow\cB\\
&\hspace{1.6cm}B\mapsto B,\\
&g:B_{1}\rightarrow B_{2}\mapsto g:B_{1}\rightarrow B_{2}.
\end{align*}
Since $Root\sigma$ and $\cB$ have the same $\cR$-2-module structure,
then $R$ is an $\cR$-homomorphism.

Also, for any $B\in obj(Root\sigma),\ (\sigma\ast
R)(B)=\sigma(B)=B=0=1_{0}(B)$, i.e. $\sigma\ast R=1_{0}.$

For $\cD\in obj(\cR$-2-Mod) and an $\cR$-homomorphism
$G:\cD\rightarrow\cB$, such that $\sigma\ast G=1_{0}$, there exist
an $\cR$-homomorphism
\begin{align*}
&\hspace{1cm}G^{'}:\cD\longrightarrow Root\sigma\\
&\hspace{1.7cm}D\mapsto G(D),\\
&d:D_{1}\rightarrow D_{2}\mapsto d:D_{1}\rightarrow D_{2}
\end{align*}
and a 2-morphism $\alpha:G\Rightarrow R\circ G^{'}$ given by
$\alpha_{D}\triangleq 1_{G(D)}:G(D)\rightarrow (R\circ
G^{'}(D))=R(G(D))=G(D)$.

For any $\cC\in obj(\cR$-2-Mod), there is a functor
\begin{align*}
&\Psi:Hom(\cC,Root\sigma)\longrightarrow Hom(\cC,\cB)\\
&\hspace{2.8cm}H\mapsto R\circ H,\\
&\hspace{1cm}\tau:H_{1}\Rightarrow H_{2}\mapsto R\ast\tau
\end{align*}
where $R\ast\tau=1_{R}\ast\tau$ is the horizontal composition of
2-morphisms. For all $H_{1},H_{2}:\cC\rightarrow Root\sigma$, and
$\beta:R\circ H_{1}\Rightarrow R\circ H_{2}$, there exists a unique
$\chi:H_{1}\Rightarrow H_{2}$ given by $\chi_{C}\triangleq
\beta_{C}:H_{1}C\rightarrow H_{2}C$, so $R$ is fully faithful.
\end{itemize}

\noindent\textbf{Step 2.} $Root\sigma$ and $Im^{2}_{pl}F$ are
equivalent in ($\cR$-2-Mod).

Recall
\begin{align*}
&\hspace{3cm}\widehat{M}_{F}:Im^{2}_{pl}F\longrightarrow \cB\\
&\hspace{3.5cm}(A,\varphi,B)\mapsto B,\\
&g:(A_{1},\varphi_{1},B_{1})\rightarrow(A_{2},\varphi_{2},B_{2})\mapsto
g:B_{1}\rightarrow B_{2}.
\end{align*}
We have $(\sigma\ast\widehat{M}_{F})_{(A,\varphi,B)}=\sigma_{B}=B,\
\forall\ (A,\varphi,B)\in obj(Im^{2}_{pl}F)$, so
$\sigma\ast\widehat{M}_{F}=1_{0}$, and under the definition of root,
there is an $\cR$-homomorphism
\begin{align*}
&\hspace{3.2cm}\Phi:Im^{2}_{pl}F\longrightarrow Root\sigma\\
&\hspace{3.5cm}(A,\varphi,B)\mapsto B,\\
&g:(A_{1},\varphi_{1},B_{1})\rightarrow(A_{2},\varphi_{2},B_{2})\mapsto
g:B_{1}\rightarrow B_{2}.
\end{align*}
$\Phi$ is well-defined, i.e. $\Phi(A,\varphi,B)=B$ is an object of
$Root\sigma$, since there exist $A\in obj(\cA)$, and
$\widetilde{\varphi}:B\rightarrow FA+0$ under the composition
$$
B\xrightarrow[]{\varphi^{-1}}FA\xrightarrow[]{r_{FA}^{-1}}FA+0,
$$
i.e. $\widetilde{\varphi}=r_{FA}^{-1}\circ
\varphi^{-1}=(\varphi\circ r_{FA})^{-1}$.

There is an $\cR$-homomorphism
\begin{align*}
&\Phi^{-1}:Root\sigma\longrightarrow Im^{2}_{pl}F\\
&\hspace{1.7cm}B\mapsto (A,\varphi,B),\\
&g:B_{1}\rightarrow B_{2}\mapsto g:B_{1}\rightarrow B_{2}
\end{align*}
where $(A,\varphi,B)\in obj(Im^{2}_{pl}F)$, given by the
following way:\\
For $B\in obj(Root\sigma)$, there exist $A^{'}\in obj(\cA)$ and
$g^{'}:B\rightarrow FA^{'}+0$, so $A$ is the just $A^{'}$ of $\cA$,
and $\varphi$ is the composition
$FA\xrightarrow[]{r_{FA}^{-1}}FA+0\xrightarrow[]{(g^{'})^{-1}}B$.
Since, $\Phi^{-1}$ restricting on the morphisms is identity, then
$\Phi$ is an $\cR$-homomorphism.

For any $B\in obj(Root\sigma),\ g\in Mor(Root\sigma)$, we have
\begin{align*}
&(\Phi\circ\Phi^{-1})(B)=\Phi(A,\varphi,B)=B,\\
&(\Phi\circ\Phi^{-1})(g)=g.
\end{align*}
Also, for any $(A,\varphi,B)\in obj(Im^{2}_{pl}F)$, and $g\in
Mor(Im^{2}_{pl}F)$, we have
\begin{align*}
&(\Phi^{-1}\circ\Phi)(A,\varphi,B)=\Phi^{-1}(B)=B,\\
&(\Phi^{-1}\circ\Phi)(g)=g.
\end{align*}
Then $\Phi\circ\Phi^{-1}=1_{Root\sigma},\
\Phi^{-1}\circ\Phi=1_{Im^{2}_{pl}F}$.
\end{proof}

\begin{Prop}
Every 1-morphism $F:\cA\rightarrow\cB$ in ($\cR$-2-Mod) factors as
the following composite, where $E_{F}$ is full and surjective,
$\Omega_{F}$ is an equivalence, and $M_{F}$ is faithful.
$$
\cA\xrightarrow[]{E_{F}}Im^{1}F\xrightarrow[]{\Omega_{F}}Im^{2}F
\xrightarrow[]{M_{F}}\cB.
$$
\end{Prop}
\begin{proof}
\noindent\textbf{Step 1.} The $\cR$-2-module $Im^{1}F$ is described
in the following way:
\begin{itemize}
\item  Category $Im^{1}F$ consists of:

$\cdot$\ Objects are those of $\cA$.

$\cdot$\ Morphisms are the equivalent classes of morphisms in $\cA$,
for two morphisms $f_{1},f_{2}:A_{1}\rightarrow A_{2}$ are equal in
$Im^{1}F$, if $Ff_{1}=Ff_{2}$ in $\cB$, denote by $[f_1]$. The
composition and identities are those of $\cA$ up to equivalence.

\item $Im^{1}F$ is a symmetric 2-group.

The unit object is just the unit object 0 of $\cA$.

There is a bifunctor
\begin{align*}
&\hspace{0.4cm}+:Im^{1}F\times Im^{1}F\longrightarrow Im^{1}F\\
&\hspace{2.2cm}(A_{1},A_{2})\mapsto A_{1}+A_{2},\\
&(A_{1}\xrightarrow[]{[f_{1}]}
A_{1}^{'},A_{2}\xrightarrow[]{[f_{2}]} A_{2}^{'})\mapsto
A_{1}+A_{2}\xrightarrow[]{[f_{1}+f_{2}]} A_{1}^{'}+A_{2}^{'}
\end{align*}
where $A_{1}+A_{2},\ [f_{1}+f_{2}]$ are given under the monoidal
structure of $\cA$. If $f_{1},f_{1}^{'}:A_{1}\rightarrow A_{2}$ are
equal in $Im^{1}F$, i.e. $Ff_{1}=Ff_{1}^{'}$. From the following
commutative diagram
\begin{center}
\scalebox{0.9}[0.85]{\includegraphics{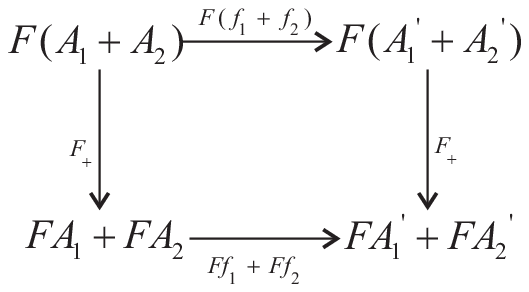}}
\end{center}
we have $F([f_{1}+f_{2}])=F([f_{1}^{'}+f_{2}])$, then $f_{1}+f_{2},\
f_{1}^{'}+f_{2}$ are equal in $Im^{1}F$.

Moreover, there are natural isomorphisms:
\begin{align*}
&<A_{1},A_{2},A_{3}>:A_{1}+A_{2}+A_{3}\rightarrow
A_{1}+A_{2}+A_{3},\\
&\hspace{2.2cm}l_{A}:0+A\rightarrow A,\\
&\hspace{2.2cm}r_{A}:A+0\rightarrow A,\\
&\hspace{1.6cm}c_{A_{1},A_{2}}:A_{1}+A_{2}\rightarrow A_{2}+A_{1},
\end{align*}
given by the natural isomorphisms in $\cA$. Since $\cA$ is a
symmetric monoidal category, so is $Im^{1}F$.

Since any morphism $[f]:A\rightarrow A^{'}$ in $Im^{'}F$ is in fact
a morphism in $\cA$ up to equivalence, and $\cA$ is a groupoid,
there exists $f^{*}:A^{'}\rightarrow A$ in $\cA$, such that
$[f^{*}]\circ [f]=[(f^*\circ f)]=1$ in $Im^{'}F$.

For any object $A\in obj(Im^{1}F)$, $A\in\cA$, there exist $A^{*}\in
obj(\cA)$, and $\eta_{A}:A^{*}+A\rightarrow 0$, so there are
$A^{*}\in obj(Im^{1}F)$, and $\eta_{A}:A^{*}+A\rightarrow 0$.

\item $Im^{1}F$ is an $\cR$-2-module.

The $\cR$-2-module structure is induced from the $\cR$-2-module
structure of $\cA$.
\end{itemize}

\noindent\textbf{Step 2.} The $\cR$-2-module $Im^{2}F$ is given by
following data:
\begin{itemize}
\item The category $Im^{2}F$ consists of:

$\cdot$\ Objects are the triple $(A,\varphi,B)$, where $A\in
obj(\cA),\ B\in obj(\cB),\ \varphi:FA\rightarrow B$ in $\cB$.

$\cdot$\ Morphism from $(A_{1},\varphi_{1},B_{1})$ to
$(A_{2},\varphi_{2},B_{2})$ is the equivalent class of a pair
$(f,g)$, denote by $[f,g]$, where $f:A_{1}\rightarrow A_{2},\
g:B_{1}\rightarrow B_{2}$, such that $g\circ
\varphi_{1}=\varphi_{2}\circ Ff,$ and for two morphisms
$(f,g),(f^{'},g^{'}):(A_{1},\varphi_{1},B_{1})\rightarrow
(A_{2},\varphi_{2},B_{2})$ are equal if $g=g^{'}$, or $Ff=Ff^{'}$.

$\cdot$\  Composition of morphisms
$(A_{1},\varphi_{1},B_{1})\xrightarrow[]{[f_{1},g_{1}]}(A_{2},\varphi_{2},B_{2})\xrightarrow[]{[f_{2},g_{2}]}$
 is given by $[f,g]=[f_{2},g_{2}]\circ[f_{1},g_{1}]$, where $f=f_{2}\circ f_{1},\ g=g_{2}\circ
g_{1}$, such that
$$
\varphi_{3}\circ F(f_{2}\circ f_{1})=\varphi_{3}\circ F(f_{2})\circ
F(f_{1})=g_{2}\circ\varphi_{2}\circ F(f_{1})=g_{2}\circ
g_{1}\circ\varphi_{1}=(g_{2}\circ g_{1})\circ\varphi_{1}.
$$
Also the above composition is well-defined, if $(f_{1},g_{1}),\
(f_{1}^{'},g_{1}^{'})$ are equal, then
$(f_{2},g_{2})\circ(f_{1},g_{1}),\
(f_{2},g_{2})\circ(f_{1}^{'},g_{1}^{'})$ are equal, because of
$g_{2}\circ g_{1}=g_{2}\circ g_{1}^{'}$.

\item $Im^{2}F$ is a symmetric 2-group.

The unit object is $(0,F_{0},0)$, where the first $0$ is the unit
object of $\cA$, the last $0$ is the unit object of $\cB$ and
$F_{0}:F0\rightarrow 0$.

There is a bifunctor
\begin{align*}
&\hspace{2cm}+:Im^{2}F\times Im^{2}F\longrightarrow Im^{2}F\\
&\hspace{0.7cm}((A_{1},\varphi_{1},B_{1}),(A_{2},\varphi_{2},B_{2}))\mapsto
(A_{1},\varphi_{1},B_{1})+(A_{2},\varphi_{2},B_{2})=(A,\varphi,B),\\
&(A_{1},\varphi_{1},B_{1})\xrightarrow[]{[f,g]}(A_{2},\varphi_{2},B_{2}),
(A_{1}^{'},\varphi_{1}^{'},B_{1}^{'})\xrightarrow[]{[f^{'},g^{'}]}(A_{2}^{'},\varphi_{2}^{'},B_{2}^{'})\mapsto
[f+f^{'},g+g^{'}]
\end{align*}
where $\varphi$ is the composition
$F(A_{1}+A_{2})\xrightarrow[]{F_{+}}FA_{1}+FA_{2}\xrightarrow[]{\varphi_{1}+\varphi_{2}}B_{1}+B_{2}$,
$[f+f^{'},g+g^{'}]$ makes the following diagram commutes:
\begin{center}
\scalebox{0.9}[0.85]{\includegraphics{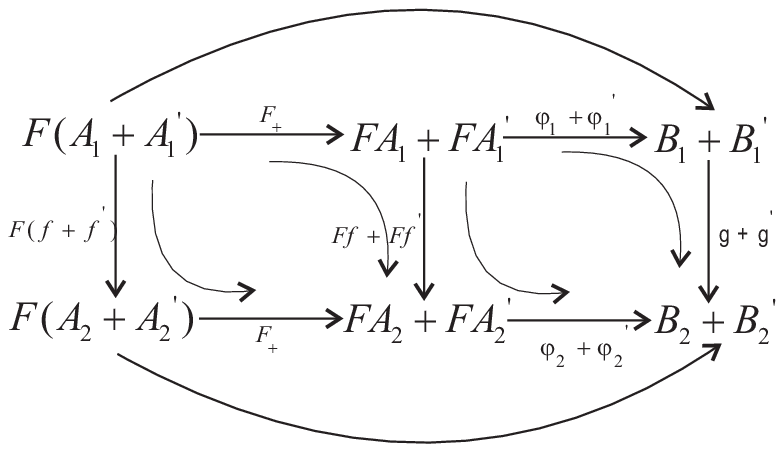}}
\end{center}
Similar as the above steps, the definition of addition of morphisms
is well-defined.

Moreover, there are natural isomorphisms:
\begin{align*}
&<(A_{1},\varphi_{1},B_{1}),(A_{2},\varphi_{2},B_{2}),(A_{3},\varphi_{3},B_{3})>\triangleq
(<A_{1},A_{2},A_{3}>,<B_{1},B_{2},B_{3}>),\\
&l_{(A,\varphi,B)}\triangleq (l_{A},l_{B}),\\
&r_{(A,\varphi,B)}\triangleq (r_{A},r_{B}),\\
&c_{(A_{1},\varphi_{1},B_{1}),(A_{2},\varphi_{2},B_{2})}\triangleq
(c_{A_{1},A_{2}},c_{B_{1},B_{2}}).
\end{align*}
Using the usual methods, we can check the above natural isomorphisms
are well-defined, and satisfy the conditions of symmetric 2-group.
\item $Im^{2}F$ is an $\cR$-2-module.

There is a bifunctor
\begin{align*}
&\hspace{4.2cm}\cdot:\cR\times Im^{2}F\longrightarrow Im^{2}F\\
&\hspace{4.4cm}(r,(A,\varphi,B))\mapsto
r\cdot(A,\varphi,B)\triangleq(r\cdot
A,\widetilde{r\cdot\varphi},r\cdot B),\\
&(r_{1}\xrightarrow[]{\alpha}r_{2},(A_{1},\varphi_{1},B_{1})\xrightarrow[]{[f,g]}(A_{2},\varphi_{2},B_2{})\mapsto
r\cdot(A_{1},\varphi_{1},B_{1})\xrightarrow[]{r\cdot[f,g]}r\cdot(A_{2},\varphi_{2},B_2{})
\end{align*}
where $\widetilde{r\cdot\varphi}:F(r\cdot
A)\xrightarrow[]{F_{2}}r\cdot FA\xrightarrow[]{r\cdot\varphi}r\cdot
B$, from $r\cdot$ is a functor and $F$ is an $\cR$-homomorphism, we
have the following commutative diagram:
\begin{center}
\scalebox{0.9}[0.85]{\includegraphics{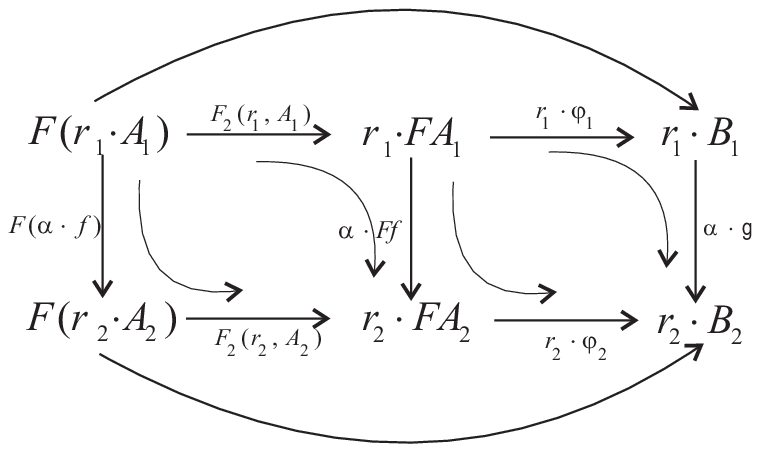}}
\end{center}
Also, if $(f,g),(f^{'},g^{'})$ are equal, then $\varphi\cdot(f,g),\
\varphi\cdot(f^{'},g^{')}$ are equal too.

Moreover, there are natural isomorphisms:
\begin{align*}
&a_{(A_{1},\varphi_{1},B_{1}),(A_{2},\varphi_{2},B_2{})}^{r}\triangleq
(a_{A_{1},A_{2}}^{r},a_{B_{1},B_{2}}^{r}),\\
&\hspace{2cm}b_{r_{1},r_{2}}^{(A,\varphi,B)}\triangleq(b_{r_{1},r_{2}}^{A},b_{r_{1},r_{2}}^{B}),\\
&\hspace{1.3cm}b_{r_{1},r_{2},(A,\varphi,B)}\triangleq(b_{r_{1},r_{2},A},b_{r_{1},r_{2},B}),\\
&\hspace{2.1cm}i_{(A,\varphi,B)}\triangleq(i_{A},i_{B}),\\
&\hspace{3cm}z_{r}\triangleq(z_{r},z_{r}).
\end{align*}
By the usual methods, we can check the above natural isomorphisms
are well-defined, and Fig.18.--31. commute.
\end{itemize}
\noindent\textbf{Step3.} $Im^{1}F,\ Im^{1}F$ are equivalent as
$\cR$-2-modules.

Define a functor
\begin{align*}
&\Omega_{F}:Im^{1}F\longrightarrow Im^{2}F\\
&\hspace{1.7cm}A\mapsto (A,1_{FA},FA),\\
&[f]:A_{1}\rightarrow A_{2}\mapsto
[f,Ff]:(A_{1},1_{FA_{1}},FA_{1})\rightarrow(A_{2},1_{FA_{2}},FA_{2})
\end{align*}
the above definition is well-defined. In fact, if
$f,f^{'}:A_{1}\rightarrow A_{2}$ are equal in $Im^{1}F$, i.e.
$Ff=Ff^{'}$, then we have $(f,Ff),\ (f^{'},Ff^{'})$ are equal in
$Im^{2}F$, i.e. $\Omega_{F}(f),\ \Omega_{F}(f^{'})$ are equal.

For any identity morphism $1_{A}:A\rightarrow A$ in $Im^{1}F,\
\Omega_{F}(1_{A})=(1_{A},F(1_{A})=(1_{A},1_{FA}):(A,1_{FA},FA)\rightarrow
(A,1_{FA},FA)$ is the identity morphism in $Im^{2}F$.

For any morphisms
$A_{1}\xrightarrow[]{[f_{1}]}A_{2}\xrightarrow[]{[f_{2}]}A_{3}$ in
$Im^{1}F$, we have $ \Omega_{F}([f_{2}\circ
f_{1}])=\Omega([f_{2}\circ f_{1}])=[f_{2}\circ f_{1},F(f_{2}\circ
f_{1})]=[f_{2}\circ f_{1},Ff_{2}\circ Ff_{1}]=[f_{2},Ff_{2}]\circ
[f_{1},Ff_{1}]=\Omega_{F}([f_{2}])\circ\Omega_{F}([f_{1}]).$ Then
$\Omega_{F}$ is a functor.

There are natural isomorphisms:
\begin{align*}
&(\Omega_{F})_{+}\triangleq(1_{A_{1}+A_{2}},F_{+}):\Omega_{F}(A_{1}+A_{2})=(A_{}+A_{},1_{F(A_{1}+A_{2})},F(A_{1}+A_{2}))
\rightarrow\Omega_{F}(A_{1})+\Omega_{F}(A_{2})\\
&\hspace{1cm}=(A_{1},1_{FA_{1}},FA_{1})+(A_{2},1_{FA_{2}},FA_{2})=(A_{}+A_{},(1_{FA_{1}}+1_{FA_{2}})\circ
F_{+},FA_{1}+FA_{2}),\\
&(\Omega_{F})_{0}\triangleq(1_{0},F_{0}):\Omega_{F}(0)=(0,1_{F0},F0)\rightarrow
(0,F_{0},0),\\
&(\Omega_{F})_{2}\triangleq (1_{r\cdot A},F_{2}):\Omega_{F}(r\cdot
A)=(r\cdot A,1_{F(r\cdot A)},F(r\cdot A)\rightarrow
r\cdot\Omega_{F}(A)=r\cdot(A,1_{FA},FA)\\
&\hspace{1cm}=(r\cdot A,r\cdot1_{FA}\circ F_{2},r\cdot FA)=((r\cdot
A,F_{2},r\cdot FA).
\end{align*}
After basic calculations, the above natural isomorphisms are
well-defined, and
$((\Omega_{F})_{+},(\Omega_{F})_{0},(\Omega_{F})_{2})$ is an
$\cR$-homomorphism, i.e. Fig.32--36.commute.

Define a functor
\begin{align*}
&\hspace{2.9cm}\Omega_{F}^{-1}:Im^{2}F\longrightarrow Im^{1}F\\
&\hspace{3.5cm}(A,\varphi,B)\mapsto A,\\
&(A_{1},\varphi_{1},B_{1})\xrightarrow[]{[f,g]}(A_{2},\varphi_{2},B_{2})\mapsto
A_{1}\xrightarrow[]{[f]} A_{2}.
\end{align*}
If $(f,g),\ (f^{'},g^{'})$ are equal in $Im^{2}F$, i.e. $Ff=Ff^{'}$,
then $f,\ f^{'}$ are equal in $Im^{1}F$.

From the definition of $\Omega_{F}^{-1}$, we see that
$\Omega_{F}^{-1}$ maps morphisms of $Im^{2}F$ to the first part of
them, so $\Omega_{F}^{-1}$ is an $\cR$-homomorphism with
$(\Omega_{F}^{-1})_{+}=id,\ (\Omega_{F}^{-1})_{0}=id,\
(\Omega_{F}^{-1})_{2}=id$.

For any $A\in obj(Im^{1}F)$, and any representative morphism
$f:A_{1}\rightarrow A_{2}$ in $Im^{1}f$, we have
\begin{align*}
&(\Omega_{F}^{-1}\circ\Omega_{F})(A)=\Omega_{F}^{-1}(A,1_{FA},FA)=A,\\
&(\Omega_{F}^{-1}\circ\Omega_{F})(f)=\Omega_{F}^{-1}(f,Ff)=f.
\end{align*}
So $\Omega_{F}^{-1}\circ\Omega_{F}=1_{Im^{1}F}$.

There is a morphism of $\cR$-homomorphisms:
\begin{align*}
&\tau:\Omega_{F}\circ\Omega_{F}^{-1}\Rightarrow
1:Im^{2}F\longrightarrow Im^{2}F,\\
&\tau_{(A,\varphi,B)}\triangleq
(1_{A},\varphi):(\Omega_{F}\circ\Omega_{F}^{-1})(A,\varphi,B)=\Omega_{F}(A)=(A,1_{FA},FA)\rightarrow
(A,\varphi,B).
\end{align*}
For any morphism $(f,g):(A_{1},\varphi_{1},B_{1})\rightarrow
(A_{2},\varphi_{2},B_{2})$ in $Im^{2}F$, with
$g\circ\varphi_{1}=\varphi_{2}\circ Ff$,
$(\Omega_{F}\circ\Omega_{F}^{-1})(f,g)=\Omega_{F}(f)=(f,Ff)$, we
have the following commutative diagram:
\begin{center}
\scalebox{0.9}[0.85]{\includegraphics{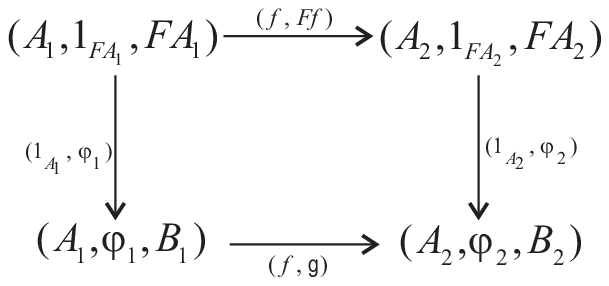}}
\end{center}
So $\tau$ is a natural transformation.
\begin{align*}
&(\Omega_{F}\circ\Omega_{F}^{-1})((A_{1},\varphi_{1},B_{1})+(A_{2},\varphi_{2},B_{2}))
=(A_{1}+A_{2},1_{F(A_{1}+A_{2})},F(A_{1}+A_{2}),\\
&(\Omega_{F}\circ\Omega_{F}^{-1})(A_{1},\varphi_{1},B_{1})
+(\Omega_{F}\circ\Omega_{F}^{-1})(A_{2},\varphi_{2},B_{2})
=(A_{1}+A_{2},F_{+},FA_{1}+FA_{2}),\\
&\tau_{(A_{1},\varphi_{1},B_{1})+(A_{2},\varphi_{2},B_{2})}=(1_{A_{1}+A_{2}},(\varphi_{1}+\varphi_{2})\circ
F_{+}),\\
&\tau_{(A_{i},\varphi_{i},B_{i})}=(1_{A_{i}},\varphi_{i}),\
\mbox{for $i=1,2$.}
\end{align*}
Thus Fig.7 commutes. We can also get commutative diagrams Fig.37. in
the similar way. Then $\tau$ is a morphisms of $\cR$-homomorphisms.

Since $\tau_{(A,\varphi,B)}=(1_{A},\varphi)$ and $Im^{2}F$ is a
groupiod, so $\tau$ is an isomorphism.

\noindent\textbf{Step 4.} There is a full and surjective
$\cR$-homomorphism $E_{F}$.

Define a functor
\begin{align*}
&\hspace{0.8cm}E_{F}:\cA\longrightarrow Im^{1}F\\
&\hspace{1.7cm}A\mapsto A\\
&f:A_{1}\rightarrow A_{2}\mapsto [f]:A_{1}\rightarrow A_{2}
\end{align*}
Obviously, $E_{F}$ is an $\cR$-homomorphism.

\noindent\textbf{Step 5.} There is a faithful $\cR$-homomorphism
$M_{F}$.

Define a functor
\begin{align*}
&\hspace{3cm}M_{F}:Im^{2}F\longrightarrow \cB\\
&\hspace{3.5cm}(A,\varphi,B)\mapsto B,\\
&(A_{1},\varphi_{1},B_{1})\xrightarrow[]{[f,g]}(A_{2},\varphi_{2},B_{2})\mapsto
B_{1}\xrightarrow[]{g} B_{2}.
\end{align*}
Obviously, $M_{F}$ is a faithful $\cR$-homomorphism.

\noindent\textbf{Step 6.} $F=M_{F}\circ\Omega_{F}\circ E_{F}.$

For any $A\in obj(\cA)$ and any morphism $f\in Mor(\cA)$, we have
$$
(M_{F}\circ\Omega_{F}\circ
E_{F})(A)=(M_{F}\circ\Omega_{F})(A)=M_{F}(A,1_{FA},FA)=FA,
$$
$$
(M_{F}\circ\Omega_{F}\circ
E_{F})(f)=(M_{F}\circ\Omega_{F})(f)=M_{F}(f,Ff)=Ff.
$$
\end{proof}

\begin{Prop}
For each $\cR$-homomorphism $F:\cA\rightarrow\cB$ in ($\cR$-2-Mod),
$M_{F}:Im^{2}F\rightarrow\cB$ is the kernel of the cokernel of F.
\end{Prop}
\begin{proof}
The cokernel of F is $(CokerF,\ p_{F},\ \pi_{F})$, where $CokerF$ is
a category whose object is the object of $\cB$, morphism is
$(f,A):B_{1}\rightarrow B_{2}$, with $A\in obj(\cA),\
f:B_{1}\rightarrow B_{2}+FA$, $p_{F}(B)=B,\
(\pi_{F})_{A}=[(\pi_{F})_{A}^{'},A]$, where
$(\pi_{F})_{A}^{'}:FA\rightarrow FA+0\rightarrow FA+F0$.
\begin{itemize}
\item Let us describe the $Kerp_{F}$ in the following way:

$\cdot$\ Objects are the triple $(B,A,b)$, where $B\in obj(\cB),\
A\in obj(\cA),\ b:B\rightarrow FA+0$.

$\cdot$\ Morphism from $(B,A,b)$ to $(B^{'},A^{'},b^{'})$ is the
equivalence class of a pair $(g,f)$, denote by $[g,f]$, where
$g:B\rightarrow B^{'}$, $f:A\rightarrow A^{'}$, such that $
(Ff+1_{0})\circ b=b^{'}\circ g$, and for two morphisms
$(g,f),(g^{'},f^{'}):(B,A,b)\rightarrow (B^{'},A^{'},b^{'})$ are
equal, if $g=g^{'}$.

$\cdot$\ Composition of morphisms. Given morphisms
$(B_{1},A_{1},b_{1})\xrightarrow[]{[g_{1},f_{1}]}(B_{2},A_{2},b_{2})\\
\xrightarrow[]{[g_{2},f_{2}]}(B_{3},A_{3},b_{3})$,
$[(g_{2},f_{2})\circ(g_{1},f_{1})]\triangleq [g_{2}\circ
g_{1},f_{2}\circ f_{1}]$. If $(g_{1},f_{1}),\
(g_{1}^{'},f_{1}^{'}):(B_{1},A_{1},b_{1})\rightarrow(B_{1},A_{1},b_{1})$
are equal, i.e. $g_{1}=g_{1}^{'}$, then
$(g_{2},f_{2})\circ(g_{1},f_{1}),\
(g_{2},f_{2})\circ(g_{1}^{'},f_{1}^{'})$ are equal.

There is a bifunctor
\begin{align*}
&\hspace{2cm}+:Kerp_{F}\times Kerp_{F}\longrightarrow Kerp_{F}\\
&\hspace{1.2cm}((B_{1},A_{1},b_{1}),(B_{2},A_{2},b_{2}))\mapsto
(B_{1},A_{1},b_{1})+(B_{2},A_{2},b_{2})\triangleq
(B_{1}+B_{2},A_{1}+A_{2},b),\\
&\hspace{2.7cm}((g_{1},f_{1}),(g_{2},f_{2}))\mapsto (g,f)\triangleq
(g_{1}+g_{2},f_{1}+f_{2})
\end{align*}
Using the similar methods in Theorem 2, $Kerp_{F}$ is an
$\cR$-2-module.
\item
There is a 2-morphism $\epsilon:p_{F}\circ M_{F}\Rightarrow 0$,
given by $\epsilon_{(A,\varphi,B)}\triangleq
[l_{FA}\circ\varphi^{-1},A]:(p_{F}\circ
M_{F})(A,\varphi,B)=p_{F}(B)=B\rightarrow 0$. By the universal
property of the kernel, there is an $\cR$-homomorphism

\begin{align*}
&\hspace{3.5cm}\Theta:Im^{2}F\longrightarrow Kerp_{F}\\
&\hspace{3.6cm}(A,\varphi,B)\mapsto (B,A,b),\\
&[f,g]:(A,\varphi,B)\rightarrow (A^{'},\varphi^{'},B^{'})\mapsto
[g,f]:(B,A,b)\rightarrow(B^{'},A^{'},b^{'})
\end{align*}
where $b,\ b^{'}$ are the compositions
$B\xrightarrow[]{\varphi^{-1}}FA\xrightarrow[]{r_{FA}}FA+0,\
B^{'}\xrightarrow[]{(\varphi^{'})^{-1}}FA^{'}\xrightarrow[]{r_{FA^{'}}}FA^{'}+0$,
respectively, such that
$$
(Ff+1_{0})\circ b=(Ff+1_0)\circ r_{FA}\circ
\varphi_{-1}=r_{FA^{'}}\circ Ff\varphi_{-1}=r_{FA^{'}}\circ
(\varphi^{'})^{-1}\circ g=b^{'}\circ g.
$$
If $(f,g),(f^{'},g^{'}):(A,\varphi,B)\rightarrow
(A^{'},\varphi^{'},B^{'})$ are equal, i.e. $g=g^{'}$, then
$\Theta(f,g),\Theta(f^{'},g^{'})$ are equal in $Kerp_{F}$.

\item There is an $\cR$-homomorphism
\begin{align*}
&\hspace{2.7cm}\Theta^{-1}:Kerp_{F}\longrightarrow Im^{2}F\\
&\hspace{3.5cm}(B,A,b)\mapsto(A,B,\varphi),\\
&[g,f]:(B,A,b)\rightarrow(B^{'},A^{'},b^{'})\mapsto[f,g]:(A,B,\varphi)\rightarrow(A^{'},B^{'},\varphi^{'})
\end{align*}
where $\varphi=b^{-1}\circ r_{FA}^{-1}:FA\rightarrow FA+0\rightarrow
B$, and $(f,g)$ satisfy
$$
g\circ\varphi=g\circ b^{-1}\circ
r_{FA}^{-1}=(b^{'})^{-1}\circ(Ff+1_0)\circ
r_{FA}^{-1}=(b^{'})^{-1}\circ r_{FA^{'}}^{-1}\circ
Ff=\varphi^{'}\circ Ff.
$$
Similarly, if $(g,f),\ (g^{'},f^{'})$ are equal in $Im^{2}F$, then
$\Theta^{-1}(g,f),\ \Theta^{-1}(g^{'},f^{'})$ are equal in
$Kerp_{F}$, too. Moreover, $\Theta^{-1}$ is an $\cR$-homomorphism.

For any $(A,\varphi,B)\in obj(Im^{2}F),\ (f,g)\in Mor(Im^{2}F)$, we
have
$$
(\Theta^{-1}\circ\Theta)(A,\varphi,B)=\Theta^{-1}(B,A,r_{FA}^{-1}\circ\varphi^{-1})=(A,(r_{FA}^{-1}\circ\varphi^{-1})^{-1}\circ
r_{FA}^{-1},B)=(A,\varphi,B),
$$
$$
(\Theta^{-1}\circ\Theta)(f,g)=\Theta^{-1}(g,f)=(f,g).
$$
For any $(B,A,b)\in obj(Kerp_{F}),
 (g,f)\in Mor(Kerp_{F})$, we have
$$
(\Theta\circ\Theta^{-1})(B,A,b)=\Theta(A,B,b^{-1}\circ
r_{FA}^{-1})=(B,A,r_{FA}^{-1}\circ(b^{-1}\circ
r_{FA}^{-1}))=(B,A,b),
$$
$$
(\Theta\circ\Theta^{-1})(g,f)=\Theta(f,g)=(g,f).
$$
Thus $\Theta:Im^{2}F\rightarrow Kerp_{F}$ is an equivalent.
\end{itemize}
\end{proof}

\begin{Prop}
For each $\cR$-homomorphism $F:\cA\rightarrow\cB$ in ($\cR$-2-Mod),
$E_{F}:\cA\rightarrow Im^{1}F$ is the coroot of the pip of F.
\end{Prop}
\begin{proof}
The pip of F is an $\cR$-2-module $PipF$, together with 2-morphism
$\pi_{F}:0\Rightarrow 0:PipF\rightarrow\cA$, whose component at $a$
is $a$ itself, where $a:0\rightarrow 0$ is an object of $PipF$.

\begin{itemize}
\item Let us describe $Coroot\pi_{F}$ in the following way:

$\cdot$\ Objects are the objects of $\cA$.

$\cdot$\ Morphisms are the morphisms of $\cA$ up to equivalent
relations, for two morphisms $f,f^{'}:A_{1}\rightarrow A_{2}$ which
are equal in $Coroot\pi_{F}$, if there exists $a:0\rightarrow0$ in
$\cA$, such that  $Fa=1_{F0}$ and
$$
f\circ r_{A_{1}}=r_{A_{2}}\circ (f^{'}+a).
$$
$\cdot$\ Composition of morphisms and identity morphisms are the
composition and identity in $\cA$ up to equivalence, and also
well-defined. In fact, given morphism $f_{2}:A_{2}\rightarrow A_{3}$
and equal morphisms $f_{1},f_{1}^{'}:A_{1}\rightarrow A_{2}$ in
$Coroot\pi_{F}$, for equivalence of $f_{1},\ f_{1}^{'}$, there
exists $a:0\rightarrow 0$, such that $Fa=1_{F0},\
f_{1}=r_{A_{1}}\circ (f_{1}^{'}+a)\circ r_{A_{1}}^{-1}$. Then there
exists $a:0\rightarrow 0$, such that
$$
f_{2}\circ f_{1}\circ r_{A_{1}}=f_{2}\circ r_{A_{2}}\circ
(f_{1}^{'}+a)=r_{A_{3}}\circ (f_{2}^{'}+1_0)\circ
(f_{1}^{'}+a)=r_{A_{3}}\circ (f_{2}^{'}\circ f_{1}^{'}+a).
$$
Then $f_{2}\circ f_{1},\ f_{2}\circ f_{1}^{'}$ are equal in
$Coroot\pi_{F}$.

$Coroot\pi_{F}$ is an $\cR$-2-module from the $\cR$-2-module
structure of $\cA$ up to equivalence.

There is an $\cR$-homomorphism
\begin{align*}
&\hspace{1cm}R:\cA\longrightarrow Coroot\pi_{F}\\
&\hspace{1.7cm}A\mapsto A,\\
&f:A_{1}\rightarrow A_{2}\mapsto f:A_{1}\rightarrow A_{2}
\end{align*}
such that for any $a:0\rightarrow 0$ in $PipF,\
(R\ast\pi_{F})_{a}=R(a)=a=(1_{0})_{a}$, i.e. $R\ast\pi_{F}=1_{0}$.

For $\cK\in obj(\cR$-2-Mod), and an $\cR$-homomorphism
$G:\cA\rightarrow \cK$, such that $G\ast \pi_{F}=1_{0}$, there exist
an $\cR$-homomorphism
\begin{align*}
&G^{'}:Coroot\pi_{F}\longrightarrow\cK\\
&\hspace{2.1cm}A\mapsto GA,\\
&\hspace{0.3cm}f:A_{1}\rightarrow A_{2}\mapsto
G(f):G(A_{1})\rightarrow G(A_{2})
\end{align*}
and a 2-morphism
\begin{align*}
&\alpha:G^{'} \circ R\Rightarrow G\\
&\hspace{1.4cm}A\mapsto \alpha_{A}\triangleq id_{GA}:(G^{}\circ
R)(A)=G^{'}(A)=GA\rightarrow GA.
\end{align*}
From the given $G$, we know that $G^{'}$ is an $\cR$-homomorphism,
$\alpha$ is a 2-morphism.

For every $\cC\in obj(\cR$-2-Mod), there is an $\cR$-homomorphism
\begin{align*}
&-\circ R:Hom(Coroot\pi_{F},\cC)\longrightarrow Hom(\cA,\cC)\\
&\hspace{4.1cm}H\mapsto H\circ R,\\
&\hspace{2.4cm}\tau:H_{1}\Rightarrow H_{2}\mapsto \tau\ast
R\triangleq \tau\ast1_{R}:H_{1}\circ R\Rightarrow H_{2}\circ R
\end{align*}
such that for any objects $H_{1},\ H_{2}$ in
$Hom(Coroot\pi_{F},\cC)$, and a 2-morphism $\beta:H_{1}\circ
R\Rightarrow H_{2}\circ R$, there is a 2-morphism
$\tau:H_{1}\Rightarrow H_{2}$ given by
$\tau_{A}\triangleq\beta_{A}:H_{1}A\rightarrow H_{2}A$. Also, if
$\tau_{1},\ \tau_{2}:H_{1}\Rightarrow H_{2}:
Coroot\pi_{F}\rightarrow\cC$, such that $\tau_{1}\ast R=\tau_{2}\ast
R:H_{1}\circ R\Rightarrow H_{2}\circ R$, then, for any $A\in
obj(\cA)$, $(\tau_{1})_{A}=(\tau_{1}\ast R)_{A}=(\tau_{2}\ast
R)_{A}=(\tau_{2})_{A}$, i.e. $\tau_{1}=\tau_{2}$. So $R$ is full and
surjective, and then $(Coroot\pi_{F},R)$ is the coroot of $\pi_{F}$.

\item There is an equivalence between $Coroot\pi_{F}$ and $Im^{1}F$.

From the definition of $E_{F}:\cA\rightarrow Im^{1}F$, we have
$E_{F}\ast\pi_{F}=1_{0}$, from the universal property of coroot of
$\pi_{F}$, there is an $\cR$-homomorphism
\begin{align*}
&\Theta:Coroot\pi_{F}\longrightarrow Im^{1}F\\
&\hspace{2cm}A\mapsto A\\
&\hspace{0.2cm}f:A_{1}\rightarrow A_{2}\mapsto f:A_{1}\rightarrow
A_{2}.
\end{align*}
Also there is an $\cR$-homomorphism
\begin{align*}
&\Theta^{-1}:Im^{1}F\longrightarrow Coroot\pi_{F}\\
&\hspace{2cm}A\mapsto A,\\
&\hspace{0.2cm}f:A_{1}\rightarrow A_{2}\mapsto f:A_{1}\rightarrow
A_{2}.
\end{align*}
In \cite {18}, the author proved that $\Theta$ and $\Theta^{-1}$ are
well-defined homomorphism of symmetric 2-groups. Since
$Coroot\pi_{f},\ Im^{1}F$ have the same $\cR$-2-module structure and
$\Theta\circ\Theta^{-1}=1,\Theta^{-1}\circ\Theta=1,$ then
$Coroot\pi_{f},\ Im^{1}F$ are equivalent $\cR$-2-modules.
\end{itemize}
\end{proof}
In the sense of 2-abelian $Gpd$-category, and from Propositions
6-11, we have
\begin{Thm}
The $Gpd$-category ($\cR$-2-Mod) is a 2-abelian $Gpd$-category.
\end{Thm}

\noindent Fang Huang, Shao-Han Chen, Wei Chen\\
Department of Mathematics\\
 South China University of
 Technology\\
 Guangzhou 510641, P. R. China

\noindent Zhu-Jun Zheng\\
Department of Mathematics\\
 South China University of
 Technology\\
 Guangzhou 510641, P. R. China \\
 and\\
Institute of Mathematics\\
Henan University\\  Kaifeng 475001, P. R.
China\\
E-mail: zhengzj@scut.edu.cn
\end{document}